\documentclass[letter,11pt]{article}

\usepackage[linktocpage=true]{hyperref}
\usepackage[margin=1.2 in, top=1 in, bottom= 1.2 in]{geometry}

\usepackage[normalem]{ulem}

\usepackage{tikz}
\usetikzlibrary{hobby}
\usepackage{tikz-cd}

\usepackage{lipsum}
\newcommand\blfootnote[1]{%
  \begingroup
  \renewcommand\thefootnote{}\footnote{#1}%
  \addtocounter{footnote}{-1}%
  \endgroup
}

\makeatletter
\g@addto@macro\normalsize{
  \setlength\abovedisplayskip{8pt}
  \setlength\belowdisplayskip{8pt}
  \setlength\abovedisplayshortskip{8pt}
  \setlength\belowdisplayshortskip{8pt}
  }
\makeatother

\usepackage{tocloft}

\usepackage[english]{babel}
\usepackage{amsfonts}
\usepackage{mathrsfs}
\usepackage{bbm}
\usepackage{latexsym}
\usepackage{math dots}
\usepackage{amssymb}
\usepackage{mathtools}

\usepackage{enumitem}
\setlist{nolistsep}
\usepackage{amsthm}
\usepackage[capitalize]{cleveref}
\crefformat{footnote}{#2\footnotemark[#1]#3} 

\newcommand\eqnitem[1][]{%
  \ifx\relax#1\relax  \item \else \item[#1] \fi
  \abovedisplayskip=0pt\abovedisplayshortskip=0pt~\vspace*{-\baselineskip}}

\newtheoremstyle{plain}{3mm}{3mm}{\slshape}{}{\bfseries}{.}{.5em}{}
\newtheoremstyle{definition}{2mm}{2mm}{}{}{\bfseries}{.}{.5em}{}
\theoremstyle{plain}
	
\newtheorem{theorem}{Theorem}

\newtheorem{lemma}[theorem]{Lemma}

\newtheorem{corollary}[theorem]{Corollary}

\newtheorem{question}[theorem]{Question}
\theoremstyle{definition}
\newtheorem{definition}[theorem]{Definition}
\newtheorem{remark}[theorem]{Remark}

\theoremstyle{plain}
\newcounter{MainTheoremCounter}

\newtheorem{Maintheorem}[MainTheoremCounter]{Theorem}

\theoremstyle{plain}
\newtheorem*{namedthm}{\namedthmname}
\newcounter{namedthm}
	

\usepackage{chngcntr}
\counterwithin{theorem}{section}

\numberwithin{equation}{section}

\allowdisplaybreaks

\usepackage{xcolor}
\definecolor{Color2}{rgb}{0.78, 0.11, 0.0}
\hypersetup{citecolor = black,colorlinks,
			linkcolor = black,
			urlcolor = Color2}

\usepackage{titlesec}
\titleformat{\section}
  {\Large\center\bfseries}
  {\thesection.}{.7em}{}
\titlespacing*{\section}{0pt}{3.5ex plus 0ex minus 0ex}{1.5ex plus 0ex}
\titleformat{\subsection}
  {\large\center\bfseries}
  {\thesubsection.}{.7em}{}
\titlespacing*{\subsection}{0pt}{3.5ex plus 0ex minus 0ex}{1.5ex plus 0ex}
\titleformat{\subsubsection}
  {\center\bfseries}
  {\thesubsubsection.}{.7em}{}
\titlespacing*{\subsubsection}{0pt}{3.5ex plus 0ex minus 0ex}{1.5ex plus 0ex}

\addto\captionsenglish{}

\usepackage{titling}
\newcommand{\Folner}{F\o{}lner}

\newcommand{\Szemeredi}{Szemer\'{e}di}

\newcommand{\eps}{\epsilon}
\newcommand{\N}{\mathbb{N}}
\newcommand{\Z}{\mathbb{Z}}
\newcommand{\R}{\mathbb{R}}

\newcommand{\Q}{\mathbb{Q}}
\newcommand{\T}{\mathbb{T}}
\newcommand{\Nz}{\N_{\geq 0}}
\newcommand{\defeq}{\vcentcolon=}

\newcommand\restr[2]{{ \left.\kern-\nulldelimiterspace #1 \right|_{#2}}}

\renewcommand{\epsilon}{\varepsilon}
\renewcommand{\leq}{\leqslant}
\renewcommand{\geq}{\geqslant}
\renewcommand{\setminus}{\backslash}

\usepackage[normalem]{ulem}

\newcommand{\spaceD}{W}
\newcommand{\Tmap}{\mathcal{T}}
\newcommand{\Smap}{\mathcal{S}}

\newcommand{\Gdot}{G^{\bullet}}
\newcommand{\Hdot}{H^{\bullet}}
\newcommand{\Kdot}{K^{\bullet}}
\newcommand{\Gammadot}{\Gamma^{\bullet}}
\newcommand{\poly}{{\normalfont\text{Poly}}\,}
\newcommand{\perpoly}[1]{\mathcal{P}_{#1}}

\newcommand{\ccd}{C}
\newcommand{\ddset}{Y}
\newcommand{\diag}[1]{\underline{#1}}
\newcommand{\mapspace}{\mathfrak{F}}
\newcommand{\uex}{\underline{e_X}}
\newcommand{\ncc}{\#_{\text{cc}}}
\newcommand{\cx}{\ncc(X)}

\newcommand{\compc}{\#}

\newcommand{\dcc}{d_{\text{cc}}}
\newcommand{\prolong}{D}
\newcommand{\lsc}{{\sc lsc}}
\newcommand{\usc}{{\sc usc}}
\newcommand{\INT}{\mathcal{I}}
\newcommand{\JINT}{\mathcal{J}}
\newcommand{\KINT}{\mathcal{K}}

\newcommand{\nn}{\N_{N,1}}
\newcommand{\zn}{\Z_{N,1}}
\newcommand{\fcx}{2^X}
\newcommand{\fcy}{2^Y}
\newcommand{\fcspaced}{2^{\spaceD}}
\newcommand{\orbit}[2]{\overline{#1^{\N} #2 }}
\newcommand{\orbitz}[2]{\overline{#1^{\Z} #2 }}
\newcommand{\longorbit}[2]{\overline{(#1)^{\N} #2 }}
\newcommand{\longorbitz}[2]{\overline{(#1)^{\Z} #2 }}

\newcommand{\one}{\mathbbm{1}}

\newcommand{\sa}{\mathrm{SA}}

\begin{document}

\title{\bfseries Simultaneous approximation in nilsystems and the multiplicative thickness of return-time sets}
\author{Daniel Glasscock}

\date{}
\maketitle

\blfootnote{2020 \emph{Mathematics Subject Classification.}  Primary: 37B20. Secondary: 37B05, 05D10.}
\blfootnote{\emph{Key words and phrases. } Topological multiple recurrence, simultaneous approximation, rational points and polynomials in nilpotent Lie groups and nilmanifolds, prolongation relation, multiplicatively thick sets, return-time sets, nil-Bohr sets, additively syndetic sets, van der Waerden's theorem.}

\begin{abstract}
In the topological dynamical system $(X,T)$, a point $x$ simultaneously approximates a point $y$ if there exists a sequence $n_1$, $n_2$, \dots of natural numbers for which $T^{n_i} x$, $T^{2n_i}x$, \dots, $T^{k n_i} x$ all tend to $y$. In 1978, Furstenberg and Weiss showed that every system possesses a point which simultaneously approximates itself (a multiply recurrent point) and deduced refinements of van der Waerden's theorem on arithmetic progressions.
In this paper, we study the denseness of the set of points that are simultaneously approximated by a given point.  We show that in a minimal nilsystem, all points simultaneously approximate a $\delta$-dense set of points under a necessarily restricted set of powers of $T$.  We tie this theorem to the multiplicative combinatorial properties of return-time sets, showing that all nil-Bohr sets and typical return-time sets in a minimal system are multiplicatively thick in a coset of a multiplicative subsemigroup of the natural numbers.  This yields an inhomogeneous multiple recurrence result that generalizes Furstenberg and Weiss' theorem and leads to new enhancements of van der Waerden's theorem.  This work relies crucially on continuity in the prolongation relation (the closure of the orbit-closure relation) developed by Auslander, Akin, and Glasner; the theory of rational points and polynomials on nilmanifolds developed by Leibman, Green, and Tao; and the machinery of topological characteristic factors developed recently by Glasner, Huang, Shao, Weiss, and Ye.
\end{abstract}


\clearpage

\small
\setcounter{tocdepth}{2}
\tableofcontents
\thispagestyle{empty}
\normalsize


\clearpage

\section{Introduction}
\label{sec_intro}

Let $X$ be a compact metric space and $\mapspace$ be a collection of continuous maps from $X$ to itself.  A point \emph{$x$ simultaneously approximates a point $y$ under $\mapspace$} if for all finite sub-collections $F \subseteq \mapspace$ and all $\eps > 0$, there exists $m \in \N \defeq \{1,2,\ldots\}$ such that for all $f \in F$, the distance between $f^m x$ and $y$ is less than $\eps$.  In this paper, we describe 1) the denseness of the set of points that are simultaneously approximated by a point under iterates of a nilrotation on a nilmanifold; 2) the multiplicative quality of the set of times of return of a point to a non-empty, open set under iterates of a single continuous map; and 3) combinatorial corollaries of both 1) and 2), including enhancements to van der Waerden's theorem on arithmetic progressions and new multiplicative configurations in translates of difference sets and their higher-order recurrence analogues.

\subsection{Simultaneous approximation in topological dynamical systems}

The phenomenon of points simultaneously approximating themselves is well-studied. Following Furstenberg \cite[pg. 9]{furstenberg_book_1981}, a point in $X$ is \emph{multiply recurrent for $\mapspace$} precisely when it simultaneously approximates itself under $\mapspace$.  When the maps in $\mapspace$ are mutually commuting, multiply recurrent points are guaranteed to exist.  This theorem of Furstenberg and Weiss is one of the foundational results in the study of recurrence in topological dynamics as it pertains to the present work.

\begin{theorem}[{cf. \cite[Thm. 1.4]{furstenberg_weiss_1978} and \cite[Thm. 2.6]{furstenberg_book_1981}}]
\label{thm_furst_mult_recurrence_intro}
Let $\mapspace$ be a countable collection of commuting, continuous maps of a compact metric space $X$ to itself.  There exists a point which is multiply recurrent for $\mapspace$.
\end{theorem}

More broadly, the phenomenon of recurrence in dynamical systems appeared first in the work of Poincar\'e in the 19$^{\text{th}}$ century (see \cite{furstenberg_1981,bergelson_2000}) and remains an active area of research in topological, measurable, homogeneous, and linear dynamics.  One reason for the topic's longevity is its fruitful relationship with combinatorics and additive number theory, from Bogolyubov's work \cite{bogolyubov_1939} linking the combinatorics of iterated difference sets with Bohr almost periodic structure to Furstenberg's proof \cite{furstenberg_1977} of Szemer\'edi's theorem via recurrence in ergodic theory, and beyond \cite{furstenberg_book_1981,frantzikinakis_mccutcheon_2009}.

While much is known about multiple recurrence in various general settings, relatively little appears to be written about the phenomenon of points simultaneously approximating other points, even in the simplest setting of a \emph{(topological dynamical) system} $(X,T)$, a compact metric space $X$ together with a continuous map $T: X \to X$.  We aim to address this here, and, following in the tradition described above, we aim to derive combinatorial corollaries from the dynamical results.

Simultaneous approximation appears as a corollary of a result of Glasner's \cite[Cor. 2.5]{glasner_1994} in the following way.  If $(X,T)$ is a minimal, weakly mixing system, then there exists a residual set of points in the product $X^d$ that have a dense $T \times T^2 \times \cdots \times T^d$-orbit.\footnote{A \emph{residual} set is one that contains a dense $G_\delta$ set.  A system $(X,T)$ is \emph{minimal} if every point has a dense orbit, \emph{ergodic} if a residual set of points have a dense orbit, and \emph{weakly mixing} if the product system $(X^2, T \times T)$ is ergodic. If $(X,T)$ is weakly mixing, then so is the system $(X^d, T \times T^2 \times \cdots \times T^d)$, and weakly mixing systems are ergodic.}  It was shown by Glasner that, in fact, a residual set of points $(x, \ldots, x)$ of the diagonal in $X^d$ have a dense $T \times T^2 \times \cdots \times T^d$-orbit.  This can be rephrased as simultaneous approximation as follows.

\begin{theorem}[{cf. \cite[Cor. 2.5]{glasner_1994}}]
\label{thm_glasner}
Let $(X,T)$ be a minimal, weakly mixing system, and define $\mapspace = \big\{T^i \ | \ i \in \N \big\}$.  There exists a residual set of points that simultaneously approximate all points in $X$ under $\mapspace$.
\end{theorem}

It is an exercise to check that if $(X,T)$ is a rotation on a compact group, then a point $x$ simultaneously approximates only itself under $\{T,T^2\}$.  Therefore, if a point $x$ simultaneously approximates $y$ under $\{T, T^2\}$, then $x$ and $y$ lie in the same fiber over the \emph{(topological) Kronecker factor of $(X,T)$}, the largest equicontinuous factor of $(X,T)$.  A system is weakly mixing if and only if its Kronecker factor is the one-point system.  Thus, the weakly mixing assumption in \cref{thm_glasner} is necessary for its conclusion.

\cref{thm_furst_mult_recurrence_intro} gives that the phenomenon of ``self'' simultaneous approximation is present in all systems, while \cref{thm_glasner} gives that ``complete'' simultaneous approximation is typical in sufficiently mixing systems.  Our first main result, \cref{maintheorem_nilsystems_have_simult_approx}, aims at developing the space between these two theorems by specifying when the set of points simultaneously approximated by a point is $\delta$-dense.

To ease the discussion, we introduce some notation.  Let $(X,T)$ be a system.  For $x \in X$ and $A \subseteq \N$, we define the set $\sa(x,A)$ to be the set of points $y \in X$ that are simultaneously approximated by $x$ under $\{T^n \ | \ n \in A\}$:
\begin{align}
\label{eqn_def_of_sa_set_intro}
    \sa(x,A) \defeq \left\{ y \in X \ \middle| \ \begin{gathered} \forall \text{ finite } F \subseteq A, \ \forall \eps > 0, \\ \exists m \in \N, \ \forall f \in F, \ d_X(T^{fm} x, y) < \eps \end{gathered} \right\}.
\end{align}
If $(X,T)$ is invertible, the definition of $\sa(x,A)$ makes sense for any $A \subseteq \Z$. Which points belong to the closed set $\sa(x,A)$ is a function of $x$, $A$, and the particular dynamical properties of the system $(X,T)$.  We will focus primarily on the denseness of $\sa(x,A)$ in $X$ as a measure of its quality.

Explicit computations are instructive in simple examples.  In an irrational rotation $T: x \mapsto x + \alpha$ of the 1-torus $\T \defeq \R / \Z$, for all $A \subseteq \Z$, the set $\sa(0,A)$ can be shown to be a closed subgroup of $\T$ and $\sa(x,A) = x + \sa(0,A)$.  More precisely, it is shown in \cref{lem_description_of_sa_set_on_torus} that
\[\sa(x,A) = x + \big\{0, 1/N, \ldots, (N-1)/N \big\}, \text{ where $N = \gcd(A-A) / \gcd A$},\]
where $A - A \defeq \{ a - a' \ | \ a, a' \in A\}$, and $\gcd$ denotes the greatest common divisor (of a set of possibly infinitely many integers). It is shown in \cref{lem_delta_dense_implies_A_in_prog,lemma_special_case_of_main_thm} that if $\gcd A = 1$, then $\sa(x,A)$ is a roughly $N^{-1}$-dense subset of $\T$ if and only if $A$ is contained in an infinite arithmetic progression of the form $N \Z + n$, where $N \in \N$ and $n \in \Z$ is coprime to $N$.  The full argument is elementary and is spelled out in \cref{section_simul_approx_for_irrational_rotations}.

If $\pi: (X,T) \to (Y,T)$ is a factor map of systems (a continuous surjection for which $\pi \circ T = T \circ \pi$), it is easy to see that $\pi \sa(x,A) \subseteq \sa(\pi x,A)$. Thus, the arithmetic condition $A \subseteq N \Z + n$ necessary for the roughly $N^{-1}$-denseness of $\sa(x,A)$ in an irrational rotation applies equally well to any system for which the irrational rotation is a factor. The following definition is motivated, then, by the guess that such an arithmetic condition on $A$ may be sufficient for the $\delta$-denseness of $\sa(x,A)$ in a general system.

\begin{definition}
\label{def_dense_simult_approx}
Let $(X,T)$ be an invertible system.
\begin{itemize}
    \item A point $x \in X$ possesses the \emph{dense simultaneous approximation property} if for all $\delta > 0$, there exists $N \in \N$ such that for all $n \in \Z$ coprime to $N$, the set $\sa(x,N \Z + n)$ is $\delta$-dense in $X$.
    \item The system $(X,T)$ possesses the \emph{almost everywhere (a.e.) dense simultaneous approximation property} if there exists a residual set $\Omega \subseteq X$ such that for all $\delta > 0$, there exists $N \in \N$ such that for all $x \in \Omega$ and $n \in \Z$ coprime to $N$, the set $\sa(x, N \Z + n)$ is $\delta$-dense in $X$.
\end{itemize}
\end{definition}

Both properties in \cref{def_dense_simult_approx} are present, by computations similar to those in \cref{section_simul_approx_for_irrational_rotations}, in rotations of compact abelian groups and, by \cref{thm_glasner}, in invertible, minimal, weakly mixing systems.  Our first main result demonstrates that minimal nilsystems also exhibit both of the properties.  A \emph{nilsystem} $(X,T)$ is a compact quotient $X = G / \Gamma$ of a nilpotent Lie group $G$ by a discrete, cocompact subgroup $\Gamma$ together with a \emph{nilrotation} $T: x \mapsto gx$, where $g \in G$.  Nilsystems, which are always invertible and are defined more precisely in \cref{sec_nilsystems}, generalize compact group rotations and have played a central role in the development of the structure theories of topological and measurable dynamical systems \cite{host_kra_book_2018} and in higher-order Fourier analysis \cite{tao_hofa_book}.

\begin{Maintheorem}
\label{maintheorem_nilsystems_have_simult_approx}
Let $(X,T)$ be a minimal nilsystem.
\begin{enumerate}[label=(\Roman*)]
    \item \label{item_points_in_nilsystems_sat_dense_approx} All points $x \in X$ satisfy the dense simultaneous approximation property.
    \item \label{item_nilsystem_sats_ae_dense_approx} The system $(X,T)$ satisfies the a.e. dense simultaneous approximation property.
\end{enumerate}
\end{Maintheorem}

The proof of \cref{maintheorem_nilsystems_have_simult_approx} will follow from a more general result, \cref{thm_main_fullversion}, that provides a bevy of additional information.  Generalizing the argument for an irrational rotation of the torus given in \cref{lemma_special_case_of_main_thm}, we show that certain rational points in the product nilmanifold $X^{d+1}$, under a nilrotation and a projection, yield a $\delta$-dense subset of simultaneously approximable points in $X$.  While nilsystems are known to inherit certain properties from their Kronecker factors,\footnote{A nilsystem is minimal if and only if its Kronecker factor is minimal \cite[Ch. 11, Thm. 6]{host_kra_book_2018}. Two nilsystems are measurably disjoint if and only if their Kronecker factors are measurably disjoint \cite[Thm. 3.1 \& 3.5]{berg_1971}. If $X = G / \Gamma$ with $G$ connected, the orthocomplement in $L^2(X)$ of the square-integrable functions on the Kronecker factor has Lebesgue spectrum \cite[Ch. V, Thm. 4.2]{auslander_green_hahn_1963_book}.} the example in \cref{rmk_kronecker_is_not_characteristic} shows that the Kronecker factor is not ``characteristic'' for simultaneous approximation.  The argument behind \cref{maintheorem_nilsystems_have_simult_approx} ultimately relies on the theory of rational points in nilmanifolds developed by Leibman \cite{leibman_rational_2006} -- adopting some of the quantitative ideas from Green and Tao \cite{green_tao_quantitative_behaviour_2012} -- and the theory of polynomial orbits in nilmanifolds \cite[Ch. 14]{host_kra_book_2018} developed by those aforementioned and others.

The arguments required to establish the a.e. simultaneous approximation property for nilsystems occupy a significant portion of the paper.  Under the assumption of minimality, the existence of one multiply recurrent point in \cref{thm_furst_mult_recurrence_intro} easily implies the existence of a residual set of multiply recurrent points.  The same is not true for simultaneous approximation: while it is true that $T^n\sa(x,A) = \sa(T^nx, A)$ when $T$ is invertible, the quality of $\delta$-denseness is generally not preserved under the map $T^n$.  Instead, to establish \ref{item_nilsystem_sats_ae_dense_approx}, we make significant use of the \emph{prolongation relation} \cite{auslander_guerin_1997}, the closure of the orbit-closure relation defined more precisely in \cref{sec_prolongation_relation}. We describe the residual set appearing in \ref{item_nilsystem_sats_ae_dense_approx} explicitly as the set of points $x \in X$ for which the point $(x,\ldots,x)$ is a point of continuity of the $(T^0 \times T^1 \times \cdots \times T^d)$-orbit-closure map.  It was shown by Akin and Glasner \cite{akin_glasner_1998} that the orbit closure at such a point coincides with the point's prolongation class.  We demonstrate in \cref{thm_continuity_of_intersections} that the map that sends a point to the intersection of its prolongation classes under two commuting, distal homeomorphisms is continuous.  This continuity is then a key tool in proving \ref{item_nilsystem_sats_ae_dense_approx}.

We show in \cref{sec_dense_approx_in_skew_prods} that there are non-nilsystems that have points possessing the dense simultaneous approximation property, but we do not know the extent to which the properties in \cref{def_dense_simult_approx} are exhibited in general systems.  In an effort to excite further investigation in this direction, we elaborate on this question and others in \cref{sec_open_ends}.

\subsection{Multiplicative thickness of return-time sets}
\label{sec_intro_thickness}

Let $(X,T)$ be a system, $x \in X$, and $U \subseteq X$.  The \emph{return-time set} $R(x,U)$ is the set of times $m \in \N$ at which the point $x$ visits the set $U$:
\[R(x,U) = \big\{ m \in \N \ \big| \ T^m x \in U \big\}.\]
Dynamical qualities of a system are often described in terms of combinatorial qualities of its return-time sets.  For example, the system $(X,T)$ is minimal if and only if for all $x \in X$ and all non-empty, open $U \subseteq X$, the set $R(x,U)$ is \emph{additively syndetic}: there exists $N \in \N$ for which the set $R(x,U)$ has non-empty intersection with every set of $N$ consecutive positive integers.  In fact, the notion of additive syndeticity has its roots in topological dynamics \cite[Ch. 2]{gottschalk_hedlund_1955}.

The notion of multiple recurrence appearing in \cref{thm_furst_mult_recurrence_intro} can be equivalently formulated in terms of the multiplicative quality of return-time sets.  Indeed, a point $x$ is multiply recurrent under $\{T^n \ | \ n \in \N\}$ if and only if for all open neighborhoods $U$ of $x$, the set $R(x,U)$ is \emph{multiplicatively thick} in the semigroup $(\N, \cdot)$: for all finite subsets $F \subseteq \N$, there exists $m \in \N$ such that $mF \subseteq R(x,U)$.  Thickness is a combinatorial ``dual'' of syndeticity, a notion made precise by the algebra of Furstenberg-families \cite[Ch. 9]{furstenberg_book_1981}.  In this framing of multiple recurrence, the following is an equivalent form of the single-transformation version of a theorem of Furstenberg and Weiss \cite[Thm. 1.5]{furstenberg_weiss_1978}, a straightforward consequence of \cref{thm_furst_mult_recurrence_intro}.

\begin{theorem}[{cf. \cite[Thm. 1.5]{furstenberg_weiss_1978}}]
\label{thm_furst_weiss}
Let $(X,T)$ be a minimal system.  There exists a residual, $T$-invariant set $\Omega \subseteq X$ for which the following holds.  For all non-empty, open $U \subseteq X$ and all $x \in \Omega \cap U$, the set $R(x,U)$ is multiplicatively thick in $(\N,\cdot)$.
\end{theorem}

The notion of simultaneous approximation can also be described naturally in terms of the multiplicative quality of return-time sets: the point $x$ simultaneously approximates $y$ under $\{T^n \ | \ n \in A\}$, $A \subseteq \N$, if and only if for all finite $F \subseteq A$ and all open neighborhoods $U$ of $y$, there exists $m \in \N$ such that $mF \subseteq R(x,U)$.  Thus, we aim to generalize \cref{thm_furst_weiss} by demonstrating the multiplicative thickness of all return-time sets $R(x,U)$, not just those for which $U$ is a neighborhood of $x$.

When $U$ is not a neighborhood of $x$, there are some natural local obstructions that must be addressed.  The set of odd natural numbers, $2 \N - 1$, is the return-time set of one point to the other in a rotation on two points.  In the hierarchy of notions of largeness that has developed hand-in-hand with ergodic Ramsey theory \cite{hindman_jones_strauss_2019}, the set $2 \N - 1$ is of a poor multiplicatively quality in the semigroup $(\N,\cdot)$: it has zero mean with respect to all dilation-invariant means on $\N$ (ie., it has zero \emph{upper Banach density in the semigroup $(\N,\cdot)$}) because every such mean is supported on the even natural numbers.  The odd numbers, nevertheless, are multiplicatively closed and thus have a rich multiplicative quality relative to themselves as the semigroup $(2\N - 1,\cdot)$.

We define a \emph{congruence subsemigroup of $(\N,\cdot)$} to be a subsemigroup of $(\N,\cdot)$ of the form
\[\nn \defeq \big\{ n \in \N \ \big| \ n \equiv 1 \pmod N \big\} = \big\{1, N + 1, 2N + 1, \ldots \big\}, \qquad N \in \N.\]
A \emph{coset} of the congruence subsemigroup $\nn$ is a set of the form $I \nn \defeq \{I, I(N+1), I(2N+1), \ldots \}$, where $I \in \N$.  Every infinite arithmetic progression contains a coset of a congruence subsemigroup: $N \N + I$ contains $(N+I) \nn$.  As multiplicatively rich sets, cosets of congruence subsemigroups allow us to describe some of the multiplicative qualities of infinite arithmetic progressions and, hence, of return-time sets in rotations on finitely many points.

\cref{maintheorem_ae_return_set_is_mult_thick} shows that the obstruction to the multiplicative thickness of return-time sets presented by rotations on finitely many points is, in fact, the only obstruction in general systems.  More precisely, the typical return-time set $R(x,U)$ is \emph{multiplicatively thick in a coset of a congruence subsemigroup $I\nn$}: for all finite $F \subseteq \nn$, there exists $m \in \nn$ such that $ImF \subseteq R(x,U)$.  In fact, \cref{maintheorem_ae_return_set_is_mult_thick} generalizes \cref{thm_furst_weiss}, as we will explain in the next subsection.

\begin{Maintheorem}
\label{maintheorem_ae_return_set_is_mult_thick}
Let $(X,T)$ be a minimal system.  There exists a residual, $T$-invariant set $\Omega \subseteq X$ for which the following holds.  For all non-empty, open $U \subseteq X$, there exists $N \in \N$ such that for all $x \in \Omega$, there exists $I \in \{1, \ldots, N\}$ such that the set $R(x,U)$ is multiplicatively thick in $I\nn$.
\end{Maintheorem}

This theorem will follow from a more general result, \cref{theorem_advanced_main_theorem_on_returntime_sets}, which elaborates on relationships between the point $x$, the set $U$, and the coset of the congruence subsemigroup $I\nn$ appearing in the conclusion.  For example, the translated return-time set $R(x,U) - t = R(T^t x,U)$ is multiplicatively thick in $(I-t)\nn$, and, if $(X,T)$ is totally minimal,\footnote{A system $(X,T)$ is \emph{totally minimal} if for all $n \in \Z \setminus \{0\}$, the system $(X,T^n)$ is minimal.} the word ``coset'' can be eliminated so that the return-time set $R(x,U)$ is multiplicatively thick in a congruence subsemigroup.

Koutsogiannis, Richter, and the author \cite{glasscock_koutsogiannis_richter_2019} initiated an investigation into the multiplicative combinatorial properties of return-time sets in topological dynamical systems.  Though the techniques are substantially different, the present work relies on several lemmas from that paper and generalizes several of the main results.

The main body of work that goes into proving \cref{maintheorem_nilsystems_have_simult_approx} proves a version of \cref{maintheorem_ae_return_set_is_mult_thick} for nilsystems and for inverse limits of nilsystems.  (\cref{mainthm_nilbohr_sets_are_mult_thick} shows that all nil-Bohr sets are multiplicatively thick in a congruence subsemigroup, the combinatorial consequences of which are discussed in the next section.)  The extension from \cref{maintheorem_ae_return_set_is_mult_thick} for inverse limits of nilsystems to general systems is accomplished by appealing to the machinery of topological characteristic factors, initiated by Glasner \cite{glasner_1994} and developed recently by Glasner, Huang, Shao, Weiss, and Ye \cite{glasner_huang_shao_weiss_ye_2020}.  Roughly speaking, we prove that the \emph{infinite-step pro-nilfactor} is characteristic for the property of return-time sets being multiplicatively thick in a coset of a congruence subsemigroup.

\subsection{Combinatorial applications}
\label{sec_intro_comb_applications}

Since the pioneering work of Furstenberg \cite{furstenberg_1977} and Furstenberg and Weiss \cite{furstenberg_weiss_1978}, measurable and topological dynamics have proven to be remarkably effective in addressing certain problems from additive combinatorics and Ramsey theory.  It is in this vein that we will show how Theorems \ref{maintheorem_nilsystems_have_simult_approx} and \ref{maintheorem_ae_return_set_is_mult_thick} and the scaffolding that supports them can be used to give new combinatorial results in two directions: enhancements to van der Waerden's theorem on arithmetic progressions and multiplicative structures in translates of difference sets and their higher-order recurrence analogues. Finally, in \cref{sec_disjointness}, we give a result concerning the disjointness between additive and multiplicative dynamical systems.

We begin with an reformulation of \cref{thm_furst_weiss} that is more convenient for our combinatorial purposes.  An elementary argument shows that Theorems \ref{thm_furst_weiss} and \ref{thm:dvdw} are equivalent; see \cite[Sec. 5]{glasscock_koutsogiannis_richter_2019}.

\begin{theorem}[{cf. \cite[Theorem 1.5]{furstenberg_weiss_1978}}]
\label{thm:dvdw}
Let $(X,T)$ be a minimal system. For all non-empty, open $U \subseteq X$ and $n_1, \ldots, n_k \in \N$, there exists $m \in \N$ such that
\begin{align}
\label{eqn:topologicalmultiplerecurrence}
U \cap T^{-n_1m} U \cap T^{-n_2m} U \cap \cdots \cap T^{-n_k m} U \neq \emptyset.
\end{align}
\end{theorem}

A topological correspondence principle is a device that links syndetic subsets $A \subseteq \N$ and addition with non-empty, open subsets $U$ of a compact space $X$ and a continuous transformation $T: X \to X$; in particular, expressions of the form in \eqref{eqn:topologicalmultiplerecurrence} are tied to expressions of the form
\begin{align}
\label{eqn:combinatorialmultiplerecurrence}
A \cap \big( A -n_1m \big) \cap \big(A -n_2m \big) \cap \cdots \cap \big(A-n_k m \big).
\end{align}
Note that when $n_1 = 1$, $n_2 = 2$, \dots, $n_k = k$, every element of the set in \eqref{eqn:combinatorialmultiplerecurrence} starts an arithmetic progression of length $k+1$ in $A$.

In one of the earliest results in Ramsey Theory, van der Waerden \cite{vanderWaerden_1927} proved that one piece of any finite partition of $\N$ necessarily contains arbitrarily long arithmetic progressions.  It is an equivalent formulation due to Kakeya and Morimoto \cite[Theorem I]{Kakeya_Morimoto_1930} that every syndetic subset of $\N$ contains arbitrarily long arithmetic progressions.  Thus, combined with a correspondence principle, \cref{thm:dvdw} implies van der Waerden's theorem.  This dynamical framing allowed Furstenberg and Weiss to describe combinatorial refinements of van der Waerden's theorem and several other foundational results in Ramsey Theory.

Just as Theorems \ref{thm_furst_weiss} and \ref{thm:dvdw} are equivalent, the following theorem is an equivalent form of \cref{maintheorem_ae_return_set_is_mult_thick}.  We demonstrate the equivalence between Theorems \ref{maintheorem_ae_return_set_is_mult_thick} and \ref{thm_comb_equivalent_to_main_thm} in \cref{sec_refinement_on_vdw}.

\begin{Maintheorem}
\label{thm_comb_equivalent_to_main_thm}
Let $(X,T)$ be a minimal system.  For all non-empty, open $U \subseteq X$, there exists $N \in \N$ such that for all non-empty, open $V \subseteq X$, there exists $I \in \{1, \ldots, N\}$ such that for all $n_1, \ldots, n_k \in I\nn$, there exists $m \in \nn$ such that
\begin{align}
    \label{eqn_set_eqn_in_comb_equiv_thm}
    V \cap T^{-n_1m} U \cap T^{-n_2m} U \cap \cdots \cap T^{-n_km} U \neq \emptyset.
\end{align}
\end{Maintheorem}

\cref{thm_comb_equivalent_to_main_thm} is a generalization of \cref{thm:dvdw}.  To derive the latter from the former, apply \cref{thm_comb_equivalent_to_main_thm} with $V = X$ and with dilates $I, I(N n_1 + 1), \ldots, I(N n_k + 1)$, and note that \eqref{eqn_set_eqn_in_comb_equiv_thm} implies that \eqref{eqn:topologicalmultiplerecurrence} holds with $INm$ as $m$.  Because Theorems \ref{maintheorem_ae_return_set_is_mult_thick} and \ref{thm_comb_equivalent_to_main_thm} are equivalent, this verifies the claim in the previous subsection that \cref{thm_furst_weiss} follows from \cref{maintheorem_ae_return_set_is_mult_thick}.

Combined with a topological correspondence principle, \cref{thm_comb_equivalent_to_main_thm} leads to an enhancement of van der Waerden's theorem.  Given an set $A \subseteq \N$ and $k \in \N$, call a positive integer $b$ an \emph{$\nn$-starter} if there exists $m \in \N$ such that
\[\big\{ b + m, \ b + (N+1)m, \ b + (2N + 1)m, \ \dots, \ b + (kN+1)m \big\} \subseteq A.\]
Thus, the integer $b$ is an $\nn$-starter if the set $A-b$ contains a dilate of an initial subset of the semigroup $(\nn,\cdot)$.
\cref{mainthm_combinatorial_main} gives that the $\nn$-starters for an additively piecewise syndetic set form an additively thick set.  An \emph{additively thick} set is one that contains arbitrarily long intervals, and an \emph{additively piecewise syndetic} set is the intersection of an additively syndetic set and an additively thick set.

\begin{Maintheorem}
\label{mainthm_combinatorial_main}
\label{main_vdw_type_theorem}
Let $A \subseteq \N$ be additively piecewise syndetic.  There exists $N \in \N$ for which the following holds.  For all additively syndetic $B \subseteq \N$, there exists $I \in \{1, \ldots, N\}$ such that for all $n_1, \ldots, n_k \in I\nn$, there exists $m \in \nn$ such that
\begin{align}
\label{eqn_set_intersection_in_combo_corollary}
    B \cap (A-n_1m) \cap \big(A-n_2m\big) \cap \cdots \cap \big(A-n_km \big) \neq \emptyset.
\end{align}
In particular, for all $k \in \N$, the set of $\nn$-starters
\begin{align}
\label{eqn_nn_starters}
    \big\{ b \in \N \ \big | \ \exists m \in \N, \ b + m, \ b + (N+1)m, \ \ldots, \ b+(k N +1)m \in A \big\}
\end{align}
is additively thick.
\end{Maintheorem}

\cref{mainthm_combinatorial_main} implies van der Waerden's theorem.  Indeed, this follows immediately from the fact that if $b$ is an $\nn$-starter with dilate $m$, then $b+m$ is the first element of an arithmetic progression in the set $A$.  Even the simplest periodic sets $A$ demonstrate that there are arithmetic obstructions to the set of starters of arithmetic progressions contained in $A$ being additively thick.  \cref{mainthm_combinatorial_main} shows that this obstruction is resolved by passing to dilates of finite subsets of the right congruence subsemigroup $\nn$.

\cref{mainthm_combinatorial_main} is derived from \cref{thm_comb_equivalent_to_main_thm} in \cref{sec_refinement_on_vdw}. \cref{thm_comb_equivalent_to_main_thm} is derived from \cref{maintheorem_ae_return_set_is_mult_thick}, which relies on all of the machinery in this paper.  Thus, while the statement of \cref{mainthm_combinatorial_main} is quite simple, the proof presented here is quite complicated. It would be meaningful to find alternative approaches to the theorem.  Of particular interest would be a combinatorial understanding of the quantity $N$ in relation to the piecewise syndetic set $A$.\\

The second combinatorial application concerns translates of difference sets and their higher-order recurrence analogues. Sum and difference sets arise naturally as supports of convolutions, and as such, their structure has been the subject of much investigation in additive combinatorics.  Combining several elementary results in Ramsey Theory, it can be shown that when a set $A \subseteq \Z$ has positive \emph{additive upper (Banach) density},
\begin{align}
\label{eqn_add_up_banach_den}
    d^*(A) \defeq \limsup_{N \to \infty} \max_{n \in \Z} \frac{ \big|A \cap \{n, \ldots, n+N-1\}\big|}{N},
\end{align}
its difference set, $A-A = \{a - a' \ | \ a, a' \in A\}$, is multiplicatively thick in the semigroup $(\Z\setminus\{0\},\cdot)$.\footnote{First, note that the set $A-A$ is a $\Delta^*$ set, meaning that it has non-empty intersection with every set of the form $E-E$, where $E \subseteq \Z$ is infinite.  This can be seen by the pigeonhole principle: each of the sets $A-e$, $e \in E$, has the same positive density along the same sequence of intervals, and hence there must exist distinct $e,e' \in E$ for which $(A-e) \cap (A-e') \neq \emptyset$.  This implies that $e-e' \in A-A$, whereby $(A-A) \cap (E-E) \neq \emptyset$.  A similar argument shows that for all $n \in \Z \setminus \{0\}$, the set $(A-A)/n$ is a $\Delta^*$ set.  Next, the intersection of two $\Delta^*$ sets is a $\Delta^*$ set.  Indeed, following the definitions, it is equivalent to show that if $E \subseteq \Z$ is infinite and $E-E = B_1 \cup B_2$, then at least one of the $B_i$'s contains $D-D$ for some infinite $D \subseteq \Z$.  This is a well-known consequence of Ramsey's Theorem.  Finally, to see that $A-A$ is multiplicatively thick in $(\Z\setminus\{0\},\cdot)$, we need only to show that for $n_1, \ldots, n_k \in \Z\setminus\{0\}$, the set $\cap_{i=1}^k (A-A)/n_i$ is non-empty.  But this set is $\Delta^*$, hence non-empty, as desired.}  Succinctly, the set $A-A$ can be shown to be multiplicatively rich by virtue of its additive largeness; this idea and its consequences were explored in \cite{bergelson_glasscock_2018,bergelson_glasscock_2020}. Similar multiplicative conclusions for translates $A-A+t$ of the difference set, however, do not seem so easily achievable via such combinatorial arguments.

In light of the work in this paper, another description of the set $A-A$ is more useful in understanding the multiplicative qualities of its translates.  A \emph{Bohr} set is one that contains the return-time set of a point to a non-empty, open set in a compact group rotation.  A \emph{Bohr set up to zero additive density} is a set of the form $B \cap D$, where $B$ is a Bohr set and $D$ is a set for which $d^*(\Z \setminus D) = 0$.  It was shown by \Folner{} \cite{folnerbogoliouboff1954,folnernoteonbogoliouboff1954} that if $d^*(A) > 0$, then the set $A-A$ is a Bohr set up to zero additive density.  As a consequence of the work leading up to \cref{maintheorem_ae_return_set_is_mult_thick}, we have that Bohr sets and their nilpotent analogues, \emph{nil-Bohr} sets (see \cref{def_nilbohr_set}), up to zero additive density are multiplicatively thick in a coset of a congruence subsemigroup in $\N$ and in $\Z$ (as defined at the beginning of \cref{sec_comb_structure_in_nilbohr_sets}).

\begin{Maintheorem}
\label{mainthm_nilbohr_sets_are_mult_thick}
All nil-Bohr sets, including those up to zero additive density, are multiplicatively thick in a coset of a congruence subsemigroup in $\N$ and in $\Z$.
\end{Maintheorem}

\cref{mainthm_nilbohr_sets_are_mult_thick} implies that sets of the form $A - A + t$, where $A \subseteq \Z$ has positive additive density, are multiplicatively thick in a coset of a congruence subsemigroup.  To describe the higher-order analogues of this result, it is helpful to see the connection between difference sets and recurrence.  Note that $n \in A-A$ if and only if $A \cap (A-n) \neq \emptyset$, whereby elements of the difference set are exactly the times at which the set $A$ ``returns to itself.''  When $\lambda$ is a translation-invariant mean on $\Z$, the \emph{correlation sequence} $n \mapsto \lambda \big(A \cap (A-n) \big)$ is positive definite.  It follows by combining classical results of Herglotz and Weiner that the sequence can be written as a sum of an almost periodic sequence and a null-sequence, defined in \cref{sec_comb_applications}.  This explains the origin of the ``Bohr up to zero density'' structure in \Folner{}'s result.

Much more recently, \emph{multi-correlation sequences} of the form $n \mapsto \lambda \big( A \cap (A-n) \cap \cdots \cap (A-kn) \big)$ were described by Bergelson, Host, and Kra \cite{bergelson_host_kra_2005} as a sum of a nilsequence and a null-sequence.  Appealing to \cref{mainthm_nilbohr_sets_are_mult_thick}, it follows that sets of the form
\[\big\{n \in \Z \ \big| \ A \cap (A-n) \cap \cdots \cap (A-kn) \neq \emptyset \big\} + t,\]
where $A \subseteq \Z$ has positive additive density, are multiplicatively thick in a coset of a congruence subsemigroup. Appealing to more general decomposition results extends the application, leading to curious new affine configurations in sets of positive density. The full extent of the application is given in \cref{cor_correlation_sequences_are_mult_thick} and the discussion following it.

Multiplicative thickness in a coset of a congruence subsemigroup follows from
\cref{mainthm_nilbohr_sets_are_mult_thick} wherever nil-Bohr structure is found.  Algebraic skew-product constructions of Furstenberg \cite{furstenberg_book_1981} show that non-empty sets of the form
\begin{align}
\label{eqn_poly_returns_to_an_inteveral}
    \big\{n \in \Z \ \big| \ \{ p(n) \} \in I \big\},
\end{align}
where $p \in \R[x]$ is a polynomial, $\{ \cdot \}: \R \to [0,1)$ denotes the factor (fractional part) map, and $I \subseteq [0,1)$ is non-empty and open, are nil-Bohr sets.  By work of Bergelson and Leibman \cite{bergelson_leibman_2007}, the same remains true when $p$ is taken from the much broader class of \emph{generalized polynomials}, described in \cref{sec_comb_applications}.  When $p(0) = 0$ and 0 is an interior point of $I$ (when $0$ and $1$ are identified), sets of the form in \eqref{eqn_poly_returns_to_an_inteveral} are known to be multiplicatively rich by virtue of their additive largeness, as mentioned above.  Dropping the assumptions on $p$ and $I$, we can appeal to \cref{mainthm_nilbohr_sets_are_mult_thick} to see that sets of the form in \eqref{eqn_poly_returns_to_an_inteveral} are multiplicatively thick in a coset of a congruence subsemigroup.  The full scope of this application is given in \cref{cor_bdd_gen_poly_mult_thick}.\\

Not all additively syndetic sets in $\N$ arise as return-time sets in minimal topological dynamical systems \cite[Sec. 8]{glasscock_koutsogiannis_richter_2019}.  Nevertheless, the class of return-time sets serves as a useful proxy when attempting to gain insight into the combinatorial nature of syndetic sets. To the extent that multiplicative configurations in return-time sets indicate the same configurations in syndetic sets, \cref{maintheorem_ae_return_set_is_mult_thick} offers tantalizing new evidence toward a positive answer to some recalcitrant problems regarding the multiplicative structure of additively syndetic sets.  It cannot, however, be generalized naively: we construct an additively syndetic subset of $\N$ in \cref{sec_add_synd_not_mult_thick} that is not multiplicatively thick in any coset of any congruence subsemigroup.

Beiglb\"{o}ck, Bergelson, Hindman, and Strauss \cite{beiglbock_bergelson_hindman_strauss_2006} asked whether or not every additively syndetic subset of $\N$ contains arbitrarily long geometric progressions.  Their question remains unanswered; in fact, it is still not known whether or not an additively syndetic set must contain two integers whose ratio is a square.  Thus, there remains a very wide gulf of possibilities for additively syndetic sets not arising from dynamics, from admitting a square ratio to multiplicative piecewise syndeticity in a coset of a congruence subsemigroup.  We lay out some concrete questions that we believe are the next in line to be considered in \cref{subsec_comb_explorations}.

\subsection{Organization of the paper}

The paper is organized as follows.  In \cref{sec_prolong}, we establish notation and terminology used throughout the work.  We derive auxiliary results on continuity of the intersections of prolongation classes (\cref{thm_continuity_of_intersections}) and on the uniform-denseness of prolongation-class-measurable sets (\cref{thm_uniform_denseness_of_msble_sets}).
    
In \cref{sec_rational_elts_in_groups}, we collect the results we need on rational elements and polynomials in nilpotent groups and nilmanifolds.  The main result is a quantitative version of the fact that rational polynomial orbits in nilmanifolds are periodic (\cref{thm_any_equivalent_poly_is_periodic}).
    
In \cref{section_simultaneous_approx_in_nilsystems}, we describe the phenomenon of simultaneous approximation in an irrational rotation of the 1-torus (\cref{section_simul_approx_for_irrational_rotations}), then we prove  \cref{maintheorem_nilsystems_have_simult_approx} in \cref{sec_proof_of_main_thm} by deriving it from the more general \cref{thm_main_fullversion}.
    
In \cref{sec_comb_structure_in_nilbohr_sets}, we demonstrate the multiplicative thickness of return-time sets in nilsystems (\cref{thm_strengthening_of_theorem_c}), inverse limits of nilsystems (\cref{thm_strengthening_of_theorem_c_for_inverse_limits}), and general systems (\cref{thm_non_invert_main_result}), and deduce several combinatorial corollaries in Sections \ref{sec_comb_applications} and \ref{sec_refinement_on_vdw}.  We prove Theorems \ref{maintheorem_ae_return_set_is_mult_thick}, \ref{thm_comb_equivalent_to_main_thm}, \ref{main_vdw_type_theorem}, and \ref{mainthm_nilbohr_sets_are_mult_thick} in Sections \ref{sec_proof_of_return_sets_mult_thick}, \ref{sec_refinement_on_vdw}, \ref{sec_refinement_on_vdw}, and \ref{sec_mult_thick_in_nilsystems}, respectively.
    
In \cref{sec_final_section}, we demonstrate that the phenomenon of simultaneous approximation occurs outside the class of nilsystems (\cref{theorem_simult_approx_is_typical}), and we construct an additively syndetic set that is nowhere multiplicatively thick (\cref{theorem_add_syndetic_not_mult_thick}).  We conclude by collecting some questions for future consideration in \cref{sec_open_ends}.

\subsection*{Acknowledgements}

The author is indebted to Andreas Koutsogiannis, Joel Moreira, and Florian Richter for helpful discussions concerning early versions of this work.  Thanks also goes to Sasha Leibman for correspondence on the nature of the uniform denseness of rational points across subnilmanifolds of a nilmanifold.

\section{Prolongations, intersections, and uniform denseness}
\label{sec_prolong}

In \cref{sec_top_set_and_systems}, we lay out the basic notation and terminology that will be used throughout the work.  The remainder of this section -- devoted to proving a continuity result regarding prolongation classes -- can be safely skipped on a first reading and be used as a reference as needed.

Denote the positive integers and integers by $\N$ and $\Z$, respectively, and write $\Nz \defeq \{0\} \cup \N$.  The notation $( x_n)_n \subseteq X$ will be used to indicate a sequence indexed by $\N$ with range in $X$.  If $X$ is a set and $T: X \to X$ is a map, a subset $\Omega \subseteq X$ is \emph{$T$-invariant} if $T \Omega \subseteq \Omega$.  Given an equivalence relation on $X$, a subset $Y \subseteq X$ is \emph{measurable} with respect to that relation if $Y$ is a union of equivalence classes. For $A \subseteq \N$ and $n \in \N$, we define
\[A - n \defeq \big\{ m \in \N \ \big| \ m + n \in A \big\} \quad \text{ and } \quad A/n \defeq \big\{ m \in \N \ \big| \ mn \in A \big\}.\]
When $A \subseteq \Z$, the sets $A-n$ and $A/n$ are defined in the same way with ``$m \in \N$'' replaced by ``$m \in \Z$''.  Any ambiguity caused by the fact that $\N \subseteq \Z$ should be resolved by context or will not matter.

\subsection{Topology, set maps, and systems}
\label{sec_top_set_and_systems}

Throughout this section, let $(X,d_X)$ and $(Y,d_Y)$ be compact metric spaces.  The open ball of radius $\delta > 0$ centered at $x \in X$ is denoted $B(x,\delta)$. A subset of $X$ is \emph{$\delta$-dense} if it has non-empty intersection with every $\delta$-ball in $X$.  Given subsets $W, X' \subseteq X$, we will abuse this terminology slightly by saying that \emph{$W$ is $\delta$-dense in $X'$} to mean that the set $W \cap X'$ is $\delta$-dense as a subset of the metric space $X'$.

Cartesian products $X^d$ of finitely many copies of $X$ will be endowed with the maximum ($L^\infty$) metric, $d_{X^d}$.  Denote by $\Delta( \cdot)$ and $\diag{\ \cdot \ }$ the diagonal injection $X \to X^d$, so that $\Delta(x) = \diag{x} = (x, \ldots, x)$ and $\Delta(X) = \diag{X} = \{ \Delta(x) \ | \ x \in X\}$. The former notation is clearer while the latter is less obtrusive, so we use both depending on the setting.  Though the dimension $d$ is suppressed in the notation, it should always be clear from context.

Any time a partition $X = X_1 \cup \cdots \cup X_C$ is given, we will consider the indices modulo $C$, so that $X_{C+1} = X_1$.  Such a partition is \emph{clopen} if each set $X_i$ is both closed and open.  We denote by $\compc(x)$ the unique element in $\{1, \ldots, C\}$ for which $x \in X_{\compc(x)}$.  Denote by $\ncc(X) \in \N \cup \{\infty\}$ the number of connected components of $X$.  If $\ncc(X) < \infty$, the connected components of $X$ form a clopen partition of $X$.\\

Denote by $\fcx$ the set of non-empty, closed subsets of $X$.  It is a compact metric space when endowed with the \emph{Hausdorff metric}, $d_H$, that measures the distance between two non-empty, closed sets $F, H \in \fcx$ as
\[d_H(F,H) \defeq \max \left ( \max_{f \in F} d_X(f,H), \ \max_{h \in H} d_X(F,h) \right),\]
where $d_X(f,H) = \min_{h \in H} d_X(f,h)$.

\begin{definition}[{cf. \cite[Sec. 18, I.]{kuratowski_1966}}]
A map $\varphi: Y \to \fcx$ is \emph{lower semicontinuous} (\lsc{}) if for all open $U \subseteq X$, the set $\{y \in Y \ | \ \varphi(y) \cap U \neq \emptyset\}$ is open in $Y$.  The map $\varphi$ is \emph{upper semicontinuous} (\usc{}) if for all closed $F \subseteq X$, the set $\{y \in Y \ | \ \varphi(y) \cap F \neq \emptyset\}$ is closed in $Y$.
\end{definition}

It is a fact that a map $\varphi: Y \to \fcx$ is continuous (as a map between metric spaces, where $\fcx$ is equipped with the Hausdorff metric) if and only if it is both \lsc{} and \usc{}; see \cite[A.5 -- A.8]{devries_book_1993}.

We will have need for two different set maps that arise naturally from a continuous map $T: X \to X$.  The first is the map $T: \fcx \to \fcx$, defined at $F \in \fcx$ by $T(F) = \{Tf \ | \ f \in F\}$.  When $T: X \to X$ is a homeomorphism, the map $T: \fcx \to \fcx$ is a homeomorphism \cite[Sec. 17, III., Thm. 2]{kuratowski_1966}.  The second is the map $X \to \fcx$ defined by $x \mapsto T^{-1}\{x\}$.  Since $T$ maps closed sets to closed sets, it follows from \cite[Sec. 18, I., Thm. 4]{kuratowski_1966} that the map $x \mapsto T^{-1}\{x\}$ is continuous if and only if the map $T$ is open (maps open sets to open sets).\\

A \emph{(topological dynamical) system} $(X,T)$ is a compact metric space $(X,d_X)$ together with a continuous map $T: X \to X$. The system is \emph{invertible} if $T$ is a homeomorphism.  Thus, by system (resp. invertible system), we mean an action of the semigroup $(\N,+)$ (resp. group $(\Z,+)$) on a compact metric space by continuous maps.  The \emph{forward and full orbit closures of a point $x \in X$ under $T$} are
\[\orbit{T}{x} \defeq \overline{\big\{ T^t x \ \big| \ t \in \N \big\}} \quad \text{ and } \quad \orbitz{T}{x} \defeq \overline{\big\{ T^t x \ \big| \ t \in \Z \big\}},\]
where the latter makes sense only in an invertible system. The system $(X,T)$ is \emph{minimal} if $\orbit{T}{x} = X$ for all $x \in X$, and it is \emph{totally minimal} if for all $n \in \N$, the system $(X,T^n)$ is minimal.  The \emph{set of return times} of a point $x \in X$ to a set $U \subseteq X$ is
\[R_T(x,U) = \big\{m \in \N \ \big| \ T^m x \in U \big\}.\]
When $(X,T)$ is invertible, we define the return-time set $R(x,U)$ in the same way with ``$m \in \N$'' replaced by ``$m \in \Z$''. It is a fact that the system $(X,T)$ is minimal if and only if for all $x \in X$ and all non-empty, open $U \subseteq X$, the return-time set $R_T(x,U)$ is \emph{additively syndetic}: there exists $N \in \N$ such that $R_T(x,U)$ has non-empty intersection with every set of $N$-many consecutive positive integers.

The following return-time set algebra is useful and follows quickly from the definitions and the set algebra described at the top of this section:
\begin{align}
\label{eqn_return_time_set_algebra}
\begin{gathered}
      R_T(x,U) - n = R_T(T^n x, U) = R_T(x, T^{-n}U), \qquad \frac{R_T(x,U)}{n} = R_{T^n}(x,U),\\
    R_{T^{n_1} \times \cdots \times T^{n_d}} \big((x_1, \ldots, x_d), U_1 \times \cdots \times U_d \big) = \bigcap_{i=1}^d R_{T^{n_i}}(x_i,U_i).
\end{gathered}
\end{align}
When the transformation is clear from the context, we write $R(x,U)$ instead of $R_T(x,U)$.

\begin{lemma}
\label{lemma_tot_min_tot_intersective}
Let $(X,T)$ be a totally minimal, invertible system, $x \in X$, and $U \subseteq X$ be non-empty and open.  The set $R(x,U)$ is \emph{totally intersective}: for all $N \in \N$ and $n \in \Z$, the set $R(x,U) \cap (N \Z + n)$ is non-empty.
\end{lemma}

\begin{proof}
By the return-time set algebra in \eqref{eqn_return_time_set_algebra}, the set $R(x,U) \cap (N \Z + n)$ is non-empty if and only if the set $R_{T^N}(x,T^{-n}U)$ is non-empty.  But $R_{T^N}(x,T^{-n}U)$ is non-empty since $(X,T^N)$ is minimal and the set $T^{-n}U$ is non-empty and open.
\end{proof}

Key to the results in this paper is the study of certain special classes of systems, namely nilsystems and, more broadly, distal systems.  Nilsystems are defined in \cref{sec_nilsystems}.  A system $(X,T)$ is $\emph{distal}$ if for all distinct $x, y \in X$, the quantity $\inf_{t \in \Z} d_X(T^t x, T^t y)$ is non-zero.  Distal systems are invertible, a fact that follows immediately from the fact that the Ellis enveloping semigroup is a group \cite[Theorem 3.1]{furstenberg_1963}.

If $(X,T)$ is an invertible system and another homeomorphism $S: X \to X$ commutes with $T$, then the triple $(X,T,S)$ is a \emph{$\Z^2$-system}, an action of $(\Z^2,+)$ on $X$ by homeomorphisms.  The definitions, terminology, and notation for systems extends to $\Z^2$-systems in the expected way, replacing $T$ with $T, S$ and $T^t$ with $T^t S^s$.  Thus, for example, a $\Z^2$-system is \emph{distal} if for all distinct $x, y \in X$, the quantity $\inf_{t,s \in \Z} d_X(T^tS^s x, T^tS^s y)$ is non-zero.

\begin{definition}
\label{def_clopen_adapted_partitions}
Let $(X,T)$ be a system.  A partition $X = X_1 \cup \cdots \cup X_{C}$ of $X$ is \emph{$T$-adapted} if for all $i \in \Z$, $TX_i = X_{i+1}$. Recall that in any partition, the indices are considered modulo the number of partition elements, so $X_{C+1} = X_1$.
\end{definition}

\begin{lemma}
\label{lem_total_minimality_of_power_v2}
Let $(X,T)$ be a minimal system, and suppose that $X$ has finitely many connected components. There exists a $T$-adapted, clopen partition $X = X_1 \cup \cdots \cup X_{\ncc(X)}$ in which each $X_i$ is connected and each of the systems $(X_i, T^{\ncc(X)})$ is totally minimal.  In particular, the system $(X,T)$ is totally minimal if and only if $X$ is connected.
\end{lemma}

\begin{proof}
Since $T$ is continuous and invertible, it permutes the connected components of $X$.  Thus, the connected components of $X$ can be labeled so as to form a $T$-adapted, clopen partition $X = X_1 \cup \cdots \cup X_{\ncc(X)}$ in which each $X_i$ is connected.  It is not hard to show that a minimal system whose phase space is connected is totally minimal; see, for example, \cite[Thm. 3.1]{ye_1992}. Thus, each of the systems $(X_i, T^{\ncc(X)})$ is totally minimal.  This shows that if $X$ is connected, then $(X,T)$ is totally minimal.  If, on the other hand, the space $X$ is not connected, then $(X,T)$ is not totally minimal since $(X,T^{\ncc(X)})$ is not minimal.
\end{proof}

A set $\Omega \subseteq X$ is \emph{residual} if it contains a dense $G_\delta$ set.  A property holds for \emph{almost every point} or \emph{almost all points} $x \in X$ if it holds for a residual set of points.  A map $f: X \to Y$ is \emph{open} (resp. \emph{semiopen)} if the image of every open set is open (resp. has non-empty interior).

\begin{lemma}
\label{lem_factor_maps_pass_residuality}
Let $\pi: X \to Y$ be a continuous, semiopen surjection of compact metric spaces.
\begin{enumerate}[label=(\Roman*)]
    \item \label{item_factor_map_fact_one} Images of residual sets under $\pi$ are residual.
    \item \label{item_factor_map_fact_two} Preimages of residual sets under $\pi$ are residual.
    \item \label{item_factor_map_fact_four} For all $\delta > 0$, there exists $\delta' > 0$ and a $\delta$-dense set $W \subseteq X$ such that for all $w \in W$, $\pi B(w,\delta) \supseteq B(\pi w,\delta')$.
\end{enumerate}
Factor maps of minimal systems are semiopen and satisfy the properties in \ref{item_factor_map_fact_one}, \ref{item_factor_map_fact_two}, and \ref{item_factor_map_fact_four}.
\end{lemma}

\begin{proof}
Statement \ref{item_factor_map_fact_one} is shown in \cite[Lemma 2.6]{glasscock_koutsogiannis_richter_2019}.  Statement \ref{item_factor_map_fact_two} is an exercise using the fact that a continuous surjection is semiopen if and only if preimages of dense sets are dense.

To see \ref{item_factor_map_fact_four}, let $\delta > 0$.  Call a point $x \in X$ \emph{$\sigma$-good} if $\pi B(x,\delta) \supseteq B(\pi x,\sigma)$.  By the compactness of $X$, it suffices to show that the set $\cup_{\sigma > 0} \{x \in X \ | \ \text{$x$ is $\sigma$-good}\}$ is $\delta$-dense.  Suppose for a contradiction that there exists $z \in X$ such that the set has empty intersection with $B(z,\delta)$.  Since $\pi$ is semiopen, there exists $y \in Y$ and $\sigma > 0$ such that $\pi B(z,\delta/2) \supseteq B(y,\sigma)$.  Thus, there exists $z' \in B(z,\delta / 2) \cap \pi^{-1}(y)$ such that $\pi B(z',\delta) \supseteq \pi B(z,\delta / 2) \supseteq B(y,\sigma)$.  This implies that $z' \in B(z,\delta)$ is $\sigma$-good, a contradiction.

Finally, that factor maps of minimal systems are semiopen is shown in the proof of \cite[Lem. 2.9]{glasscock_koutsogiannis_richter_2019}.
\end{proof}

\subsection{The prolongation relation in distal systems}
\label{sec_prolongation_relation}

Let $(X,T)$ be an invertible system.  Following \cite{auslander_guerin_1997}, for $x \in X$, denote by $\prolong_T(x)$ the \emph{prolongation class of $x$}, defined by
\begin{align}
    \label{def_prolongation}
    \prolong_T(x) \defeq \big\{ y \in X \ \big| \ \exists (x_n)_n \subseteq X, \ \exists (k_n)_n \subseteq \Z, \ \lim_{n \to \infty} x_n = x \text{ and } \lim_{n \to \infty} T^{k_n} x_n = y \big\}.
\end{align}
Note that $\prolong_T(x)$ is closed subset of $X$ and that $\orbitz{T}{x} \subseteq \prolong_T(x)$.

As functions from $X$ into $\fcx$, the map $\orbitz{T}{}: x \mapsto \orbitz{T}{x}$ is \lsc{} \cite[Lem. 2.13]{glasscock_koutsogiannis_richter_2019} and the map $\prolong_T: x \mapsto \prolong_T(x)$ is \usc{} \cite[Section 2]{auslander_guerin_1997}.  As such, for each, the points of continuity form a residual subset of $X$ \cite{fort_1951}.  It is shown in \cite[Thm. 1, Lemmas 2, 3]{akin_glasner_1998} that $\orbitz{T}{x} = \prolong_T(x)$ precisely when $x$ is a point of continuity of the map $\orbitz{T}{}$.

Viewing $\orbitz{T}{}$ and $\prolong_T$ as relations (ie. $(x,y) \in \orbitz{T}{} \subseteq X^2$ if and only if $x \in \orbitz{T}{y}$, and similarly for $\prolong_T$), the prolongation relation is the closure of the orbit relation.  In general, neither relation is an equivalence relation.  If $(X,T)$ is distal, then $\orbitz{T}{}$ is an equivalence relation \cite[Thm. 3.2]{furstenberg_1963}, but the prolongation relation (which is always closed, reflexive, and symmetric) may fail to be transitive.\footnote{Consider the distal transformation of $\T^3$ defined by $T(x,y,z) = (x,y+2x,z)$ when $0 \leq x \leq 1/2$ and $T(x,y,z) = (x,y,z+2x)$ when $1/2 \leq x \leq 1$.  It is quick to check that for all $y, z \in \T$, the prolongation class $\prolong_T((0,y,z))$ is the union of $\{0\} \times \T \times \{z\}$ and $\{0\} \times \{y\} \times \T$.  Thus, $\prolong_T$ is not transitive: the points $(0, 1/2, 0)$ and $(0,0,1/2)$ belong to $\prolong_T((0,0,0))$, but $(0, 1/2, 0) \not\in \prolong_T((0, 0, 1/2))$.}

The prolongation relation is transitive (and, hence, an equivalence relation) in systems that have a sufficient supply of automorphisms; classes of such ``homogeneous'' or ``regular'' systems are main objects of study in \cite{akin_glasner_1998}, \cite{auslander_guerin_1997}, and \cite{auslander_markley_1998}.  Moreover, in such systems, the map $x \mapsto \prolong_T(x)$ is continuous.  We summarize what we need from the discussion above in the following lemma.

\begin{lemma}
\label{lemma_orbit_prolong_continuity_considerations}
Let $(X,T,S)$ be a minimal, distal $\Z^2$-system.
\begin{enumerate}[label=(\Roman*)]
    \item \label{item_orbit_prolong_continuity_considerations_zero} The relations $\orbitz{T}{}$ and $\prolong_T$ are equivalence relations whose equivalence classes are $T$-invariant, $\prolong_T$ is closed, and $\orbitz{T}{} \subseteq \prolong_T$.
    \item \label{item_orbit_prolong_continuity_considerations_four} The maps $T$ and $\prolong_S$ commute, and $\prolong_S$ is $T \times T$-invariant.
    \item \label{item_orbit_prolong_continuity_considerations_two} The map $\prolong_T: X \to \fcx$ is continuous.
    \item \label{item_orbit_prolong_continuity_considerations_five} For all $x, y \in X$, $\prolong_T(x) \cap \prolong_S(y) \neq \emptyset$.
\end{enumerate}
All of the previous statements hold with the maps $T$ and $S$ interchanged. 
\end{lemma}

\begin{proof}
\ref{item_orbit_prolong_continuity_considerations_zero} That $\orbitz{T}{}$ is an equivalence relation follows from \cite[Thm. 3.2]{furstenberg_1963}.  The discussion in \cite[Section 2]{auslander_guerin_1997} gives that $\prolong_T$ is closed, symmetric, and reflexive.  It follows from \cite[Cor. 8]{auslander_markley_1998} (with, in their notation, the group generated by $T$ and $S$ as ``$T$'' and the group generated by $T$ as ``$\Phi$'') that $\prolong_T$ is an equivalence relation.  It follows from the definitions that $\orbitz{T}{} \subseteq \prolong_T$ and that for all $x \in X$, $T \orbitz{T}{x} = \orbitz{T}{x}$ and $T \prolong_T(x) = \prolong_T(x)$.

\ref{item_orbit_prolong_continuity_considerations_four} Because the maps $T$ and $S$ commute, it is quick to check from the definition of the prolongation relation that for all $x \in X$, $T \prolong_S(x) = \prolong_S(Tx)$.  It follows that $x \in \prolong_S(y)$ if and only if $Tx \in \prolong_S(Ty)$, whereby $(T \times T) \prolong_S = \prolong_S$.

\ref{item_orbit_prolong_continuity_considerations_two}
Since $\prolong_T$ is a closed equivalence relation, the quotient space $X / \prolong_T$ is compact (and metrizable).  It follows from general considerations (eg. \cite[Sec. 18, I., Thm. 4]{kuratowski_1966}) that the map $\prolong_T: X \to \fcx$ is continuous if and only if the quotient map $\prolong_T: X \to X / \prolong_T$ is open.  By \ref{item_orbit_prolong_continuity_considerations_four}, the pair $(X / \prolong_T,T,S)$ is a system, and by the same fact, the map $\prolong_T: (X,T,S) \to (X / \prolong_T,T,S)$ is a factor map of minimal, distal systems.  It follows by \cite[Thm. 8.1]{furstenberg_1963} that the map $\prolong_T: X \to X / \prolong_T$ is open, and hence that the map $\prolong_T$ is continuous.

\ref{item_orbit_prolong_continuity_considerations_five} Let $x, y \in X$.  It suffices by the containment in \ref{item_orbit_prolong_continuity_considerations_zero} to show that $\orbitz{T}{x} \cap \prolong_S(y) \neq \emptyset$.  By the minimality of $(X,T,S)$, there exists $(k_n)_n, (\ell_n)_n \subseteq \Z$ such that $\lim_{n \to \infty} S^{\ell_n} T^{k_n} x = y$.  Passing to a subsequence, there exists $z \in X$ such that $\lim_{n \to \infty} T^{k_n} x = z$.  It follows that $y \in \prolong_S(z)$.  Since $\prolong_S$ is reflexive, we have that $z \in \prolong_S(y)$.  But $z \in \orbitz{T}{x}$, too, so $\orbitz{T}{x} \cap \prolong_S(y) \neq \emptyset$.

All of the previous statements hold with $T$ and $S$ interchanged by the symmetry between $T$ and $S$.
\end{proof}

Looking ahead to \cref{sec_prelim_to_main}, the minimal, distal $\Z^2$-system in which we will apply \cref{lemma_orbit_prolong_continuity_considerations} and the results in the following sections will be the orbit closure of a point on the diagonal of $X^{d+1}$ under the maps $T \times \cdots \times T$ and $T^0 \times T^1 \times \cdots \times T^d$.  Such systems are given as examples of ``homogeneous systems'' in \cite[Page 43]{akin_glasner_1998}.

\subsection{Continuity of the intersections of prolongations}

The main result in this section is the continuity of the map $\INT: X^2 \to \fcx$ defined by
\begin{align}
\label{eqn_def_of_int}
    \INT(x,y) \mapsto \prolong_T(x) \cap \prolong_S(y)
\end{align}
in minimal, distal $\Z^2$-systems $(X,T,S)$.  Intersections of continuous, set-valued maps are always \usc{} but not, in general, continuous.  Continuity in our case will follow from the ``regularity'' provided by having commuting automorphisms. Distality is key in a few steps to provide openness via the theorem of Furstenberg that factor maps of minimal, distal systems are open.

\begin{theorem}
\label{thm_continuity_of_intersections}
Let $(X,T,S)$ be a minimal, distal $\Z^2$-system.  The map $\INT: X^2 \to \fcx$ defined in \eqref{eqn_def_of_int} is continuous.
\end{theorem}

The remainder of this section consists of a proof of \cref{thm_continuity_of_intersections}.  We build up to the proof in a number of smaller steps.

\begin{lemma}
\label{lem_map_is_usc}
The map $\INT$ is well-defined (for all $(x,y) \in X^2$, the set $\INT(x,y)$ is non-empty and closed) and \usc{}.
\end{lemma}

\begin{proof}
That the set $\INT(x,y)$ is closed is immediate and that it is non-empty follows from \cref{lemma_orbit_prolong_continuity_considerations} \ref{item_orbit_prolong_continuity_considerations_five}.  That the map $\INT$ is \usc{} follows from \cref{lemma_orbit_prolong_continuity_considerations} \ref{item_orbit_prolong_continuity_considerations_two} and a more general fact \cite[Sec. 18, V., Thm. 1]{kuratowski_1966}: the intersection of two \usc{} maps is \usc{}.
\end{proof}

Denote by $X / \prolong_S$ the space of $\prolong_S$-equivalence classes.  By \cref{lemma_orbit_prolong_continuity_considerations} \ref{item_orbit_prolong_continuity_considerations_zero} and \ref{item_orbit_prolong_continuity_considerations_two}, the space $X / \prolong_S$ can be identified with the image of the map $\prolong_S: X \to \fcx$, and, as such, it is a compact metric space when endowed with the Hausdorff metric.

\begin{lemma}
\label{lem_prolong_map_is_open}
Let $z \in X$.  The map $\prolong_S: \prolong_T(z) \to X/\prolong_S$ is open.  
\end{lemma}

\begin{proof}
We will prove the lemma by recognizing $\prolong_S: \prolong_T(z) \to X/\prolong_S$ as a factor map of distal systems.

Since $T: X \to X$ is a homeomorphism, so is $T: \fcx \to \fcx$.  Thus, by \cref{lemma_orbit_prolong_continuity_considerations} \ref{item_orbit_prolong_continuity_considerations_four}, the triple $(X/\prolong_S,T,S)$ is a system.  The same facts show that $\prolong_S: (X,T,S) \to (X/\prolong_S,T,S)$ is a factor map of systems.  Since $(X,T,S)$ is minimal and distal and $\prolong_S$-equivalence classes are $S$-invariant, the system $(X/\prolong_S,T)$ is a minimal, distal system.

Since $\prolong_T(z)$ is $T$-invariant, the pair $(\prolong_T(z),T)$ is also a system.  We claim that $\prolong_S: (\prolong_T(z),T) \to (X/\prolong_S,T)$ is a factor map of systems.  Continuity of $\prolong_S$ and the intertwining property hold by  \cref{lemma_orbit_prolong_continuity_considerations} \ref{item_orbit_prolong_continuity_considerations_four} and \ref{item_orbit_prolong_continuity_considerations_two}, respectively.  To see that the map is onto, let $F \in X/\prolong_T$.  There exists $y \in X$ such that $F = \prolong_T(y)$.  By \cref{lemma_orbit_prolong_continuity_considerations} \ref{item_orbit_prolong_continuity_considerations_five}, there exists $x \in \prolong_T(z) \cap \prolong_S(y)$.  Thus $x \in \prolong_T(z)$ is such that $\prolong_S(x) = \prolong_S(y) = F$.

Since $\prolong_S: (\prolong_T(z),T) \to (X/\prolong_S,T)$ is a factor map and $(X/\prolong_S,T)$ is a minimal, distal system, it follows from \cite[Thm. 8.1]{furstenberg_1963} that the map $\prolong_S: \prolong_T(z) \to X/\prolong_S$ is open, completing the proof.
\end{proof}

The following lemma establishes continuity of $\INT$ separately in each coordinate.

\begin{lemma}
\label{lem_int_is_continuous_in_one_coord}
Let $z \in X$. The map $\INT_T: X \to \fcx$ defined by $\INT_T(y) = \INT(z,y)$ is continuous.
\end{lemma}

\begin{proof}
The map $\INT_T$ is \usc{} because the map $\INT$ is \usc{} (by \cref{lem_map_is_usc}).  To see that $\INT_T$ is \lsc{}, let $U \subseteq X$ be open, and put
\[V \defeq \big\{ y \in X \ \big| \ \INT_T(y) \cap U \neq \emptyset \big\}.\]
We must show that the set $V$ is open in $X$.  In what follows, it is helpful to refer to \cref{fig_for_one_var_continuity}.

If $V$ is empty, it is open.  Otherwise, let $y \in V$, and let $u \in \INT_T(y) \cap U = \prolong_T(z) \cap \prolong_S(y) \cap U$.  Note that $\prolong_S(u) = \prolong_S(y)$.  Define $U' = \prolong_T(z) \cap U$ so that $U'$ is an open neighborhood of $u$ in $\prolong_T(z)$.

By \cref{lem_prolong_map_is_open}, the map $\prolong_S: \prolong_T(z) \to X/\prolong_S$ is open.  It follows that $\prolong_S(U')$ is an open neighborhood of $\prolong_S(u) = \prolong_S(y)$ in $X/\prolong_S$.  Since $\prolong_S: X \to X/\prolong_S$ is continuous, there exists an open neighborhood $W \subseteq X$ of $y$ such that $\prolong_S(W) \subseteq \prolong_S(U')$.  We will show that $W \subseteq V$; since $W$ is an open neighborhood of $y$ and $y \in V$ is arbitrary, this suffices to show that $V$ is open and finish the proof of the lemma.

Let $y' \in W$. Since $\prolong_S(W) \subseteq \prolong_S(U')$, there exists $u' \in U'$ such that $\prolong_S(y') = \prolong_S(u')$.  It follows that
\[u' \in \prolong_S(u') \cap U' = \prolong_S(u') \cap \prolong_T(z) \cap U = \prolong_T(z) \cap \prolong_S(y') \cap U = \INT_T(y') \cap U.\]
This shows that $\INT_T(y') \cap U \neq \emptyset$, whereby $y' \in V$, as was to be shown.
\end{proof}

\begin{figure}
    \centering
    
    \begin{tikzpicture}

    \path[draw,use Hobby shortcut,closed=true]
    (-3.5/1.3,0/1.3) .. (-1/1.3,-2.5/1.3) .. (2/1.3,-1/1.3) .. (1/1.3,1/1.3);
    \node at (-1.9,0.7) {$U$};
    \node at (-1.7,-.6) {$U'$};
    \filldraw[black] (-0.75,-0.73) circle (2pt) node[anchor=north west]{$u$};
    \filldraw[black] (0.08,-0.63) circle (2pt) node[anchor=north west]{$u'$};
    
    \path[draw,use Hobby shortcut,closed=true]
    (-1.3,-3.9) .. (0,-3.9) .. (1.3,-3.9) .. (0.5,-2.9);
    \node at (-0.85,-3.25) {$W$};
    \filldraw[black] (-0.33,-3.6) circle (2pt) node[anchor=west]{$y$};
    \filldraw[black] (0.68,-3.6) circle (2pt) node[anchor=west]{$y'$};
    
    \draw (-4.2,-1.35) .. controls (0,-.5) .. (5,-.5);
    \draw[ultra thick] (-2.43,-1) .. controls (0,-.6) .. (1.55,-0.52);
    \filldraw[black] (-3.5,-1.22) circle (2pt) node[anchor=north]{$z$};
    \node at (3.3,-0.8) {$D_T(z) = D_T(u)$};
    
    \draw (0,-5) .. controls (-1,-1) and (-1,1) .. (-0.2,3.6);
    \node at (0.87,2.8) {$D_S(y) = D_S(u)$};

    \draw (1,-5) .. controls (0,0) and (-1,1) .. (2.5,2.5);
    \node at (2.15,1.4) {$D_S(y') = D_S(u')$};
    
    \end{tikzpicture}
    
    \caption{The diagram accompanying the proof of \cref{lem_int_is_continuous_in_one_coord}, which shows that when $z$ is fixed, the map $y \mapsto \prolong_T(z) \cap \prolong_S(y)$ is continuous.}
    \label{fig_for_one_var_continuity}
\end{figure}

By interchanging the roles of $T$ and $S$, the following analogue of \cref{lem_int_is_continuous_in_one_coord} holds: for all $z \in X$, the map $\INT_S: X \to \fcx$ defined by $\INT_S(x) = \INT(x,z)$ is continuous.

\begin{lemma}
\label{lem_jint_map_is_open}
Let $z \in X$. Define the maps $\INT_T$ and $\INT_S$ as above, and define $Z = \INT_T(X) \times \INT_S(X) \subseteq (\fcx)^2$.  The map $\JINT: X \to Z$ defined by $\JINT(x) = (\INT_T(x),\INT_S(x))$ is open.
\end{lemma}

\begin{proof}
We will prove the lemma by recognizing $\JINT: X \to Z$ as a factor map of minimal, distal systems.

First, note that the map $T: \fcx \to \fcx$ restricted to $\INT_T(X)$ is a homeomorphism of $\INT_T(X)$.  Indeed, since $T: X \to X$ is a homeomorphism, the map $T: \fcx \to \fcx$ is a homeomorphism, so we need only to show that $T$ maps $\INT_T(X)$ into $\INT_T(X)$.  For this, it suffices to show that the maps $T$ and $\INT_T$ commute.  For $y \in X$, we see by the invertibility of $T$ and \cref{lemma_orbit_prolong_continuity_considerations} \ref{item_orbit_prolong_continuity_considerations_zero} and \ref{item_orbit_prolong_continuity_considerations_four} that
\begin{align*}
    T \INT_T(y) &= T \big( \prolong_T(z) \cap \prolong_S(y) \big) \\
    &= T\prolong_T(z) \cap T\prolong_S(y) \\
    &= \prolong_T(z) \cap \prolong_S(Ty) = \INT_T(Ty).
\end{align*}
Similarly, since the maps $S$ and $\prolong_S$ commute, the map $S: \fcx \to \fcx$ restricts to a homeomorphism of $\INT_S(X)$.

Allowing a slight abuse of notation, we define the maps $T, S: Z \to Z$ by $T(F,H) = (TF, H)$ and $S(F,H) = (F,SH)$.  It follows by the previous paragraph that $T$ and $S$ so defined are commuting homeomorphisms of $Z$, whereby $(Z,T,S)$ is a $\Z^2$-system.

We claim now that $\JINT: (X,T,S) \to (Z,T,S)$ is a factor map of $\Z^2$-systems.  Continuity follows from the continuity of the maps $\INT_T$ and $\INT_S$. The intertwining property follows from the definition of the map $\JINT$ and the definition of the maps $T$ and $S$ on $Z$.  To see that $\JINT$ is onto, let $(F,H) \in Z$.  There exist $x, y \in X$ such that $F = \INT_T(y) = \prolong_T(z) \cap \prolong_S(y)$ and $H = \INT_S(x) = \prolong_T(x) \cap \prolong_S(z)$.  By \cref{lemma_orbit_prolong_continuity_considerations} \ref{item_orbit_prolong_continuity_considerations_five}, there exists $w \in \prolong_T(x) \cap \prolong_S(y)$.  Since $\prolong_S(w) = \prolong_S(y)$ and $\prolong_T(w) = \prolong_T(x)$, we see that
\[\JINT(w) = \big( \INT_T(w), \INT_S(w) \big) = \big( \prolong_T(z) \cap \prolong_S(w), \prolong_T(w) \cap \prolong_S(z) \big) = \big( \INT_T(y), \INT_S(x) \big) = (F,H),\]
as desired.

Since $\JINT: (X,T,S) \to (Z,T,S)$ is a factor map and $(X,T,S)$ is a minimal, distal system, it follows from \cite[Thm. 8.1]{furstenberg_1963} that the map $\JINT: X \to Z$ is open, completing the proof.
\end{proof}

We can finally prove \cref{thm_continuity_of_intersections}.

\begin{proof}[Proof of \cref{thm_continuity_of_intersections}]
The map $\INT$ is \usc{} by \cref{lem_map_is_usc}, so we need only to show that it is \lsc{}. Let $U \subseteq X$ be open, and put
\[V \defeq \big\{ (x,y) \in X^2 \ \big| \ \INT(x,y) \cap U \neq \emptyset \big\}.\]
We must show that the set $V$ is open in $X$. In what follows, it is helpful to refer to \cref{fig_for_two_var_continuity}.

If $V$ is empty, it is open.  Otherwise, let $(x,y) \in V$, and let $z \in \INT(x,y) \cap U$.  Define the map $\JINT$ and the space $Z$ with respect to the point $z$ as in the statement of \cref{lem_jint_map_is_open}.  Since the maps $\INT_T$ and $\INT_S$ are continuous, the map $\KINT: X^2 \to Z$ defined by $\KINT(x',y') = \big( \INT_T(y'), \INT_S(x') \big)$ is continuous. Note that
\begin{align*}
    \JINT(z) &= \big( \INT_T(z), \INT_S(z) \big)\\
    &= \big( \prolong_T(z) \cap \prolong_S(z), \prolong_T(z) \cap \prolong_S(z) \big) \\
    &= \big( \prolong_T(z) \cap \prolong_S(y), \prolong_T(x) \cap \prolong_S(z) \big) \\
    &= \big( \INT_T(y), \INT_S(x) \big) = \KINT(x,y).
\end{align*}

By \cref{lem_jint_map_is_open}, the map $\JINT$ is open, so the set $\JINT(U)$ is an open neighborhood of $\JINT(z) = \KINT(x,y)$ in $Z$.   Since $\KINT$ is continuous, there exists an open neighborhood $W \subseteq X^2$ of $(x,y)$ such that $\KINT(W) \subseteq \JINT(U)$.  We will show that $W \subseteq V$; since $W$ is an open neighborhood of $(x,y)$ and $(x,y) \in V$ is arbitrary, this suffices to show that $V$ is open and finish the proof of the lemma.

Let $(x',y') \in W$.  Since $\KINT(W) \subseteq \JINT(U)$, there exists $z' \in U$ such that $(\INT_T(y'),\INT_S(x')) = \KINT(x',y') = \JINT(z') = (\INT_T(z'),\INT_S(z'))$.  It follows that
\[\prolong_T(z) \cap \prolong_S(y') = \INT_T(y') = \INT_T(z') = \prolong_T(z) \cap \prolong_S(z').\]
Since $\prolong_S$ is an equivalence relation, this shows that $\prolong_S(y') = \prolong_S(z')$.  Similarly, we see that $\prolong_T(x') = \prolong_T(z')$.  Therefore,
\[z' \in \prolong_T(z') \cap \prolong_S(z') \cap U = \prolong_T(x') \cap \prolong_S(y') \cap U = \INT(x',y') \cap U.\]
This shows that $\INT(x',y') \cap U \neq \emptyset$, and hence that $(x',y') \in V$, as was to be shown.
\end{proof}

\begin{figure}
    \centering
    
    \begin{tikzpicture}
    
    \path[draw,use Hobby shortcut,closed=true]
    (-3.5/1.7,0/1.7) .. (-1/1.7,-2.5/1.7) .. (2/1.7,-1/1.7) .. (1/1.7,1/1.7);
    \node at (-1.7,0) {$U$};
    \filldraw[black] (-0.75,-0.7) circle (2pt) node[anchor=north west]{$z$};
    \filldraw[black] (0.06,0.12) circle (2pt) node[anchor=south west]{$z'$};

    \path[draw,use Hobby shortcut,closed=true]
    (-1.3,-3.85) .. (0,-3.85) .. (1.3,-3.85) .. (0.5,-2.85);
    \node at (-0.85,-3.25) {$\pi_2 W$};
    \filldraw[black] (-0.33,-3.6) circle (2pt) node[anchor=west]{$y$};
    \filldraw[black] (0.68,-3.6) circle (2pt) node[anchor=west]{$y'$};

    \draw (-4.2,-1.35) .. controls (0,-.5) .. (4.5,-.5);
    \draw (-4,2) .. controls (-1,0) .. (4.5,0);
    \path[draw,use Hobby shortcut,closed=true]
    (1.95,-.4) .. (3.1,0.6) .. (3.7,-0.4) .. (2.95,-1.4);
    \node at (2.8,-1.1) {$\pi_1 W$};
    \filldraw[black] (3,-0.5) circle (2pt) node[anchor=north]{$x$};
    \filldraw[black] (3,0) circle (2pt) node[anchor=south]{$x'$};
    \node at (-2.5,-1.5) {$D_T(x) = D_T(z)$};
    \node at (-2.35,2) {$D_T(x') = D_T(z')$};

    \draw (0,-5) .. controls (-1,-1) and (-1,1) .. (-0.2,3.6);
    \node at (0.87,2.8) {$D_S(y) = D_S(z)$};

    \draw (1,-5) .. controls (0,0) and (-1,1) .. (2.5,2.5);
    \node at (2.15,1.4) {$D_S(y') = D_S(z')$};
    
    \end{tikzpicture}
    
    \caption{The diagram accompanying the proof of \cref{thm_continuity_of_intersections}, which shows that the map $(x,y) \mapsto \prolong_T(x) \cap \prolong_S(y)$ is continuous.}
    \label{fig_for_two_var_continuity}
\end{figure}

\begin{remark}
Nowhere in the proof of \cref{thm_continuity_of_intersections} did we use the definition of the prolongation relation.  Thus, the theorem holds in a greater generality for equivalence relations satisfying some or all of the properties laid out in \cref{lemma_orbit_prolong_continuity_considerations}.  We will not have need for this observation, so any useful generalizations along these lines are left to the interested reader.
\end{remark}

\subsection{Uniform denseness of prolongations in prolongations}
\label{sec_uniform_denseness_in_prolongation_classes}

The main result in this section, \cref{thm_uniform_denseness_of_msble_sets}, is a corollary of \cref{thm_continuity_of_intersections} that will be useful later on.

Let $X$ be a compact metric space and $F \in \fcx$.  Recall that $\ncc(F)$ denotes the number of connected components of $F$.  Denote by $\dcc(F) \in [0,\infty]$ the infimum of $d_X(x,y)$ as $x$ and $y$ range over distinct connected components of $F$; if $F$ is connected, put $\dcc(F) = \infty$.

\begin{lemma}
\label{lem_num_and_dist_conn_comp}
Let $(F_n)_n \subseteq \fcx$ and $F \in \fcx$ be such that $\lim_{n \to \infty} F_n = F$, and suppose that $\ncc(F) < \infty$.
\begin{enumerate}[label=(\Roman*)]
    \item \label{item_num_conn_comp_in_limit} $\ncc(F) \leq \liminf_{n \to \infty} \ncc(F_n)$.
    \item \label{item_dist_conn_comp_in_limit} If $\ncc(F_1) = \ncc(F_2) = \cdots = \ncc(F)$, then $\lim_{n \to \infty} \dcc(F_n) = \dcc(F)$.
\end{enumerate}
\end{lemma}

\begin{proof}
\ref{item_num_conn_comp_in_limit} Define $C = \ncc(F)$.  There exist disjoint open sets $U_1, \ldots, U_C \subseteq X$ such that $F = \cup_{i=1}^C (F \cap U_i)$.  It follows by the definition of the Hausdorff metric that for all sufficiently large $n$, the set $F_n \subseteq \cup_{i=1}^C U_i$, so that $\ncc(F_n) \geq C$.

\ref{item_dist_conn_comp_in_limit} Define $C = \ncc(F)$.  If $C = 1$, the conclusion holds immediately since $\dcc(F_n) = \dcc(F) = \infty$.  Suppose $C > 1$ so that $\delta \defeq \dcc(F) \in (0, \infty)$.  Let $0 < \eps <  \delta / 2$.  There exist disjoint, open sets $U_1, \ldots, U_C \subseteq X$ such that $F = \cup_{i=1}^C (F \cap U_i)$ and such that the distance between any two distinct sets is at least $\delta - \eps$.  It follows by the definition of the Hausdorff metric that for all sufficiently large $n$, the set $F_n \subseteq \cup_{i=1}^C U_i$.  Since $F_n$ has $C$-many connected components, each $F_n \cap U_i$ is connected, so $\dcc(F_n) > \delta - \eps$.  This shows that $\liminf_{n \to \infty} \dcc(F_n) \geq \delta$.

Since $\dcc(F) = \delta$, there exists $f_i \in F \cup U_i$ and $f_j \in F \cap U_j$, $i \neq j$, such that $\delta \leq d(f_i, f_j) < \delta + \eps/3$.  For all sufficiently large $n$, the set $F_n \subseteq \cup_{i=1}^C U_i$ and there exist $f_i', f_j' \in F_n$ such that $d(f_i,f_i') < \eps/3$ and $d(f_j,f_j') < \eps/3$. Since $\eps < \delta / 2$, $f_i' \in F_n \cap U_i$ and $f_j' \in F_n \cap U_j$, showing that $\dcc(F_n) \leq \delta + \eps$ and finishing the proof.
\end{proof}
 
\begin{lemma}
\label{lemma_constant_number_of_ccs}
Let $(X,T,S)$ be a minimal, distal $\Z^2$-system, and suppose that there exists $x_0 \in X$ for which $\ncc(\prolong_S(x_0)) < \infty$. There exists $\ccd \in \N$ and $\delta_0 > 0$ such that for all $x \in X$,
\[\ncc(\prolong_S(x)) = \ccd \quad \text{ and } \quad \dcc(\prolong_S(x)) > \delta_0.\]
\end{lemma}

\begin{proof}
For $N \in \N$, define $C_N = \{x \in X \ | \ \ncc(\prolong_S(x)) \leq N\}$.  It follows from \cref{lem_num_and_dist_conn_comp} \ref{item_num_conn_comp_in_limit} that the set $C_N$ is closed. Moreover, it follows from \cref{lemma_orbit_prolong_continuity_considerations}
\ref{item_orbit_prolong_continuity_considerations_zero} and \ref{item_orbit_prolong_continuity_considerations_four} and the fact that homeomorphisms preserve connectedness that the set $C_N$ is $T$- and $S$-invariant.  Therefore, by the minimality of $(X,T,S)$, every set $C_N$ is either empty or all of $X$.  By assumption, some set $C_N$ is non-empty.  Since $C_1 \subseteq C_2 \subseteq \cdots$, we have that there exists $\ccd \in \N$ such that all $S$-prolongation classes have $\ccd$-many connected components.

If $\ccd = 1$, then $\dcc(\prolong_S(x)) = \infty$ for all $x \in X$, so any positive value of $\delta_0$ will suffice.  Suppose $\ccd > 1$ so that $\dcc(\prolong_S(x)) \in (0,\infty)$ for all $x \in X$.  It follows from the previous paragraph and \cref{lem_num_and_dist_conn_comp} \ref{item_dist_conn_comp_in_limit} that the map $\dcc \circ \prolong_S: X \to (0,\infty)$ is continuous.  Since $X$ is compact, there exists $\delta_0 > 0$ such that for all $x \in X$, $\dcc(\prolong_S(x)) > \delta_0$, as desired.
\end{proof}

For the following theorem, recall the definition of ``measurable'' from the beginning of \cref{sec_prolong} and the definition of the map $\INT: X^2 \to \fcx$ as $\INT(x,y) = \prolong_T(x) \cap \prolong_S(y)$ from \eqref{eqn_def_of_int}.

\begin{theorem}
\label{thm_uniform_denseness_of_msble_sets}
Let $(X,T,S)$ be a minimal, distal $\Z^2$-system, and suppose that there exists $x \in X$ for which $\ncc(\prolong_S(x)) < \infty$. For all $\delta > 0$, there exists $\eps > 0$ such that for all $y \in X$, every $\eps$-dense, $\prolong_T$-measurable subset of $X$ is $\delta$-dense in every connected component of $\prolong_S(y)$.
\end{theorem}

\begin{proof}
Let $\delta > 0$.  Let $\delta_0 > 0$ be from \cref{lemma_constant_number_of_ccs} so that for all $y \in X$, the distance between any two points in distinct connected components of $\prolong_S(y)$ is greater than $\delta_0$.  Note that if a set $Z \subseteq X$ is $\min(\delta, \delta_0)$-dense in $\prolong_S(y)$, then it is $\delta$-dense in every connected component of $\prolong_S(y)$.  Therefore, by replacing $\delta$ with $\min(\delta,\delta_0)$, it suffices to show that there exists $\eps > 0$ such that for all $y \in X$, every $\eps$-dense, $\prolong_T$-measurable subset of $X$ is $\delta$-dense in $\prolong_S(y)$.

First we will show that there exists $\eps > 0$ such that for all $x, y \in X$ with $d_X(x,y) < \eps$,
\begin{align}
    \label{eqn_close_intersections}
    \INT(x,y) \cap B(x,\delta) \cap B(y,\delta) \neq \emptyset.
\end{align}
Indeed, by \cref{thm_continuity_of_intersections}, the map $\INT$ is continuous.  Since $X$ is compact, the map $\INT$ is uniformly continuous, so there exists $\eps < \delta / 2$ such that for all $(x,y), (x',y') \in X^2$ with $d_{X^2}\big((x,y),(x',y')\big) < \eps$,
\[d_H \big( \INT(x,y), \INT(x',y') \big) < \delta / 2.\]
If $x, y \in X$ satisfy $d_X(x,y) < \eps$, then $d_{X^2}\big((x,y),(x,x)\big) < \eps$, so $d_H \big( \INT(x,y), \INT(x,x) \big) < \delta / 2$. Combining this with the definition of the Hausdorff metric and the fact that $x \in \INT(x,x)$, we see that the set $\INT(x,y) \cap B(x, \delta / 2)$ is non-empty.  Now \eqref{eqn_close_intersections} follows from the fact that $B(x,\delta / 2) \subseteq B(x,\delta) \cap B(y, \delta)$.

We claim that the $\eps$ from the previous paragraph suffices for the conclusion of the theorem. Let $y \in X$, and let $X' \subseteq X$ be an $\eps$-dense, $\prolong_T$-measurable subset of $X$.  To show that $X'$ is $\delta$-dense in $\prolong_S(y)$, we will show that for all $y' \in \prolong_S(y)$, the set $X' \cap \prolong_S(y') \cap B(y',\delta) \neq \emptyset$.

Let $y' \in \prolong_S(y)$.  Since $X'$ is $\eps$-dense, there exists $x' \in X' \cap B(y',\eps)$.  It follows from \eqref{eqn_close_intersections} that $\prolong_T(x') \cap \prolong_S(y') \cap B(x',\delta) \cap B(y',\delta) \neq \emptyset$.  Note that $\prolong_T(x') \subseteq X'$ since $x' \in X'$ and $X'$ is $\prolong_T$-measurable.  Therefore, the set $X' \cap \prolong_S(y') \cap B(y',\delta)$ is non-empty, as was to be shown.
\end{proof}

\section{Rationality in nilpotent groups}
\label{sec_rational_elts_in_groups}

The main results in this section concern the relationship between rationality and periodicity in nilpotent groups and nilmanifolds.  For example, we show in \cref{sec_periodicity_of_rat_polys} that $\Gamma$-rational-valued polynomials in a nilpotent group $G$ are periodic modulo $\Gamma$.
Much has already been written in this direction, notably by Leibman \cite{leibman_rational_2006} and Green and Tao \cite[Appendix A]{green_tao_quantitative_behaviour_2012}.
While we have need for a quantitative approach similar to Green and Tao's, the formulations we will need are not stated explicitly in their work.  In developing the material for our needs, we aim to keep the setting as general as possible by relying on the Lie structure as little as possible. Some recent examples of nilpotent groups without Lie structure appearing in the subject \cite{jamneshan_shalom_tao_2021} demonstrate that there may be future value in keeping a general setting.

Let $G$ be a group and $\Gamma \leq G$ be a subgroup.  An element $g \in G$ is \emph{$Q$-rational (with respect to $\Gamma$)}, $Q \in \N$, if there exists $n \in \{1, \ldots, Q\}$ such that $g^n \in \Gamma$; it is \emph{rational (with respect to $\Gamma$)} if it is $Q$-rational for some $Q$.  The set of $Q$-rational and rational elements of $G$ are denoted by $\Q_Q(G,\Gamma)$ and $\Q(G,\Gamma)$, respectively.

For the product group $G^d$ and the subgroup $\Gamma^d \leq G^d$, it is quick to check that
\[\Q_Q(G^d,\Gamma^d) \subseteq \Q_{Q}(G,\Gamma)^d \subseteq \Q_{Q^d}(G^d,\Gamma^d),\]
whereby $\Q(G^d,\Gamma^d) = \Q(G,\Gamma)^d$.

\subsection{Rational elements in nilpotent groups}

Let $G$ be a group.  The \emph{$i^{\text{th}}$-iterated commutator subgroup of $G$}, denoted $G_i$, is defined inductively by $G_1 = G$ and $G_{i+1} = [G_i,G]$, where the commutator bracket is defined by $[g,h] = ghg^{-1}h^{-1}$. The sequence $(G_i)_i$ is a decreasing sequence of normal subgroups of $G$ called the \emph{lower central series of $G$}.  If $G_{d+1} = \{e_G\}$, $d \in \Nz$, then the group $G$ is \emph{$d$-step nilpotent}.

If $G$ is nilpotent, then the rational elements with respect to a subgroup $\Gamma$ form a group.  This was proved by Leibman \cite[Lemma 1.3]{leibman_rational_2006}.  \cref{lemma_rationality_of_subgroup} below is a quantitative version of that result, proved in the same way but keeping track of the height of the rationals involved.

\begin{lemma}
\label{lem_useful_facts_for_leibman_lemma}
Let $G$ be a nilpotent group with lower central series $(G_i)_i$.
\begin{enumerate}[label=(\Roman*)]
    \item \label{itm_nil_fact_one} $[G_j,G_k] \subseteq G_{j+k}$.
    \item \label{itm_nil_fact_three} For all $g_1, \ldots, g_\ell \in G$, $(g_1 \cdots g_\ell)^n G_2 = g_1^n \cdots g_\ell^n G_2$.
    \item \label{itm_nil_fact_two} For all $g_1 \in G_j$ and $g_2 \in G_k$, if $j+k \geq i+1$, then $[g_1^n, g_2^m] G_{i+2} = [g_1,g_2]^{nm} G_{i+2}$.
    \item \label{itm_nil_fact_four} For all subgroups $H \leq G$, $\big[G_2 H, G_{i+1} (H \cap G_i)\big] \subseteq G_{i+2}H$.
\end{enumerate}
\end{lemma}

\begin{proof}
\ref{itm_nil_fact_one} This is a standard fact; see, eg., \cite[Ch. 2, Lemma 3]{host_kra_book_2018}.

\ref{itm_nil_fact_three} and \ref{itm_nil_fact_two} These follow from \ref{itm_nil_fact_one} and simple inductions.  For \ref{itm_nil_fact_two}, note, for example, that $[g_1,g_2]^2 = [g_1,g_2]g_1g_2g_1^{-1}g_2^{-1} = \big[ [g_1,g_2], g_1\big] [g_1^2,g_2] \in [g_1^2,g_2] G_{2j + k} \subseteq [g_1^2,g_2] G_{i+2}$.

\ref{itm_nil_fact_four}  Since $G_2$ is a normal subgroup of $G$, we have $G_2 H = \{h g \ | \ h \in H, \ g \in G_2\}$, and similarly for $G_{i+1} (H \cap G_i)$.  Therefore, the group $\big[G_2 H, G_{i+1} (H \cap G_i)\big]$ is generated by elements of the form $[h_1 g_1, h_2 g_2] = h_1 g_1 h_2 g_2 g_1^{-1} h_1^{-1} g_2^{-1} h_2^{-1}$ with $g_1 \in G_2$, $g_2 \in G_{i+1}$, $h_1 \in H$, and $h_2 \in H \cap G_{i}$.  Since $g_1 \in G_2$ and $h_2 \in G_i$, $g_1 h_2 \in h_2 g_1 G_{i+2}$.  Since $h_1^{-1} \in G_1$ and $g_2^{-1} \in G_{i+1}$, $h_1^{-1} g_2^{-1} \in g_2^{-1} h_1^{-1} G_{i+2}$.  Therefore, $[h_1 g_1, h_2 g_2] \in h_1 h_2 [g_1, g_2] h_1^{-1} h_2^{-1} G_{i+2} = [h_1,h_2]G_{i+2} \subseteq G_{i+2}H$.
\end{proof}

\begin{lemma}
\label{lemma_rationality_of_subgroup}
For all $Q \in \N$ and $d \in \Nz$, there exists $Q' \in \N$ such that the following holds.  Let $G$ be a $d$-step nilpotent group, and let $\Gamma \leq G$ be a subgroup. The subgroup of $G$ generated by $\Q_Q(G,\Gamma)$ is contained in $\Q_{Q'}(G,\Gamma)$.
\end{lemma}

\begin{proof}
Let $Q \in \N$ and $d \in \Nz$.  Put $Q' = (Q!)^{d (d+1) / 2}$.  Let $G$ be a $d$-step nilpotent group, and let $\Gamma \leq G$.  Let $H$ be the subgroup generated by $\Q_Q(G,\Gamma)$.  Define $H^{\wedge m} \defeq \{h^m \ | \ h \in H\}$. To show that $H \subseteq \Q_{Q'}(G,\Gamma)$, it suffices to show that $H^{\wedge Q'} \subseteq \Gamma$.

Let $(H_i)_i$ be the lower central series of $H$; since $G$ is $d$-step nilpotent, so is $H$, so $H_{d+1} = \{e_G\}$. We will show by induction on $i$ that for all $i \in \N$,
\begin{align}
\label{eqn_first_stop_in_leibman_lemma}
    H_i^{\wedge Q!^i} \subseteq H_{i+1} (H_i \cap \Gamma).
\end{align}

Let $i = 1$, and let $h \in H_1$.  By the definition of $H_1 = H$, there exist $r_1, \ldots, r_k \in \Q_Q(G,\Gamma)$ and $e_1, \ldots, e_k \in \Z$ such that $h = r_1^{e_1} \cdots r_k^{e_k}$.  By \cref{lem_useful_facts_for_leibman_lemma} \ref{itm_nil_fact_three}, $h^{Q'} H_2 = r_1^{Q'e_1} \cdots r_k^{Q'e_k} H_2 \subseteq (\Gamma \cap H) H_2$.   Therefore, $H_1^{\wedge Q!} \subseteq H_{2} (H_1 \cap \Gamma)$, verifying the base case.

Suppose \eqref{eqn_first_stop_in_leibman_lemma} holds for some $i \in \N$.  Using $H_{i+1} = [H_i,H] = [H,H_i]$, \cref{lem_useful_facts_for_leibman_lemma} \ref{itm_nil_fact_two} and \ref{itm_nil_fact_four} (with $H$ as $G$ and $H \cap \Gamma$ as $H$), and the induction hypothesis,
\begin{align*}
    H_{i+1}^{\wedge Q!^{i+1}} \subseteq [H_1, H_i]^{\wedge Q!^{i+1}} H_{i+2} &\subseteq [H_1^{\wedge Q!},H_i^{\wedge Q!^i}] H_{i+2} \\
    &\subseteq [H_{2} (H_1 \cap \Gamma), H_{i+1} (H_i \cap \Gamma)] H_{i+2} \subseteq H_{i+2} (H_1 \cap \Gamma).
\end{align*}
Thus, $H_{i+1}^{\wedge Q!^{i+1}} \subseteq H_{i+2} (H_1 \cap \Gamma)$.
Intersecting both sides with $H_{i+1}$, we see that $H_{i+1}^{\wedge Q!^{i+1}} \subseteq H_{i+2} (H_{i+1} \cap \Gamma)$, as desired.

Now we will show by induction on $i$ that for all $i \in \N$,
\begin{align}
\label{eqn_second_stop_in_leibman_lemma}
    H^{\wedge Q!^{i(i+1)/2}} \subseteq H_{i+1}(H \cap \Gamma);
\end{align}
setting $i=d$ gives the desired result and finishes the proof of the lemma.  For $i=1$, we need that $H^{\wedge Q!} \subseteq H_{2}(H \cap \Gamma)$, which has already been shown in \eqref{eqn_first_stop_in_leibman_lemma}.

Suppose \eqref{eqn_second_stop_in_leibman_lemma} holds for some $i \in \N$.  Using the induction hypothesis, \cref{lem_useful_facts_for_leibman_lemma} \ref{itm_nil_fact_three}, and \eqref{eqn_first_stop_in_leibman_lemma},
\begin{align*}
    H^{\wedge Q!^{(i+1)(i+2)/2}} = \big(H^{\wedge Q!^{i(i+1)/2}} \big)^{\wedge Q!^{i+1}} &\subseteq \big(H_{i+1}(H \cap \Gamma) \big)^{\wedge Q!^{i+1}} \\
    &\subseteq H_{i+1}^{\wedge Q!^{i+1}}(H \cap \Gamma)^{\wedge Q!^{i+1}} H_{i+2} \subseteq H_{i+2} (H_1 \cap \Gamma).
\end{align*}
This completes the proof of \eqref{eqn_second_stop_in_leibman_lemma} and the proof of the lemma.
\end{proof}

\subsection{Periodicity of rational polynomials modulo a subgroup}
\label{sec_periodicity_of_rat_polys}

Let $G$ be a group. A \emph{filtration of degree $d$ on $G$}, $d \in \Nz$, is a sequence $\Gdot = (G^{(i)})_{i=0}^{d+1}$ of subgroups of $G$ such that
\[G = G^{(0)} \supseteq G^{(1)} \supseteq \cdots \supseteq G^{(d)} \supseteq G^{(d+1)} = \{e_G\}\]
and such that for all $0 \leq i, j \leq i+j \leq d+1$, $[G^{i},G^{j}] \subseteq G^{(i+j)}$.  A \emph{filtration on $G$} is a filtration of some degree $d$.  If $G$ is a $d$-step nilpotent group, setting $G^{(0)} = G$ and $G^{(i)} = G_i$ yields the lower central series filtration of $G$ (cf. \cref{lem_useful_facts_for_leibman_lemma} \ref{itm_nil_fact_one}), a filtration of degree $d$.

A filtration $\Gdot$ on $G$ passes naturally to its subgroups and quotient groups (cf. \cite[Ch. 14, Sec. 1.3]{host_kra_book_2018}).
\begin{itemize}
    \item The \emph{induced filtration} on $H \leq G$, defined by $(G^{(i)} \cap H)_{i=0}^{d+1}$, is a filtration of degree $d$.
    \item If $K \leq G$ is normal, the \emph{quotient filtration} on $G/K$, defined by $(K G^{(i)} / K)_{i=0}^{d+1}$ is a filtration of degree $d$.  If $K = G^{(d)}$, then the quotient filtration is of degree $d-1$ since $K G^{(d)} / K = \{e_{G/K}\}$.
\end{itemize}

Denote by $G^\Z$ the space of bi-infinite, $G$-valued sequences (functions from $\Z$ to $G$). For $P \in G^\Z$ and $N \in \Z$, define $\partial_N P, \partial'_N P \in G^\Z$ by
\begin{align}
    \notag
    (\partial_NP)(n) &= P(n+N)P(n)^{-1},\\
    \label{eqn_def_of_reverse_diff}
    (\partial'_NP)(n) &= P(n+N)^{-1}P(n).
\end{align}
We denote the $i$-fold composition $\partial_N \circ \cdots \circ \partial_N$ by $\partial_N^{\circ i}$ and similarly for $\partial'_N$.

\begin{definition}
Let $\Gdot = (G^{(i)})_{i=0}^{d+1}$ be a degree-$d$ filtration on $G$.  A sequence $P \in G^\Z$ is a \emph{polynomial sequence with coefficients in $\Gdot$} if for all $i \in \{0, \ldots, d+1\}$, the sequence $\partial_1^{\circ i}P$ is $G^{(i)}$-valued. The set of polynomials with coefficients in $\Gdot$ is denoted $\poly(\Gdot)$.
\end{definition}

We will make frequent use of the following facts about polynomials with coefficients in $\Gdot$, a filtration on $G$ of degree $d \geq 1$.
\begin{enumerate}[label=(\Roman*)]
    \item \label{item_poly_fact_one} \cite[Ch. 14, Thms 8, 10]{host_kra_book_2018} A sequence $P \in G^\Z$ is a member of $\poly(\Gdot)$ if and only if there exist $a_0, \ldots, a_d \in G$, with $a_i \in G^{(i)}$, such that for all $n \in \Z$,
    \[P(n) = a_0^{\binom{n}{0}} a_1^{\binom{n}{1}} a_2^{\binom{n}{2}} \cdots a_d^{\binom{n}{d}}.\]
    The elements $a_0, \ldots, a_d$ are called \emph{the coefficients of $P$}.  The coefficients of $P$, and hence all values of $P$, are determined by any $(d+1)$-many consecutive values of $P$.
    
    \item \label{item_poly_fact_two} \cite[Ch. 14, Thm. 5 \& Prop. 7]{host_kra_book_2018} The set $\poly(\Gdot)$, endowed with pointwise multiplication, is a group that is closed under translation $P \mapsto (n \mapsto P(n+1))$ and inversion $P \mapsto (n \mapsto P(n)^{-1})$.
    
    \item \label{item_poly_fact_three} \cite[Ch. 14, Prop. 9]{host_kra_book_2018} Let $H \leq G$, and let $\Hdot$ be the induced filtration.  The group $\poly(\Hdot)$ is a subgroup of $\poly(\Gdot)$.  If $(d+1)$-many consecutive values of $P$ belong to $H$, then $P \in \poly(\Hdot)$.  In particular, a polynomial $P \in \poly(\Gdot)$ belongs to $\poly(\Hdot)$ if and only if it is $H$-valued; that is, $\poly(\Hdot) = \poly(\Gdot) \cap H^\Z$.
    
    \item \label{item_poly_fact_four} Let $K$ be a normal subgroup of $G$ and $\Kdot$ be the induced filtration. Define $H = G / K$, and let $\Hdot$ be the quotient filtration.  Applying the quotient map $\pi: G \to H$ coordinate-wise and identifying $H^\Z$ with $G^\Z / K^\Z$ yields the quotient map $\pi: G^\Z \to G^\Z / K^\Z$.  It is quick to check using fact \ref{item_poly_fact_one} and the definition of the quotient filtration that $\pi \poly(\Gdot) = \poly(\Hdot)$.  By fact \ref{item_poly_fact_three}, $\poly(\Kdot) = \poly(\Gdot) \cap K^\Z$, so $\poly(\Kdot)$ is a normal subgroup of $\poly(\Gdot)$ and, restricting $\pi$ to $\poly(\Gdot)$, we get the quotient map $\pi: \poly(\Gdot) \to \poly(\Gdot) / \poly(\Kdot)$.  Therefore, $\pi \poly(\Gdot) = \poly(\Hdot) = \poly(\Gdot) / \poly(\Kdot)$.
\end{enumerate}

We will be particularly interested in rational-valued polynomials.

\begin{definition}
Let $\Gamma \leq G$.  A polynomial $P \in \poly(\Gdot)$ is \emph{$Q$-rational (with respect to $\Gamma$)}, $Q \in \N$, if it is $\Q_Q(G,\Gamma)$-valued: for all $n \in \Z$, $P(n)$ is $Q$-rational (with respect to $\Gamma$).  The polynomial $P$ is \emph{rational (with respect to $\Gamma$)} if it is $Q$-rational for some $Q$.
\end{definition}

The following lemma gives a quantitative version of the fact that a polynomial is rational if and only if its coefficients are rational.

\begin{lemma}
\label{lem_coeffs_and_values}
For all $Q \in \N$ and $d \in \Nz$ there exists $Q' \in \N$ for which the following holds.  Let $\Gdot$ be degree-$d$ filtration on $G$, $\Gamma \leq G$, and $P \in \poly(\Gdot)$.
\begin{enumerate}[label=(\Roman*)]
    \item \label{item_rat_poly_conclusion_one} If $P$ is $Q$-rational, then the coefficients of $P$ are $Q'$-rational.
    \item \label{item_rat_poly_conclusion_two} If the coefficients of $P$ are $Q$-rational, then $P$ is $Q'$-rational.
    \item \label{item_rat_poly_conclusion_three} If $(d+1)$-many consecutive values of $P$ are $Q$-rational, then $P$ is $Q'$-rational.
\end{enumerate}
\end{lemma}

\begin{proof}
Let $Q \in \N$ and $d \in \Nz$, and put $Q' = Q'(Q,d)$ from \cref{lemma_rationality_of_subgroup}. \ref{item_rat_poly_conclusion_one} holds since the coefficients of $P$ are contained in the subgroup of $G$ generated by the values of $P$.  \ref{item_rat_poly_conclusion_two} holds since the values of $P$ are contained in the subgroup generated by the coefficients of $P$.  To see \ref{item_rat_poly_conclusion_three}, let $H$ be the subgroup of $G$ generated by the $Q$-rational elements.  Since $(d+1)$-many consecutive values of $P$ are $H$-valued, polynomial fact \ref{item_poly_fact_three} above gives that $P$ is $H$-valued. Since all elements of $H$ are $Q'$-rational, the polynomial $P$ is $Q'$-rational.
\end{proof}

\begin{lemma}
\label{lem_periodic_polys_modulo_subgroup_for_abelian_groups}
For all $Q \in \N$ and $d \in \Nz$, there exists $N \in \N$ such that the following holds.  If $G$ is an abelian group with a degree-$d$ filtration $\Gdot$, $\Gamma \leq G$, and $P \in \poly(\Gdot)$ is $Q$-rational with respect to $\Gamma$, then
\[\partial'_NP \in \poly(\Gammadot),\]
where $\Gammadot$ denotes the induced filtration on $\Gamma$. 
\end{lemma}

\begin{proof}
Let $Q \in \N$ and $d \in \Nz$.  Suppose $d = 0$, and define $N = 1$.  Let $G$ be a group with a filtration $\Gdot = (G^{(i)})_{i=0}^1$ of degree $0$; thus $G^{(0)} = G$ and $G^{(1)} = \{e_G\}$.  Let $\Gamma \leq G$ and $P \in \poly(\Gdot)$ be $Q$-rational with respect to $\Gamma$.  Since $\partial_1 P$ is $G^{(1)}$-valued, the sequence $P$ is constant.  Thus $\partial'_N P$ is the constant $e_G$ sequence, whereby $\partial'_NP \in \poly(\Gammadot)$.

Suppose $d \geq 1$. Let $Q'$ be from \cref{lem_coeffs_and_values}, and define $N = d! Q'!$.  Let $G$, $\Gdot$, $\Gamma$, and $P \in \poly(\Gdot)$ be as in the statement of the lemma.  To conclude that $\partial'_NP \in \poly(\Gammadot)$, it suffices by polynomial fact \ref{item_poly_fact_three} above to show that $\partial'_NP$ is $\Gamma$-valued. 

Let $a_0, \ldots, a_d \in G$ be the coefficients of $P$ (cf. polynomial fact \ref{item_poly_fact_one} above).
Since $G$ is abelian, for all $n \in \Z$,
\[(\partial'_NP)(n) = P(n+N)^{-1}P(n) = a_1^{\binom{n}{1} - \binom{n+N}{1}} a_2^{\binom{n}{2} - \binom{n+N}{2}} \cdots a_d^{\binom{n}{d} - \binom{n+N}{d}}.\]
For each $i \in \{1, \ldots, d\}$, note that $Q'!$ divides $N / i!$, which divides $\binom{n}{i} - \binom{n+N}{i}$.  By \cref{lem_coeffs_and_values}, each coefficient $a_i$ is $Q'$-rational, so $a_i^{\binom{n+N}{i} - \binom{n}{i}} \in \Gamma$.  It follows that $(\partial'_NP)(n) \in \Gamma$ for all $n \in \Z$, as was to be shown.
\end{proof}

We move now to prove the same result without the assumption that $G$ is abelian. We will have need for the following useful telescoping identity for $P \in G^\Z$ and $a, b, n \in \Z$:
\begin{align}
    \label{eqn_useful_identity_for_derivatives}
    \partial'_{ab}P(n) = \partial'_a P \big(a(b-1) + n \big) \  \partial'_a P \big(a(b-2) + n \big) \ \cdots \ \partial'_a P \big(a + n \big) \ \partial'_a P \big(n \big).
\end{align}

\begin{theorem}
\label{thm_periodic_polys}
    For all $Q \in \N$ and $d \in \Nz$, there exists $N \in \N$ such that the following holds.  If $G$ is a group with a degree-$d$ filtration $\Gdot$, $\Gamma \leq G$, and $P \in \poly(\Gdot)$ is $Q$-rational with respect to $\Gamma$, then
    \[\partial'_NP \in \poly(\Gammadot),\]
    where $\Gammadot$ denotes the induced filtration on $\Gamma$. 
\end{theorem}

\begin{proof}
The proof proceeds by induction on $d \in \Nz$. If $d=0$, set $N=1$, and use the same reasoning at the beginning of \cref{lem_periodic_polys_modulo_subgroup_for_abelian_groups} to see that $\partial'_NP \in \poly(\Gammadot)$.

Suppose $d \geq 1$ and that the theorem holds for $d-1$. Let $Q \in \N$.  Let $N_0$ be from the induction hypothesis for $d-1$ and $Q$.  Let $Q'$ be from \cref{lemma_rationality_of_subgroup} for $d$ and $Q$. Let $N_1$ be from \cref{lem_periodic_polys_modulo_subgroup_for_abelian_groups} for $Q'$ and $d$.  Finally, define $N = N_0N_1 Q'!$.

Let $G$ be a group with a degree-$d$ filtration $\Gdot$, $\Gamma \leq G$, and $P \in \poly(\Gdot)$ be $Q$-rational with respect to $\Gamma$.  Define $H = G / G^{(d)}$, and let $\pi: G \to H$ be the quotient map. Consider the following filtrations: a) the quotient filtration $\Hdot$, a filtration of degree $d-1$ on $H$; b) the induced filtration $G^{(d)\bullet}$, a filtration of degree $d$ on $G^{(d)}$; and c) the induced filtration $(\pi \Gamma)^\bullet$ on the subgroup $\pi \Gamma \leq H$ induced from $\Hdot$, a filtration of degree $d-1$. It follows by polynomial fact \ref{item_poly_fact_four} above that the quotient map $\pi: \poly(\Gdot) \to \poly(\Hdot)$ has kernel $\ker(\pi) = \poly(G^{(d)\bullet})$.  It is quick to check from polynomial facts \ref{item_poly_fact_one} and \ref{item_poly_fact_three} that $\pi \poly(\Gammadot) = \poly\big((\pi \Gamma)^\bullet \big)$.

Since $P \in \poly(\Gdot)$ is $Q$-rational with respect to $\Gamma$ and $\pi$ is a homomorphism, the polynomial $\pi P \in \poly(\Hdot)$ is $Q$-rational with respect to $\pi \Gamma$. Thus, by the induction hypothesis, we have that $\partial_{N_0}' \pi P \in \poly\big((\pi \Gamma)^\bullet \big)$. Using \eqref{eqn_useful_identity_for_derivatives} to write $\partial_{N_0N_1}' \pi P$ as a product of translates of $\partial_{N_0}' \pi P$, it follows by polynomial fact \ref{item_poly_fact_two} above that  $\partial_{N_0N_1}' \pi P \in \poly\big((\pi \Gamma)^\bullet \big)$.

Since $\partial_{N_0N_1}' \pi P \in \poly\big((\pi \Gamma)^\bullet \big) = \pi \poly(\Gammadot)$, there exists $\rho \in \poly(\Gammadot)$ such that
\[\pi \rho = \partial_{N_0N_1}' \pi P = \pi \partial_{N_0N_1}' P.\]
It follows that $q \defeq \rho^{-1} \partial_{N_0N_1}' P \in \ker(\pi) = \poly(G^{(d)\bullet})$.  Note that $\rho^{-1}$ is a $1$-rational polynomial and $\partial_{N_0N_1}' P$ is the product of two $Q'$-rational polynomials, all with respect to $\Gamma$.  It follows by \cref{lemma_rationality_of_subgroup} that the polynomial $q \in \poly(G^{(d)\bullet})$ is $Q'$-rational with respect to $\Gamma \cap G^{(d)}$.

Now $G^{(d)}$ is an abelian group in the center of $G$, $\Gamma \cap G^{(d)} \leq G^{(d)}$, and $G^{(d)\bullet}$ is a filtration of degree $d$.  Since $q \in \poly(G^{(d)\bullet})$ is $Q'$-rational with respect to $\Gamma \cap G^{(d)}$, by \cref{lem_periodic_polys_modulo_subgroup_for_abelian_groups},
\[\partial'_{N_1} q \in \poly\big( (\Gamma \cap G^{(d)} )^{\bullet} \big),\]
where $(\Gamma \cap G^{(d)} )^{\bullet}$ denotes the induced filtration on $\Gamma \cap G^{(d)}$.
This implies that $q$ is $N_1$-periodic modulo $\Gamma \cap G^{(d)}$: for all $n_1, \ldots, n_\ell \in \Z$ that are congruent modulo $N_1$, there exists $\gamma \in \Gamma \cap G^{(d)}$ such that $q(n_1) \cdots q(n_\ell) = q(n_1)^\ell \gamma$.

Combining the definition of $N$ as $Q'!N_0N_1$, that $\partial_{N_0N_1}' P = \rho q$, the identity in \eqref{eqn_useful_identity_for_derivatives}, and that $G^{(d)}$ is in the center of $G$, it follows that for all $n \in \Z$, there exists $\gamma \in \Gamma$ such that
\begin{align*}
    \partial'_{N}P(n) & = \partial'_{N_0N_1}P \big( (Q'! - 1)N_0N_1 + n \big) \cdots \partial'_{N_0N_1}P (n) \\
    &= \rho \big( (Q'! - 1)N_0N_1 + n \big) \cdots \rho (n) \cdot q \big( (Q'! - 1)N_0N_1 + n \big) \cdots q (n) \\
    &= \rho \big( (Q'! - 1)N_0N_1 + n \big) \cdots \rho (n) \cdot q(n)^{Q'!} \gamma.
\end{align*}
Since $\rho$ is $\Gamma$-valued and $q(n)^{Q'!} \in \Gamma$ (since $q$ is $Q'$-rational), $\partial'_{N}P(n) \in \Gamma$. This shows that $\partial'_{N}P$ in $\Gamma$-valued, which by polynomial fact \ref{item_poly_fact_three} above allows us to conclude that $\partial'_{N}P \in \poly(\Gammadot)$.
\end{proof}

\subsection{Rational polynomials in nilpotent Lie groups and nilmanifolds}
\label{sec_rationals_in_nilmanifolds}

Let $G$ be a $d$-step nilpotent Lie group and $\Gamma \leq G$ be a discrete, cocompact subgroup.  The space $X = G / \Gamma$ is a \emph{$d$-step nilmanifold}; with the quotient topology, it is compact.  The quotient map $G \to X$ is denoted by $\pi$, and we define $e_X = \pi (e_G)$.  We denote the left action of $G$ on $X$ by juxtaposition: for $g \in G$ and $x \in X$, $gx \defeq \pi(gh)$, where $h \in G$ is such that $\pi(h) = x$.  Note, then, that $\pi(gh) = g \pi(h)$ for all $g, h \in G$.

A nilmanifold can have many different representations as a homogeneous space of a nilpotent Lie group.  As our interests lie in the space $X$ and dynamics on it, we will assume that the nilpotent Lie groups appearing in this paper are spanned by their connected component and finitely many other commuting elements.  This assumption appears throughout the recent literature on this topic; see, for example, \cite[Ch. 10, Sec. 1.3]{host_kra_book_2018}.

The spaces of bi-infinite, $G$- and $X$-valued sequences $G^{\Z}$ and $X^{\Z}$ are endowed with the product topology.  The space $X^\Z$ is naturally identified with $G^\Z / \Gamma^\Z$ and, hence, as a quotient of $G^\Z$.  We denote the quotient map $G^\Z \to X^\Z$ by $\pi$ as well; under the identification, this is the same as applying the map $\pi: G \to X$ to each coordinate, so there is no risk for confusion.  The coordinate-wise left action of $G$ on $G^\Z$ and $X^\Z$ will be denoted by juxtaposition: for $g \in G$ and a sequence $p$ in $G^\Z$ or $X^\Z$, $(gp)(n) = gp(n)$.

We equip $G$ with a right-invariant Riemannian metric $d_G$.  Because $d_G$ is right invariant, it descends to a metric $d_X$ on the nilmanifold $X = G / \Gamma$ that is consistent with the quotient topology.  Having $d_G$, we fix a right-invariant metric $d_{G^{\Z}}$ on $G^{\Z}$ that is consistent with the product topology; it descends to a metric $d_{X^{\Z}}$ on $X^{\Z}$ that is consistent with the quotient topology.

A subgroup $H \leq G$ is \emph{rational (with respect to $\Gamma$)} if $\Gamma \cap H$ is cocompact in $H$.  The rational subgroups of $G$ are closed in $G$.  There are a number of equivalent characterizations of rational subgroups: see \cite[Ch. 10, Lemma 14]{host_kra_book_2018} and the discussion just after \cref{def_rat_in_nilmanifolds} below.  A \emph{subnilmanifold of $X$} is a set of the form $Hx$, where $H$ is a rational subgroup of $G$ and $x \in X$.  The terminology is justified by realizing $Hx$ as a quotient of $H$ by a discrete, cocompact subgroup of $H$, and thus as a nilmanifold itself.

If $\Gdot$ is a \emph{rational filtration} (every subgroup in it is a rational subgroup of $G$),  the space of polynomials $\poly(\Gdot) \subseteq G^\Z$ is a closed, $G$-invariant subgroup of $G^\Z$ (cf. \cite[Ch. 14, Prop. 12]{host_kra_book_2018}).  When $G$ possesses a discrete, cocompact subgroup $\Gamma$ (as it always will in this paper), the lower central series filtration on $G$ is rational; see \cite[Th. 5.1.8 (a), Cor. 5.2.2]{corwin_greenleaf_1990} or \cite[Ch. 10, Thm. 23]{host_kra_book_2018}.

\begin{definition}
\label{def_poly_orbits}
Let $\Gdot$ be a rational filtration on $G$, and denote the induced filtration on $\Gamma$ by $\Gammadot$.  The \emph{set of polynomial orbits with coefficients in $\Gdot$}, $\poly(X, \Gdot) \subseteq X^\Z$, is defined to be the image of $\poly(\Gdot) \subseteq G^\Z$ under the quotient map $\pi: G^\Z \to X^\Z$.  In other words, a polynomial orbit with coefficients in $\Gdot$ is a sequence $p \in X^\Z$ of the form $p(n) = P(n) e_X$, where $P \in \poly(\Gdot)$.  The group $\poly(\Gdot)$ acts on $\poly(X,\Gdot)$ on the left, and the stabilizer of the constant $e_X$ sequence is $\poly(\Gdot) \cap \Gamma^\Z = \poly(\Gammadot)$.  Thus, the set $\poly(X, \Gdot)$ is naturally identified with the quotient $\poly(\Gdot) / \poly(\Gammadot)$.
\end{definition}

Key to the main results in this paper is the fact that $\poly(X,\Gdot)$ is a compact, finite dimensional space.  It can be realized as a subnilmanifold of $X^{d+1}$ \cite[Ch. 14, Sec. 2.2]{host_kra_book_2018}, but we will not need the entire description.  We collect just what we need to know in the following theorem and describe part of the space in \cref{lem_poly_rep} in the next section.

\begin{theorem}[{\cite[Ch. 14, Prop. 12 \& Sec. 2.2]{host_kra_book_2018}}]
\label{lem_topology_on_poly_space}
Let $G$ be a nilpotent Lie group, $\Gamma \leq G$ be a discrete, cocompact subgroup, and $X = G/\Gamma$.  Let $\Gdot$ be a rational filtration on $G$ of degree $d \in \Nz$, and let $\Gammadot$ be the induced filtration on $\Gamma$.
\begin{enumerate}[label=(\Roman*)]
    \item The group $\poly(\Gammadot)$ is a cocompact subgroup of $\poly(\Gdot)$, whereby $\poly(X, \Gdot)$ is compact.
    \item For all $\eps > 0$ and $N \in \N$, there exists $\delta > 0$ such that for all $p, q \in \poly(X,\Gdot)$, if $d_X\big(p(i), q(i) \big) < \delta$ for all $i \in \{0, 1, \ldots, d\}$, then $d_X \big(p(i), q(i) \big) < \eps$ for all $i \in \{-N, \ldots, N\}$.
    \item A polynomial orbit in $\poly(X,\Gdot)$ is determined uniquely by its values at $(d+1)$-many consecutive integers.
\end{enumerate}
\end{theorem}

Recall the operator $\partial'_N: \poly(\Gdot) \to \poly(\Gdot)$ defined in \eqref{eqn_def_of_reverse_diff}: for $P \in \poly(\Gdot)$, $(\partial'_NP)(n) = P(n+N)^{-1}P(n)$.

\begin{lemma}
\label{lemma_condition_for_periodicity}
Let $\Gdot$ be a rational filtration on $G$.  For all polynomials $P \in \poly(\Gdot)$, the polynomial orbit $\pi P \in \poly(X, \Gdot)$ is $N$-periodic if and only if $\partial'_N P \in \poly(\Gammadot)$.
\end{lemma}

\begin{proof}
By polynomial fact \ref{item_poly_fact_three} from the previous section, we have that $\partial'_N P \in \poly(\Gammadot)$ if and only if for all $n \in \Z$, $P(n+N)^{-1} P(n) \in \Gamma$.  This happens if and only if $P(n+N)^{-1} P(n) \Gamma = \Gamma$.  This is equivalent to the equality between $(\pi P)(n) \defeq P(n)e_X$ and $(\pi P)(n+N) \defeq P(n+N)e_X$, ie., the $N$-periodicity of the polynomial orbit $\pi P$.
\end{proof}

The following is a rather technical result that is useful later on when showing that a certain set of periodic polynomial orbits is measurable with respect to the prolongation relation.

\begin{lemma}
\label{lem_periodic_prolong_lemma}
Let $\Gdot$ be a rational filtration on $G$. Suppose that $(g_n)_n \subseteq G$, $(P_n)_n \subseteq \poly(\Gdot)$, and $P, Q \in \poly(\Gdot)$ are such that
\begin{enumerate}[label=(\Roman*)]
    \item \label{item_periodic_prolong_lemma_two} $\lim_{n \to \infty} \pi P_n = \pi P$, and
    \item \label{item_periodic_prolong_lemma_one} $\lim_{n \to \infty} g_n P_n = Q$.
\end{enumerate}
If $N \in \N$ is such that $\partial'_N P \in \poly(\Gammadot)$, then $\partial'_N Q \in \poly(\Gammadot)$.
\end{lemma}

\begin{proof}
It follows from \ref{item_periodic_prolong_lemma_two} that there exists $(\gamma_n)_n \subseteq \poly(\Gammadot)$ such that $\lim_{n \to \infty} d_{G^\Z}(P_n, P \gamma_n) \allowbreak = 0$.  Since $\partial'_N: \poly(\Gdot) \to \poly(\Gdot)$ is continuous and $\partial'_N (g_n P_n) = \partial'_N P_n$, hypotheses \ref{item_periodic_prolong_lemma_two} and \ref{item_periodic_prolong_lemma_one} imply that for all $N \in \N$,
\begin{enumerate}[label=(\Roman*)]
\setcounter{enumi}{2}
    \item \label{item_periodic_prolong_lemma_four} $\lim_{n \to \infty} d_{G^\Z} \big( \partial'_N P_n, \partial'_N (P \gamma_n) \big) = 0$, and
    \item \label{item_periodic_prolong_lemma_three} $\lim_{n \to \infty} \partial'_N P_n = \partial'_NQ$.
\end{enumerate}

Suppose that $N \in \N$ is such that $\partial'_N P \in \poly(\Gammadot)$. The polynomial $\partial'_N (P \gamma_n)$ is a pointwise product of three polynomials in $\poly(\Gammadot)$: for $m \in \Z$, $\partial'_N (P \gamma_n) (m) = \gamma_n(m+N)^{-1}( \partial'_N P )(m) \gamma_n(m)$.  Since $\poly(\Gammadot)$ is closed under translations, inversions, and products (cf. polynomial fact \ref{item_poly_fact_two} in the previous section), we see that $\partial'_N (P \gamma_n) \in \poly(\Gammadot)$.

To show that $\partial'_N Q \in \poly(\Gammadot)$, we will show that for all $\eps > 0$, there exists $n \in \N$ such that $d_{G^\Z} \big(\partial'_N Q, \partial'_N (P \gamma_n) \big) < \eps$.  This suffices since $\partial'_N (P \gamma_n) \in \poly(\Gammadot)$ and $\poly(\Gammadot)$ is closed in $\poly(\Gdot)$. Let $\eps > 0$.  Choose $n \in \N$ such that the terms in \ref{item_periodic_prolong_lemma_four} and \ref{item_periodic_prolong_lemma_three} are both within $\eps / 2$ of their limits. We see that
\[d_{G^\Z} \big(\partial'_N Q, \partial'_N (P \gamma_n) \big) \leq d_{G^\Z} \big(\partial'_N Q, \partial'_N P_n \big) + d_{G^\Z} \big(\partial'_N P_n, \partial'_N (P \gamma_n) \big) < \eps,\]
as was to be shown.
\end{proof}

We move now to discuss the rationality of points and polynomial orbits in nilmanifolds.

\begin{definition}
\label{def_rat_in_nilmanifolds}
A point $x \in X = G / \Gamma$ is \emph{$Q$-rational (with respect to $\Gamma$)} if there exists a $Q$-rational group element $g \in G$ such that $\pi (g) = x$; it is \emph{rational} if it is $Q$-rational for some $Q \in \N$.  The set of $Q$-rational elements of $X$ is denoted $\Q_Q(X)$; thus, $\Q_Q(X) \defeq \pi \Q_Q(G, \Gamma).$  A polynomial orbit $p \in \poly(X,\Gdot)$ is $Q$-rational if all of its values are $Q$-rational. A subnilmanifold $Y$ of $X$ is \emph{rational} if it takes the form $Y = Hx$, where $H \leq G$ is rational (as defined at the top of this subsection) and $x \in X$ is rational.
\end{definition}

The word \emph{rational} is used to describe at least six different objects.  Here is a summary of its uses and the results we will need about them.
\begin{itemize}
    \item An element $g \in G$ is rational if there exists $n \in \N$ such that $g^n \in \Gamma$.
    \item A polynomial $P \in \poly(\Gdot)$ is rational if all of its values are rational.
    \item A subgroup $H \leq G$ is rational if $H \cap \Gamma$ is cocompact in $H$.  By \cite[Corollary 1.14]{leibman_rational_2006}, the subgroup $H$ is rational if and only if the rational elements of $G$ are dense in $H$.
    \item A point $x \in X$ is rational if it can be written as $x = \pi (g)$ where $g \in G$ is rational.
    \item A polynomial orbit $p \in \poly(X,\Gdot)$ is rational if all of its values are rational.
    \item A subnilmanifold $Y$ of $X$ is rational if it can be written as $H x$, where both the subgroup $H \leq G$ and the point $x \in X$ are rational. By \cite[Corollary 2.10]{leibman_rational_2006}, the subnilmanifold $Y$ is rational if and only if the rational points of $X$ are dense in $Y$.
\end{itemize}

\begin{lemma}
\label{lem_lifts_of_rational_elements_and_polys}
For all $Q \in \N$ and $d \in \Nz$, there exists $Q' \in \N$ for which the following holds.  Let $X = G / \Gamma$ be a nilmanifold and $\Gdot$ be a rational filtration on $G$ of degree $d$.  If $x \in X$ is $Q$-rational, then every element of $\pi^{-1} \{x\}$ is $Q'$-rational, where $\pi: G \to X$ denotes the quotient map.  If $p \in \poly(X, \Gdot)$ is $Q$-rational, there exists a $Q'$-rational $P \in \poly(\Gdot)$ such that $\pi P = p$.
\end{lemma}

\begin{proof}
Let $Q \in \N$ and $d \in \Nz$.  Let $Q'$ be from \cref{lemma_rationality_of_subgroup}. Let $G / \Gamma$ be a nilmanifold, and let $\Gdot$ be a rational filtration on $G$ of degree $d$. Suppose $x \in X$ is $Q$-rational.  By definition, there exists $g \in \Q_Q(G, \Gamma)$ such that $\pi g = x$.  Since $\pi^{-1} \{x\} = \{ g \gamma \ | \ \gamma \in \Gamma\}$, it follows by \cref{lemma_rationality_of_subgroup} that $\pi^{-1} \{x\} \subseteq \Q_{Q'}(G,\Gamma)$.

Let $p \in \poly(X, \Gdot)$ be $Q$-rational.  By definition, there exists $P \in \poly(\Gdot)$ such that $\pi P = p$.  Since the values of $p$ are $Q$-rational, it follows by the previous paragraph that the values of $P$ are $Q'$-rational.  Thus, the polynomial $P$ is $Q'$-rational.
\end{proof}

\begin{theorem}
\label{thm_any_equivalent_poly_is_periodic}
    For all $Q \in \N$ and $d \in \Nz$, there exists $N \in \N$ for which the following holds.  Let $X = G / \Gamma$ be a nilmanifold and $\Gdot$ be a rational filtration on $G$ of degree $d$.  If $p \in \poly(X, \Gdot)$ has $(d+1)$-many consecutive $Q$-rational values, then $p$ is $N$-periodic.
\end{theorem}

\begin{proof}
Let $Q \in \N$ and $d \in \Nz$.  Let $Q'$ be from \cref{lem_lifts_of_rational_elements_and_polys}, and let $Q''$ be from \cref{lem_coeffs_and_values} for $Q'$ and $d$.  Let $N$ be from \cref{thm_periodic_polys} for $Q''$ and $d$.

Let $G / \Gamma$ be a nilmanifold, and let $\Gdot$ be a rational filtration on $G$ of degree $d$. Suppose that $p \in \poly(X, \Gdot)$ has $(d+1)$-many consecutive $Q$-rational values.  Let $P \in \poly(\Gdot)$ be such that $\pi P = p$.  Since $(d+1)$-many consecutive values of $p$ are $Q$-rational, $(d+1)$-many consecutive values of $P$ are $Q'$-rational.  It follows from \cref{lem_coeffs_and_values} \ref{item_rat_poly_conclusion_three} that $P$ is $Q''$-rational. It follows from \cref{thm_periodic_polys} that $\partial'_{N} P \in \poly(\Gammadot)$, and then \cref{lemma_condition_for_periodicity} gives that $p$ is $N$-periodic.
\end{proof}

\subsection{Nilsystems, nil-Bohr sets, and a subspace of polynomial orbits}
\label{sec_nilsystems}

Let $G$ be a $d$-step nilpotent Lie group, $\Gamma \leq G$ be a discrete, cocompact subgroup, and $X = G/\Gamma$ be the associated $d$-step nilmanifold. The group $G$ acts on the nilmanifold $X$ via \emph{nilrotations}, $T_g: x \mapsto gx$, $g \in G$.  Fixing $g \in G$, the pair $(X,T_g)$ is a \emph{$d$-step nilsystem}.  A single, fixed nilrotation $T_g$ will generally be denoted by $T$.  If $g_1, g_2 \in G$ commute, then the system $(X,T_{g_1},T_{g_2})$ is a \emph{$d$-step $\Z^2$-nilsystem}.  Given a minimal nilsystem or $\Z^2$-nilsystem, by changing the representation of the nilmanifold, we can assume that the parent nilpotent Lie group $G$ is spanned by its identity component and the finitely many commuting elements which define the nilrotations on the system; see \cite[Ch. 11, Sec. 1.2]{host_kra_book_2018}.

For this section, it may be useful to recall the dynamical system terminology established in \cref{sec_top_set_and_systems}.  The following two facts about nilsystems will be of particular use.  
\begin{enumerate}[label=(\Roman*)]
    \item \label{item_nil_is_distal} \cite{keynes_1967} Nilsystems and $\Z^2$-nilsystems are distal.  Thus, if $(X,T)$ is a nilsystem and $x \in X$, the system $(\orbitz{T}{x},T)$ is minimal.  Distality will allow us to apply the main results in \cref{sec_prolong} to $\Z^2$-nilsystems.
    \item \label{item_nil_has_nil_orbits} \cite[Ch. 11, Thm. 9]{host_kra_book_2018} For all $g \in G$ and $x \in X$, the orbit closure $\orbitz{T_g}{x}$ is a subnilmanifold of $X$ containing $x$.  By \cite[Corollary 2.10]{leibman_rational_2006}, we conclude that if $x$ is rational, then the rational points of $X$ are dense in $\orbitz{T_g}{x}$.  Since compact manifolds have finitely many connected components and $\orbitz{T_g}{x}$ is a nilmanifold, we deduce that $\orbitz{T_g}{x}$ has finitely many connected components.
\end{enumerate}

Nil-bohr sets -- sets of return times of a point $x \in X$ to an open set $U \subseteq X$ under a nilrotation $T_g$ -- will be important objects of study in \cref{sec_comb_structure_in_nilbohr_sets}.

\begin{definition}
\label{def_nilbohr_set}
A set $A \subseteq \Z$ is a \emph{nil-Bohr set} (resp. \emph{nil-Bohr$_0$ set}) if there exists a nilsystem $(X,T)$, a non-empty, open set $U \subseteq X$, and a point $x \in X$ (resp. $x \in U$) such that $R(x,U)$ is non-empty and $R(x,U) \subseteq A$.  A \emph{Bohr set} is a nil-Bohr set from a nilsystem system $(X,T)$ of step at most $1$, a rotation on a finite dimensional torus.  A nil-Bohr set is \emph{totally minimal} if the nilsystem $(X,T)$ is totally minimal.
\end{definition}

There are two helpful modifications to the definition of nil-Bohr sets: the nilsystem $(X,T)$ can be assumed to be minimal, and the point $x$ in the definitions can be taken to be the base point $e_X$.  These modifications do not alter the class of nil-Bohr sets; this follows from fact \ref{item_nil_is_distal} above and the ``change of base point'' argument in \cite[Ch. 11, Sec. 1.2]{host_kra_book_2018}.  We will use several times the fact that totally minimal nil-Bohr sets are totally intersective, an immediate consequence of \cref{lemma_tot_min_tot_intersective}.\\

We now move to describe certain elements of the product space $X^{d+1}$ as polynomial orbits. Let $(X,T)$ be a minimal, $d$-step nilsystem.  Consider the product nilmanifold $X^{d+1}$, naturally identified with $G^{d+1} / \Gamma^{d+1}$, whose points $\vec x \in X^{d+1}$ we treat as functions $\{0, \ldots, d\} \to X$ so that $\vec x= (\vec x(0), \vec x(1), \ldots, \vec x(d))$.  Define the maps $\Tmap = T \times T \times \cdots \times T$ and $\Smap = T^0 \times T^1 \times \cdots \times T^d$, and consider the $\Z^2$-nilsystem
\begin{align}
\label{eqn_def_of_W}
    \big( \spaceD \defeq \overline{\big\{\Tmap^n \Smap^m \diag{e_X} \ \big| \ n,m \in \Z \big\}}, \Tmap, \Smap \big),
\end{align}
where $\diag{\ \cdot \ }: X \to X^{d+1}$ denotes the diagonal embedding (defined first in \cref{sec_top_set_and_systems}).
Note that $\diag{X} \subseteq W$ since $(X,T)$ is minimal and $\Tmap \diag{e_X} = \diag{T e_X}$. It follows by \cite[Corollary 2.10]{leibman_rational_2006} that $\spaceD$ is a rational subnilmanifold of $X^{d+1}$.  Since $(X^{d+1},\Tmap,\Smap)$ is distal, the system $(\spaceD,\Tmap,\Smap)$ is minimal; in fact, this holds much more generally: the system $(\spaceD,\Tmap,\Smap)$ is minimal assuming only that the system $(X,T)$ is minimal \cite[Thm. 5.1]{glasner_1994}.

Denote by $\Gdot$ is the lower central series filtration of $G$; it is a filtration of degree $d$.  The next lemma shows that $\spaceD$ can be realized as a subspace of the space of polynomial orbits.

\begin{lemma}
\label{lem_poly_rep}
There exists a continuous injection $[\ \cdot \ ]: \spaceD \to \poly(X,\Gdot)$ with the property that for all $\vec w \in \spaceD$ and $i \in \{0, 1, \ldots, d\}$, $[\vec w](i) = \vec w(i)$. Moreover, for all $\vec w \in \spaceD$, $x \in X$, and $n, k \in \Z$,
\begin{align}
\label{eqn_polynomial_relation_tmap}
    [\Tmap^n \vec w] &= g^n [\vec w], \\
\label{eqn_polynomial_relation}
    [\Smap^n \diag{x}](k) &= T^{nk} x.
\end{align}
\end{lemma}

\begin{proof}
Let $\vec w \in \spaceD$.  By the definition of $\spaceD$, there exist $(n_j)_j, (m_j)_j \subseteq \Z$ such that
\begin{align}
\label{eqn_ex_point_limits_to_vec_w}
    \lim_{j \to \infty} \Tmap^{n_j} \Smap^{m_j} \uex = \vec w.
\end{align}
For $j \in \N$, define $p_j: \Z \to X$ by $p_j(i) = g^{n_j} g^{m_j i} e_X$.  That $p_j \in \poly(X,\Gdot)$ follows from the facts in Sections \ref{sec_periodicity_of_rat_polys}. It follows from \eqref{eqn_ex_point_limits_to_vec_w} that for $i \in \{0, 1, \ldots, d\}$,
\[\lim_{j \to \infty} p_j(i) = \vec w(i).\]
It follows from \cref{lem_topology_on_poly_space} that the sequence $(p_j)_j$ converges in $\poly(X, \Gdot)$.  We define $[\vec w]$ to be the limit of this sequence: for all $i \in \Z$,
\[[\vec w](i) = \lim_{j \to \infty} p_j(i).\]

It is clear from the definition of $[\ \cdot \ ]: \spaceD \to \poly(X,\Gdot)$ that for all $\vec w \in \spaceD$ and all $i \in \{0, 1, \ldots, d\}$, $[\vec w](i) = \vec w(i)$.  In particular, if $[\vec w] = [\vec v]$, then $\vec w = \vec v$, and so the map $[\ \cdot \ ]$ is injective.
The continuity of the map $[\ \cdot \ ]$ follows immediately from \cref{lem_topology_on_poly_space}.

To see \eqref{eqn_polynomial_relation_tmap} and \eqref{eqn_polynomial_relation}, note that $[\Tmap^n \vec w]$ and $g^n [\vec w]$ are both polynomial orbits that agree on $\{0, \ldots, d\}$.  Similarly, the map $i \mapsto T^{ni} x$ is an element of $\poly(X,\Gdot)$ that agrees with the polynomial orbit $[\Smap^{n}\diag{x}]$ on the set $\{0, 1, \ldots, d\}$.  In both cases, \cref{lem_topology_on_poly_space} gives equality.
\end{proof}

The $\Smap$-orbits of points in $\spaceD$ will be of particular interest in the next section.  Since we will need to consider their connected components, it is convenient to introduce some more notation.  For $x \in X$, define the system
\begin{align}
\label{eqn_labeling_scheme_for_subnilmanifolds}
    \big(Z \defeq \orbitz{\Smap}{\diag{x}}, \Smap \big).
\end{align}
The minimal nilsystem $(Z,\Smap)$ depends on the point $x$, though $x$ is not indicated in the notation.  Define $\pi_1: X^{d+1} \to X$ by $\pi_1(\vec x) = \vec x(1)$, and note that $T \circ \pi_1 = \pi_1 \circ \Smap$.  Since $(X,T)$ is minimal, the map $\pi_1: (Z,\Smap) \to (X,T)$ is a factor map.  It follows that $\ncc(X) \mid \ncc(Z)$.

Given a $T$-adapted labeling $X = X_1 \cup \cdots \cup X_{\ncc(X)}$ of the connected components of $X$ (recall \cref{def_clopen_adapted_partitions}), we give an $\Smap$-adapted labeling $Z = Z_1 \cup \cdots \cup Z_{\ncc(Z)}$ of the connected components of $Z$ in such a way that
\begin{align}
\label{eqn_connected_component_labeling}
    \compc(\diag{x}) = \compc(x).
\end{align}
(Recall the definition of $\compc$ from \cref{sec_top_set_and_systems}.) It follows that $\pi_1 Z_i = X_i$ for all $i \in \Z$.

The following lemma shows that for all $\vec w \in Z$ and all $n \in \Z$, the dilated polynomial orbit $k \mapsto [\vec w](nk)$ is represented by an element of $Z$; it goes further to give the specific connected component in which that element lies.

\begin{lemma}
\label{lem_polys_closed_under_dilation}
Let $x \in X$, and define $(Z,\Smap)$ as in \eqref{eqn_labeling_scheme_for_subnilmanifolds}.  For all $\vec w \in Z$ and $n \in \Z$, there exists $\vec z \in Z_{n\compc(\vec w) - (n-1)\compc(\diag{x})}$ such that for all $k \in \Z$, $[\vec z](k) = [\vec w](nk)$.
\end{lemma}

\begin{proof}
Let $\vec w \in Z$ and $n \in \Z$.
Note that the points $\Smap^{\compc(\vec w) - \compc(\diag{x})} \diag{x}$ and $\vec w$ are in the same connected component of $Z$.  Since $\Smap^{\ncc(Z)}$ restricted to that component is totally minimal (cf. \cref{lem_total_minimality_of_power_v2}), there exists $(m_j)_{j} \subseteq \Z$ such that 
\begin{align}
    \notag
    \lim_{j \to \infty} \Smap^{\ncc(Z) m_j + \compc(\vec w) - \compc(\diag{x})} \diag{x} &= \vec w, \text{ whence }\\
    \label{eqn_limit_of_bracket_to_w}
    \lim_{j \to \infty} [\Smap^{\ncc(Z) m_j + \compc(\vec w) - \compc(\diag{x})} \diag{x}] &= [\vec w],
\end{align}
where the last conclusion follows from the continuity described in \cref{lem_poly_rep}. It follows from \eqref{eqn_polynomial_relation} that for all $k \in \Z$, $[\Smap^{n\ncc(Z) m_j + n(\compc(\vec w) - \compc(\diag{x}))} \diag{x}](k) = [\Smap^{\ncc(Z) m_j + \compc(\vec w) - \compc(\diag{x})} \diag{x}](nk)$.  Defining $p \in X^\Z$ by $p(k) = [\vec w](nk)$, it follows by combining \eqref{eqn_limit_of_bracket_to_w} with the previous line that
\[\lim_{j \to \infty} [\Smap^{n\ncc(Z) m_j + n(\compc(\vec w) - \compc(\diag{x}))} \diag{x}] = p.\]
Since $\poly(X,\Gdot)$ is closed in $X^\Z$, this shows that $p \in \poly(X,\Gdot)$.  It also shows that $\vec z \defeq (p(k))_{k=0}^d$ is an element of $Z$ and $\compc(\vec z) =  n (\compc(\vec w) - \compc(\diag{x})) + \compc(\diag{x})$.  Since $[\vec z]$ and $p$ are polynomial orbits that agree on $\{0, 1, \ldots, d\}$, it follows by \cref{lem_topology_on_poly_space} that $[\vec z] = p$, as was to be shown.
\end{proof}

\section{Simultaneous approximation in nilsystems}
\label{section_simultaneous_approx_in_nilsystems}

Let $(X,T)$ be an invertible system, $x \in X$, and $A \subseteq \Z$.  Recall from \eqref{eqn_def_of_sa_set_intro} in the introduction that the set $\sa(x,A)$ consists of the points in $X$ that are simultaneously approximated by $x$ under $\{T^n \ | \ n \in A\}$. We note four useful facts that follow quickly from the definition.  First,
\begin{align}
\label{eqn_description_of_sa_on_diagonal}
    \sa(x,A) = \bigcap \big\{ \Delta^{-1} \big( \ \longorbit{T^{n_1} \times \cdots \times T^{n_k}}{\diag{x}} \cap \Delta(X) \ \big) \ \big| \ k \in \N, \ n_1, \ldots, n_k \in A \big\},
\end{align}
where both $\Delta(\ \cdot \ )$ and $\diag{\ \cdot \ }$ are used to denote the diagonal injection $X \to X^k$ (as defined in \cref{sec_top_set_and_systems}).  This shows that the set $\sa(x,A)$ is closed.  Second, $T^t \sa(x,A) = \sa(T^t x,A)$.  Third, the set $\sa(x,nA)$ in the system $(X,T)$ is the same as the set $\sa(x,A)$ in the system $(X,T^n)$; thus, sometimes we make the assumption that $\gcd A = 1$ for simplicity. And fourth, if $\pi: (X,T) \to (Y,T)$ is a factor map, then $\pi \sa(x,A) \subseteq \sa(\pi x, A)$.

The main result in this section, \cref{thm_main_fullversion}, concerns set the $\sa(x,N\Z + n)$ in minimal nilsystems.  It will allow us to deduce the dense simultaneous approximation properties in \cref{maintheorem_nilsystems_have_simult_approx}, and it will also lead to proofs of Theorems \ref{maintheorem_ae_return_set_is_mult_thick}, \ref{thm_comb_equivalent_to_main_thm}, \ref{main_vdw_type_theorem}, and \ref{mainthm_nilbohr_sets_are_mult_thick} in \cref{sec_comb_structure_in_nilbohr_sets}.

\subsection{Simultaneous approximation in an irrational rotation}
\label{section_simul_approx_for_irrational_rotations}

We begin by analyzing simultaneous approximation in the simple setting of an irrational rotation of the 1-torus, $\T = \R / \Z$.  This will serve nicely to motivate the simultaneous approximation property definitions in \cref{def_dense_simult_approx}, and it will reveal the germ of the main idea behind \cref{maintheorem_nilsystems_have_simult_approx}.

Fix $\alpha \in \T \setminus \Q$, and consider the system $(\T,T: x \mapsto x + \alpha)$.  The group $\T$ acts on itself by isometries that commute with $T$, so $\sa(x,A) = x + \sa(0,A)$.  Denote by $N^{-1} \Z$ the subgroup of $\T$ generated by $1/N$, ie. $N^{-1} \Z \defeq \{0, 1/N, \ldots, (N-1)/N\}$.  The next result shows that we can explicitly determine $\sa(0,A)$ in this setting.

\begin{lemma}
\label{lem_description_of_sa_set_on_torus}
For all $x \in \T$ and non-empty $A \subseteq \Z$,
\begin{align}
    \label{eqn_best_description_of_sa}
    \sa \big(x,A \big) = x + N^{-1} \Z,
\end{align}
where $N = \gcd (A-A) / \gcd A$.
\end{lemma}

\begin{proof}
Let $x \in \T$ and $A \subseteq \Z$ be non-empty.  Since $\sa(x,A) = x + \sa(0,A)$, it suffices to prove the lemma in the case that $x=0$.  Note that $\sa(0,A)$ is the same as the set $\sa(0,A/\gcd A)$ in the system $(\T, T^{\gcd A})$, which is also an irrational rotation of $\T$.  Thus, by replacing $A$ with $A / \gcd A$, it suffices to prove \cref{lem_description_of_sa_set_on_torus} in the case that $\gcd A = 1$.  Finally, note that if $0 \in A$, then $\sa(0,A) = \{0\}$, and the result holds since $\gcd(A-A)/\gcd(A) = 1$. Thus, we proceed under the assumption that $A \subseteq \Z \setminus \{0\}$ is non-empty.

Define $N = \gcd(A-A)$, and note that there exists $n \in \Z \setminus \{0\}$ such that $A \subseteq N \Z + n$.  Furthermore, for all $a \in A$, we have that $\gcd(n,N) = \gcd(a,N) \mid a$, whereby $\gcd(n, N) \mid \gcd A$.  Since $\gcd A = 1$, we have that $N$ and $n$ are coprime.

By the definition of simultaneous approximation, we have that $\sa(0,A) = \cap_{F} \sa(0,F)$, where the intersection is over finite subsets $F \subseteq A$.  Therefore, by including more elements in a finite set $F$ if necessary, the equality in \eqref{eqn_best_description_of_sa} will follow if we show that for all finite $F \subseteq A$ for which $\gcd F = 1$ and $\gcd (F - F) = N$, $\sa (0,F ) = N^{-1} \Z$.

Let $F = \{n_1, \ldots, n_k\} \subseteq A$ be finite and such that $\gcd F = 1$ and $\gcd (F-F) = N$.  Note that since $\alpha$ is irrational,
\begin{align}
\label{eqn_suffices_to_show_for_irrat_rotation_zero}
    \longorbitz{T^{n_1} \times \cdots \times T^{n_k}}{\diag{0}} = \left \{ (n_1 t, \ldots, n_k t) \in \T^k \ \middle| \ t \in \T \right\}.
\end{align}
It follows by combining \eqref{eqn_description_of_sa_on_diagonal} and \eqref{eqn_suffices_to_show_for_irrat_rotation_zero} that
\begin{align}
\label{eqn_description_of_sa_inject}
    \Delta\big(\sa(0,F)\big) = \left \{ (n_1 t, \ldots, n_k t) \in \T^k \ \middle| \ t \in \T \right\} \cap \Delta(\T).
\end{align}
Note that for $t \in \T$ and distinct $n, m \in \Z$, the equality $n t = m t$ holds if and only if $t \in \left(n - m\right)^{-1} \Z$. Since $n^{-1} \Z \cap m^{-1}\Z = \gcd(n,m)^{-1}\Z$, we see that
\begin{align}
\label{eqn_equality_condition_for_toral_dilates}
    n_1t = n_2 t = \cdots = n_k t \quad \text{ if and only if } \quad t \in N^{-1} \Z.
\end{align}
Combining \eqref{eqn_description_of_sa_inject} and \eqref{eqn_equality_condition_for_toral_dilates} with the fact that $\gcd(n,N) = \gcd(n_i,N) = 1$, we see that
\[\Delta\big(\sa(0,F)\big) = \left \{ (n_1 t, \ldots, n_k t) \in \T^k \ \middle| \ t \in N^{-1}\Z \right\} = \Delta \big(N^{-1}\Z \big).\]
It follows that $\sa(0,F) = N^{-1}\Z$, as was to be shown.
\end{proof}

It follows from \cref{lem_description_of_sa_set_on_torus} and its proof that, roughly speaking, the set $\sa(0,A)$ is roughly $N^{-1}$-dense if and only if $A$ lies in an arithmetic progression of the shape $N \Z + n$, where $N$ is sufficiently large and $n$ is coprime to $N$.  We make this observation precise in the next two results. A strengthening of the ``if'' direction is given in \cref{lemma_special_case_of_main_thm}, which we prove in a way as to highlight the flow of the proof of the main theorem, \cref{thm_main_fullversion}, in \cref{sec_proof_of_main_thm}.

\begin{lemma}
\label{lem_delta_dense_implies_A_in_prog}
Let $A \subseteq \Z$ be such that $\gcd A = 1$, and let $\delta > 0$.  If $\sa(0,A)$ is $\delta$-dense in $\T$, then there exists $N \in \N$, $N \geq \delta^{-1}/2$, and $n \in \Z$ coprime to $N$ such that $A \subseteq N \Z + n$.
\end{lemma}

\begin{proof}
Define $N = \gcd (A-A)$.  It follows from \cref{lem_description_of_sa_set_on_torus} and its proof that $\sa(0,A) = N^{-1} \Z$ and that there exists $n \in \Z$ coprime to $N$ such that $A \subseteq N \Z + n$. Since $N^{-1} \Z$ is $(2N)^{-1}$-dense and $\sa(0,A)$ is $\delta$-dense, we have that $\delta \geq (2N)^{-1}$, whereby $N \geq \delta^{-1}/2$.
\end{proof}

For the next result, we introduce notation that aligns with the notation from \cref{lem_poly_rep}: for a point $t \in \T$, denote by $[t] \in \T^{\Z}$ the ``polynomial orbit'' defined by $[t](k) = kt$, $k \in \Z$.  The statement of \cref{lemma_special_case_of_main_thm} and use of this notation are intended to draw parallels with the statement and proof of \cref{thm_main_fullversion}.  Recall the definition of a totally minimal Bohr set from \cref{def_nilbohr_set}.

\begin{lemma}[Special case of \cref{theorem_advanced_main_theorem_on_returntime_sets}]
\label{lemma_special_case_of_main_thm}
For all $\delta > 0$, $N \in \N$ with $N \geq \delta^{-1} / 2$, and $n \in \Z$ coprime to $N$, there exists a $\delta$-dense set $Y \subseteq \T$ such that for all $y \in Y$, all finite $F \subseteq N \Z + n$, and all $\eps > 0$, the set
\begin{align}
    \label{eqn_target_set_for_grp_rotation}
    \big\{m \in \Z \ \big| \ \forall f \in F, \ d_{\T}(T^{mf}0, y) < \eps \big\}
\end{align}
is a totally minimal Bohr set; in particular, $Y \subseteq \sa(0,N\Z+n)$, so the set $\sa(0,N\Z+n)$ is $\delta$-dense in $\T$.
\end{lemma}

\begin{proof}
Let $\delta > 0$, $N \geq \delta^{-1}/2$, and $n$ be coprime to $N$. Put $Y = N^{-1}\Z$, and note that $Y$ is $\delta$-dense in $\T$.  Let $y \in Y$, $F \subseteq N \Z + n$ be finite, and $\eps > 0$.

Because $n$ and $N$ are coprime, there exists $q \in Y$ such that $nq = y$.  Since $Nq = 0$, the polynomial orbit $[q] \in \T^\Z$ is $N$-periodic.  By the continuity of the map $[ \cdot ]: \T \to \T^{\Z}$, if the point $m\alpha$ is sufficiently close to $q$, then for all $f \in F$, the point $[m\alpha](f) = T^{mf}0$ is within $\eps$ of $[q](f) = [q](Nf_0 + n) = [q](n) = y$. Such an $m$ is thus a member of the set in  \eqref{eqn_target_set_for_grp_rotation}.  Since $(X,T)$ is a totally minimal group rotation, the set of such $m$'s in \eqref{eqn_target_set_for_grp_rotation} is a totally minimal Bohr set.  Finally, note that $Y \subseteq \sa(0,N\Z+n)$ follows from the non-emptiness of the set in \eqref{eqn_target_set_for_grp_rotation}.
\end{proof}

Let us describe in broad strokes how this argument will be generalized to minimal nilsystems in \cref{sec_proof_of_main_thm}.  In the argument, we used crucially that the point $T^{m} 0$ being near an $N$-rational point implies that the point $T^{Nm} 0$ is near $0$.  This is no longer true for a nilrotation, but the idea can be salvaged by appealing to the fact that a sufficiently long finite portion of a polynomial orbit determines the entire orbit.  Thus, we will consider the set of those $m$'s for which the vector $(T^{0m} e_X, T^{1m} e_X, T^{2m} e_X, \ldots, T^{dm} e_X)$ is near an $N$-rational point of the product nilmanifold $X^{d+1}$.  Precisely which $N$-rational vector we wish to approximate is governed by \cref{lem_polys_closed_under_dilation}, which serves as the analogue to the assertion that ``there exists $q \in Y$ such that $nq = y$'' from the proof above.  And, just as we used that the $N$-rational points of the torus are $1/(2N)$-dense, we will appeal to the fact that the $N$-rational vectors of $X^{d+1}$ are sufficiently dense in certain rational subnilmanifolds.  We begin with this in the next section.

\begin{remark}
\label{rmk_kronecker_is_not_characteristic}
Straightforward, if not tedious, calculations along the lines of the ones here can be used to determine the sets $\sa(x,A)$ explicitly in other systems.  One such calculation in the minimal nilsystem $\big(\T^2, T: (x,y) \mapsto (x + \alpha, y + 2x + \alpha) \big)$ is useful in demonstrating that the Kronecker factor is not ``characteristic'' for simultaneous approximation, in the sense that the sets $\sa(x,A)$ are not generally expressible as a union of fibers above the Kronecker factor.  The Kronecker factor of $(\T^2,T)$ is the irrational rotation in the first coordinate.  We will show that for all $N \in \N$, the set $\sa \big((0,0), N \Z + 1 \big)$ is contained in $N^{-1} \Z \times (4N)^{-1}\Z$.  It will follow then by \cref{thm_main_fullversion} that the set $\sa \big((0,0), N \Z + 1 \big)$ is $\delta$-dense but not a union of fibers over the Kronecker factor.

To see that $\sa \big((0,0), N \Z + 1 \big) \subseteq N^{-1} \Z \times (4N)^{-1}\Z$, suppose $(y,z) \in \sa((0,0), N \Z + 1)$.  It follows from \cref{lem_description_of_sa_set_on_torus} that $y \in N^{-1}\Z$.  To see that $z \in (4N)^{-1}\Z$, we will show that for all $f \in N \Z + 1$, $(f^2 - 1 ) z = 0$.  This will show that $z \in \big(\gcd \{ f^2 - 1 \ | \ f \in N \Z + 1 \}\big)^{-1} \Z$, and a short calculation shows then that that $z \in (4N)^{-1} \Z$.  Let $f \in N\Z + 1$ and $\eps > 0$, and let $\eps' = \eps'(f,\eps) > 0$ be sufficiently small.  Since $(y,z) \in \sa((0,0), \{1,f\})$, there exists $m \in \Z$ such that both $T^m(0,0) = (m\alpha, m^2 \alpha)$ and $T^{mf}(0,0) = (mf \alpha, m^2 f^2 \alpha)$ are within $\eps'$ of $(y,z)$.  Since $\eps'$ is sufficiently small, the point $f^2 z$ is within $\eps$ of $z$, and so $(f^2 - 1)z$ is within $\eps$ of 0.  Since $\eps > 0$ was arbitrary, $(f^2 - 1)z = 0$, as was to be shown.
\end{remark}

\subsection{Uniform denseness of periodic polynomial orbits}
\label{sec_prelim_to_main}
\label{sec_prelim_to_main_periodic}

Let $(X,T)$ be a minimal, $d$-step nilsystem, where $X = G / \Gamma$ and $T = T_g$ for some $g \in G$.  We follow the notation and terminology laid out in \cref{sec_nilsystems}, including the system $(\spaceD,\Tmap,\Smap)$, where $\spaceD$ is a subnilmanifold of $X^{d+1}$, defined in \eqref{eqn_def_of_W}.  Key to the argument in \cref{sec_proof_of_main_thm} will be the periodicity derived from the rational elements of the product nilmanifold $X^{d+1}$.  Thus, we focus closely on the subset of $\spaceD$ of polynomial orbits considered in \cref{sec_nilsystems} that are periodic.

Let $\Gdot$ be the lower central series filtration on $G$.  Recall the continuous injection $[\ \cdot \ ]: W \to \poly(X,\Gdot)$ described in \cref{lem_poly_rep}.  For $N \in \N$, define
\begin{align}
    \label{eqn_def_of_PN}
    \perpoly{N} \defeq \big\{ \vec w \in \spaceD \ \big| \ \text{the polynomial orbit $[\vec w] \in \poly(X,\Gdot)$ is $N$-periodic} \big\}.
\end{align}
The set of $N$-periodic polynomial orbits in $\poly(X,\Gdot)$ is closed and $G$-invariant, where recall that $G$ acts on $\poly(X,\Gdot)$ coordinate-wise on the left.  Combining the continuity described in \cref{lem_poly_rep} with the fact that $[\Tmap \vec w] = g [\vec w]$, we see that $\perpoly{N}$ is a closed, $\Tmap$-invariant subset of $\spaceD$.

Denote by $\prolong_\Tmap$ and $\prolong_\Smap$ the prolongation relations on $\spaceD$ with respect to the maps $\Tmap$ and $\Smap$ defined in \eqref{def_prolongation}. Since $(\spaceD,\Tmap,\Smap)$ is a minimal, distal $\Z^2$-system, the results from \cref{sec_prolong} apply; in particular, both $\prolong_\Tmap$ and $\prolong_\Smap$ are equivalence relations on $\spaceD$.

One of the keys to proving that minimal nilsystems satisfy the a.e. dense simultaneous approximation property (part \ref{item_nilsystem_sats_ae_dense_approx} of \cref{maintheorem_nilsystems_have_simult_approx}) is to show that the set $\perpoly{N}$ is $\delta$-dense in the $\Smap$-orbit closures of almost all points in $\spaceD$.  This is accomplished in two steps: first we show that the set $\perpoly{N}$ is measurable with respect to the $\Tmap$-prolongation equivalence relation (\cref{lem_periodic_prolong_result}), then we appeal to the uniform denseness result from \cref{sec_uniform_denseness_in_prolongation_classes} (\cref{thm_uniform_denseness_of_msble_sets}).  Recall the definition of ``measurable'' from the beginning of \cref{sec_prolong}.

\begin{lemma}
\label{lem_periodic_prolong_result}
For all $N \in \N$, the set $\perpoly{N}$ is measurable with respect to the $\Tmap$-prolongation equivalence relation.
\end{lemma}

\begin{proof}
Let $N \in \N$. To show that $\perpoly{N}$ is $\Tmap$-prolongation measurable, we must show that for all $\vec w \in \perpoly{N}$, $\prolong_\Tmap(\vec w) \subseteq \perpoly{N}$.  Let $\vec w \in \perpoly{N}$ and $\vec v \in \prolong_\Tmap(\vec w)$: there exist $(\vec{w_n})_n \subseteq \spaceD$ and $(k_n)_n \subseteq \Z$ such that $\lim_{n \to \infty} \vec{w_n} = \vec w$ and $\lim_{n \to \infty} \Tmap^{k_n}\vec{w_n} = \vec v$. We will show that $\vec v \in \perpoly{N}$.

Let $P \in \poly(\Gdot)$ and $(\tilde{P_n})_n \subseteq \poly(\Gdot)$ be such that $\pi P = [\vec w]$ and $\pi \tilde{P_n} = [\vec {w_n}]$, where $\pi: \poly(\Gdot) \to \poly(X,\Gdot)$ is as in \cref{def_poly_orbits}.  By \cref{lem_topology_on_poly_space}, there exists a compact set $C \subseteq \poly(\Gdot)$ such that $C \poly(\Gammadot) = \poly(\Gdot)$.  Thus, for each $n \in \N$, there exists $c_n \in C$ and $\gamma_n \in \poly(\Gammadot)$ such that $g^{k_n} \tilde{P_n} = c_n \gamma_n$.  Define $(P_n)_n \subseteq \poly(\Gdot)$ by $P_n = \tilde{P_n} \gamma_n^{-1}$, so that $\pi P_n = [\vec {w_n}]$ and $g^{k_n} P_n = c_n \in C$.

Since $C$ is compact, by passing to a subsequence, there exists $Q \in \poly(\Gdot)$ such that $\lim_{n \to \infty} g^{k_n} P_n = Q$.  Since $\lim_{n \to \infty} \vec{w_n} = \vec w$, we have by the continuity in \cref{lem_poly_rep} that $\lim_{n \to \infty} \pi P_n = \pi P$. Thus, the hypotheses of \cref{lem_periodic_prolong_lemma} are satisfied with $g^{k_n}$ as $g_n$ and with $(P_n)_n$, $P$, and $Q$ as they are.  Since $[\vec w]$ is $N$-periodic and $\pi P = [\vec w]$, we have by \cref{lemma_condition_for_periodicity} that $\partial'_N P \in \poly(\Gammadot)$.  It follows from \cref{lem_periodic_prolong_lemma} that $\partial'_N Q \in \poly(\Gammadot)$, and hence by \cref{lemma_condition_for_periodicity} that $\pi Q$ is $N$-periodic. The proof is complete if we show that $\pi Q = [\vec v]$, since $[\vec v]$ being $N$-periodic implies that $\vec v \in \perpoly{N}$.

Since $\lim_{n \to \infty} g^{k_n} P_n = Q$, we have on the one hand that $\lim_{n \to \infty} \pi \big(g^{k_n} P_n\big) = \pi Q$.  On the other hand, it follows from \eqref{eqn_polynomial_relation_tmap} that $\pi \big(g^{k_n} P_n\big) = g^{k_n} \pi P_n = g^{k_n} [\vec{w_n}] = [\Tmap^{k_n} \vec{w_n}]$.  Since $\lim_{n \to \infty} \Tmap^{k_n}\vec{w_n} = \vec v$, we have that $\lim_{n \to \infty} [\Tmap^{k_n}\vec{w_n}] = [\vec v]$.  Therefore, $\lim_{n \to \infty} \pi \big(g^{k_n} P_n\big) = [\vec v]$, verifying that $\pi Q = [\vec v]$.
\end{proof}

\begin{theorem}
\label{theorem_periodic_polys_are_uniformly_dense}
Let $(X,T)$, $(\spaceD,\Tmap,\Smap)$, $\Gdot$, and $\prolong_\Smap$ be as above.  Define
\[\Omega = \big\{ x \in X \ \big | \ \text{the point $\diag{x} \in \spaceD$ is a point of continuity of the map $\orbitz{\Smap}{}: \spaceD \to \fcspaced$}  \big\}.\]
\begin{enumerate}[label=(\Roman*)]
    \item \label{item_main_prereq_one} The set $\Omega$ is a residual, $T$-invariant subset of $X$, and for all $x \in \Omega$, $\orbitz{\Smap}{\diag{x}} = \prolong_\Smap(\diag{x})$.

    \item \label{item_main_prereq_two} There exists $\ccd \in \N$ such that for all $x \in \Omega$, the set $\orbitz{\Smap}{\diag{x}}$ has $\ccd$-many connected components.
    
    \item \label{item_main_prereq_three} For all $x_0 \in X$ and all $\delta > 0$, there exists $N \in \N$ such that for all $x \in \Omega \cup \{x_0\}$, the set $\perpoly{N}$ (defined in \eqref{eqn_def_of_PN}) is $\delta$-dense in every connected component of $\orbitz{\Smap}{\diag{x}}$.
\end{enumerate}
\end{theorem}

\begin{proof}
\ref{item_main_prereq_one} Since $\orbitz{\Smap}{\diag{Tx}} = \Tmap \orbitz{\Smap}{\diag{x}}$ and $\Tmap$ is continuous, the set $\Omega$ is $T$-invariant. That $\Omega$ is residual will follow from \cite[Prop. 4.4]{glasscock_koutsogiannis_richter_2019}, but the inclusion of a ``zero coordinate'' in $\spaceD$ requires a bit of care when applying that result.

Define
\[Y = \overline{\big\{ (T \times \cdots \times T)^n (T^1 \times \cdots \times T^d)^m (e_X, \ldots, e_X) \ \big| \ n, m \in \Z \big\} } \subseteq X^d,\]
where $(e_X, \ldots, e_X) \in X^d$. Consider the map $Y \to \fcy$ and the map $\spaceD \to \fcspaced$ defined by
\begin{align}
    \label{eqn_set_map_forward_orbit_one}
    \vec y &\mapsto \overline{\big\{ (T^n \vec y(1), T^{2n} \vec y(2), \ldots, T^{dn} \vec y(d)) \ \big| \ n \in \N \cup \{0\} \big\}};\\
    \label{eqn_set_map_forward_orbit_two}
    \vec w &\mapsto \{\vec w(0)\} \times \overline{\big\{ (T^n \vec w(1), T^{2n} \vec w(2), \ldots, T^{dn} \vec w(d)) \ \big| \ n \in \N \cup \{0\} \big\}}.
\end{align}
It is shown in \cite[Prop. 4.4]{glasscock_koutsogiannis_richter_2019} that the points of continuity of the map in \eqref{eqn_set_map_forward_orbit_one} that lie on the diagonal in $X^d$ form a residual subset of the diagonal.  Since $X \to \fcx$, $x \mapsto \{x\}$, is continuous, the map in \eqref{eqn_set_map_forward_orbit_two} is continuous at $\vec w$ if and only if the map in \eqref{eqn_set_map_forward_orbit_one} is continuous at $(\vec w(1), \ldots, \vec w(d))$.  It follows that the points of continuity of the map in \eqref{eqn_set_map_forward_orbit_two} intersected with $\diag{X}$ form a residual subset of $\diag{X}$.  Since the system $(Y, T \times T^2 \times \cdots \times T^d)$ is distal, forward orbits are equal to full orbits, so the map in \eqref{eqn_set_map_forward_orbit_two} is the same as the map $\orbitz{\Smap}{}: \spaceD \to \fcspaced$.  Combining the previous two sentences, the points of continuity of the map $\orbitz{\Smap}{}: \spaceD \to \fcspaced$ intersected with $\diag{X}$ form a residual subset of $\diag{X}$, as was to be shown.
    
As mentioned in \cref{sec_prolongation_relation}, it follows from \cite[Thm. 1, Lemmas 2, 3]{akin_glasner_1998} that $\vec w \in \spaceD$ is a point of continuity of $\orbitz{\Smap}{}: \spaceD \to \fcspaced$ if and only if $\orbitz{\Smap}{\vec w} = \prolong_{\Smap}(\vec w)$.  Thus, for all $x \in \Omega$, $\orbitz{\Smap}{\diag{x}} = \prolong_\Smap(\diag{x})$.
    
\ref{item_main_prereq_two} Let $\ccd$ be from \cref{lemma_constant_number_of_ccs} for the system $(\spaceD,\Tmap,\Smap)$. For all $\vec w \in \spaceD$, the prolongation class $\prolong_{\Smap}(\vec w)$ has $\ccd$-many connected components.  Thus, for all $x \in \Omega$, the set $\prolong_\Smap(\diag{x}) = \orbitz{\Smap}{\diag{x}}$ has $\ccd$-many connected components.

\ref{item_main_prereq_three} Let $x_0 \in X$ and $\delta > 0$.  The quantity $N$ will be taken to be the product of $N_1$ and $N_2$, defined as follows.

Let $h \in G$ be such that $h e_X = x_0$.  Since $\orbitz{\Smap}{\diag{x_0}} = \diag{h} \orbitz{(\Smap')}{\diag{e_X}}$, where $\Smap' \defeq T_{h^{-1} g h}^0 \times T_{h^{-1} g h}^1 \times \cdots \times T_{h^{-1} g h}^d$ and juxtaposition denotes the product in the group $G^{d+1}$, the set $\orbitz{(\Smap')}{\diag{e_X}}$ is a rational subnilmanifold of $X^{d+1}$.  It follows that the set $\Q(X^{d+1})$ is dense in $\orbitz{(\Smap')}{\diag{e_X}}$, whereby the set $\diag{h} \Q(X^{d+1})$ is dense in $\orbitz{\Smap}{\diag{x_0}}$.  Therefore, there exists $Q_1 \in \N$ such that $\diag{h} \Q_{Q_1}(X^{d+1})$ is $\delta$-dense in every connected component of $\orbitz{\Smap}{\diag{x_0}}$. Let $N_1 = N_1(Q_1,d) \in \N$ be as guaranteed by \cref{thm_any_equivalent_poly_is_periodic}; it follows from that theorem and the $G$-invariance of $\perpoly{N_1}$ that $\diag{h} \Q_{Q_1}(X^{d+1}) \subseteq \perpoly{N_1}$, and hence that $\perpoly{N_1}$ is $\delta$-dense in every connected component of $\orbitz{\Smap}{\diag{x_0}}$.

To define $N_2$, let $\eps$ be from \cref{thm_uniform_denseness_of_msble_sets} for the system $(\spaceD,\Tmap,\Smap)$. Since $\spaceD$ is a rational subnilmanifold of $X^{d+1}$, the set $\Q(X^{d+1})$ is dense in $\spaceD$.  Therefore, there exists $Q_2 \in \N$ such that $\Q_{Q_2}(X^{d+1})$ is $\eps$-dense in $\spaceD$.  Let $N_2 = N_2(Q_2,d) \in \N$ be as guaranteed by \cref{thm_any_equivalent_poly_is_periodic}; it follows from that theorem that $\Q_{Q_2}(X^{d+1}) \subseteq \perpoly{N_2}$, and hence that $\perpoly{N_2}$ is $\eps$-dense in $\spaceD$.  \cref{lem_periodic_prolong_result} gives $\perpoly{N_2}$ is $\Tmap$-prolongation measurable, so it follows from \cref{thm_uniform_denseness_of_msble_sets} that for all $\vec w \in \spaceD$, the set $\perpoly{N_2}$ is $\delta$-dense in every connected component of $\prolong_{\Smap}(\vec w)$.  By \ref{item_main_prereq_one}, for all $x \in \Omega$, the set $\perpoly{N_2}$ is $\delta$-dense in every connected component of $\orbitz{\Smap}{\diag{x}}$.

Setting $N = N_1N_2$ and noting that $\perpoly{N_1} \cup \perpoly{N_2} \subseteq \perpoly{N}$, we have shown that for all $x \in \Omega \cup \{x_0\}$, the set $\perpoly{N}$ is $\delta$-dense in every connected component of $\orbitz{\Smap}{\diag{x}}$, as desired.
\end{proof}

\subsection{Simultaneous approximation in nilsystems: proof of \cref{maintheorem_nilsystems_have_simult_approx}}
\label{sec_proof_of_main_thm}

In this section, we prove \cref{{thm_main_fullversion}} and derive \cref{maintheorem_nilsystems_have_simult_approx} from it.  The basic idea in the proof of \cref{{thm_main_fullversion}} was outlined in \cref{lemma_special_case_of_main_thm} and the discussion following it.  Recall the definition of $T$-adapted from \cref{def_clopen_adapted_partitions}.

\begin{theorem}
\label{thm_main_fullversion}
Let $(X,T)$ be a minimal nilsystem.  There exists a residual, $T$-invariant set $\Omega \subseteq X$ for which the following holds.  Let $X = X_1 \cup \cdots \cup X_{\cx}$ be a $T$-adapted labeling of the connected components of $X$. For all $x_0 \in X$ and $\delta > 0$, there exist $\ccd, N \in \N$ with $\ccd \ncc(X) \mid N$ such that for all $x \in \Omega \cup \{x_0\}$ and all $i, I, n \in \Z$ with
\begin{align}
\label{eqn_conditions_in_main_thm}
    \gcd(n,N) = 1 \quad \text{ and } \quad i \equiv n I + \compc(x) \pmod{\cx},
\end{align}
there exists a $\delta$-dense set $Y \subseteq X_i$ such that for all $y \in \ddset$, all $\eps > 0$, and all finite $F \subseteq N \Z + n$, there exists a totally minimal nil-Bohr set $B \subseteq \Z$ for which
\begin{align}
\label{eqn_main_theorem_containment}
\ccd B + I \subseteq \big\{ m \in \Z \ \big| \ \forall  f \in F, \ d_{X} \big(T^{fm} x, y \big) < \eps \big\};
\end{align}
in particular, $Y \subseteq \sa(x,N \Z + n)$, so the set $\sa(x,N \Z + n)$ is $\delta$-dense in $X$.
\end{theorem}

\begin{proof}
Let $\Omega \subseteq X$ be the residual set described in \cref{theorem_periodic_polys_are_uniformly_dense}. Let $x_0 \in X$ and $\delta > 0$. Let $\ccd \in \N$ the the product of $\ncc(\orbitz{\Smap}{\diag{x_0}})$ and the $\ccd$ from \cref{theorem_periodic_polys_are_uniformly_dense} \ref{item_main_prereq_two}.  Let $N \in \N$ be from \cref{theorem_periodic_polys_are_uniformly_dense} \ref{item_main_prereq_three}. By replacing $N$ with $\ccd \ncc(X) N$ if necessary, we may assume that $\ccd \ncc(X) \mid N$; indeed, note that since $\perpoly{N} \subseteq \perpoly{\ccd \ncc(X) N}$, the conclusion of \cref{theorem_periodic_polys_are_uniformly_dense} \ref{item_main_prereq_three} still holds.

Let $x \in \Omega \cup \{x_0\}$. Define the system $(Z \defeq \orbitz{\Smap}{\diag{x}}, \Smap)$ as in \eqref{eqn_labeling_scheme_for_subnilmanifolds}, following the connected component labeling described in \eqref{eqn_connected_component_labeling}.  It follows from the discussion between \eqref{eqn_labeling_scheme_for_subnilmanifolds} and \eqref{eqn_connected_component_labeling} and the definition of $C$ that $\ncc(X) \mid \ncc(Z) \mid C$.

Let $i, I, n \in \Z$ satisfy the conditions in \eqref{eqn_conditions_in_main_thm}.  Let $n' \in \Z$ be a modulo $N$ inverse to $n$. We will show that the congruence in \eqref{eqn_conditions_in_main_thm} implies that there exists $m \in \Z$ such that
\begin{align}
\label{eqn_congruence_explaining_component_location}
    n' \big(m\cx + i- \compc(x)\big) \equiv I \pmod {\ccd}.
\end{align}
To see this, note that it follows from basic number theory that such an $m$ exists if and only if $n'(i-\compc(x)) \equiv I \pmod {\gcd(n'\cx,\ccd)}$.  Note that $\gcd(n' \cx, \ccd \big) = \cx$ since $\cx \mid \ccd$, $\ccd \mid N$, and $\gcd(n',N) = 1$.  By assumption, we have that $i \equiv nI + \compc(x) \pmod{\cx}$.  This easily rearranges to $n'(i - \compc(x)) \equiv I \pmod{ \ncc(X)}$, noting that $n$ and $n'$ are inverses modulo $\cx$ (since $\cx \mid N$).
    
Denote by $\pi_1: X^{d+1} \to X$ the ``evaluate at 1'' projection (where $d$ is the step of the nilsystem $(X,T)$), and define $\ddset = \pi_1 \big(\perpoly{N} \cap Z_{m \cx + i} \big)$.  By the way the connected components of $Z$ are labeled (described in \eqref{eqn_connected_component_labeling}), the map $\pi_1$ maps $Z_{m \cx + i}$ onto $X_{i}$.  Since $\perpoly{N}$ is $\delta$-dense in $Z_{m \cx + i}$ and $\pi_1$ is 1-Lipschitz, we see that $\ddset$ is a $\delta$-dense subset of $X_{i}$.

Let $y \in \ddset$, $\eps > 0$, and $F \subseteq N \Z + n$ be finite.  By the definition of $\ddset$, there exists $\vec w \in \perpoly{N} \cap Z_{m \cx + i}$ such that $\vec{w}(1) = y$. Since $\vec w \in \perpoly{N}$, the polynomial orbit $[\vec w]$ is $N$-periodic.  By \cref{lem_polys_closed_under_dilation} (with $n'$ as $n$), there exists $\vec z \in Z_{n' (m \cx + i) - (n'-1) \compc(\diag{x})}$ such that for all $k \in \Z$, we have $[\vec z](k) = [\vec w](n' k)$; in particular, the polynomial orbit $[\vec z]$ is $N$-periodic and $[\vec z](n) = [\vec w](1) = y$. By \eqref{eqn_congruence_explaining_component_location}, the fact that $\ncc(Z) \mid \ccd$, and the fact that $\compc(\diag{x}) = \compc(x)$, we have that $Z_{n' (m \cx + i) - (n'-1) \compc(\diag{x})} = Z_{I+\compc(\diag{x})}$.

It follows from \cref{lem_poly_rep} that there exists $\eps' > 0$ such that if $\vec u, \vec v \in Z$ satisfy $d_{Z}(\vec u, \vec v) < \eps'$, then for all $f \in F$, $d_{X}([\vec u](f), [\vec v](f)) < \eps$. Define
\[B = \big\{ b \in \Z \ \big| \ d_{Z} \big(\Smap^{b \ccd + I} \diag{x}, \vec z \big) < \eps' \big\}.\]
Since $\ncc(Z) \mid \ccd$, by \cref{lem_total_minimality_of_power_v2}, the system $(Z_{I+\compc(\diag{x})}, S^{\ccd})$ is a totally minimal nilsystem.  Since the points $\Smap^{I} \diag{x}$ and $\vec z$ both belong to $Z_{I+\compc(\diag{x})}$, the set $B$ is a non-empty, totally minimal nil-Bohr set.

To see \eqref{eqn_main_theorem_containment}, let $b \in B$.  Since $d_{Z} \big(\Smap^{b \ccd + I} \diag{x}, \vec z \big) < \eps'$, we have for all $f \in F$ that $d_{X} \big( [\Smap^{b \ccd + I} \diag{x}](f), [\vec z](f) \big) < \eps$.  By \eqref{eqn_polynomial_relation}, we have that $[\Smap^{b \ccd + I} \diag{x}](f) = T^{f(\ccd b + I)} x$. Since $F \subseteq N \Z + n$ and $[\vec z]$ is $N$-periodic, for all $f \in F$, $[\vec z](f) = [\vec z](n) = y$.  Putting the previous sentences together, we find that for all $f \in F$, $d_{X} \big(T^{f(\ccd b + I)} x, y \big) < \eps$, verifying \eqref{eqn_main_theorem_containment}.

To conclude the proof, it suffices to show that for all $x \in \Omega \cup \{x_0\}$ and $i, n \in \Z$ with $n$ coprime to $N$, the set $Y$ is contained in $\sa(x,N\Z + n)$.  Let $x \in \Omega \cup \{x_0\}$ and $i, n \in \Z$ with $n$ coprime to $N$.  Since $\gcd(n,\ncc(X)) = 1$, there exists $I \in \Z$ for which $i \equiv nI + \compc(x) \pmod{\ncc(X)}$ holds.  The containment $Y \subseteq \sa(x,N\Z+n)$ follows from the non-emptiness of the sets in \eqref{eqn_main_theorem_containment}.
\end{proof}

The union ``$\Omega \cup \{x_0\}$'' appearing in the statement of \cref{thm_main_fullversion} allows us to deduce both of the conclusions in \cref{maintheorem_nilsystems_have_simult_approx} with one statement.  Were one interested in proving only the simultaneous approximation property for the point $x_0$, the proof of \cref{thm_main_fullversion} would remain nearly the same, but the auxiliary work leading up to the theorem would be lessened significantly; in particular, none of the results concerning the prolongation relation would be necessary.

\begin{proof}[Proof of \cref{maintheorem_nilsystems_have_simult_approx}]
Let $(X,T)$ be a minimal nilsystem.  We must show that every point $x_0 \in X$ satisfies the dense simultaneous approximation property and that the system $(X,T)$ satisfies the a.e. dense simultaneous approximation property.  Let $\Omega \subseteq X$ be the residual set from \cref{thm_main_fullversion}, and let $X = X_1 \cup \cdots \cup X_{\cx }$ be a $T$-adapted labeling of the connected components of $X$.  Let $x_0 \in X$ and $\delta > 0$, and let $N \in \N$ be as given in \cref{thm_main_fullversion}.  Both the dense simultaneous approximation property for $x_0$ and the a.e. dense simultaneous approximation property for $(X,T)$ follow immediately now from the conclusion of \cref{thm_main_fullversion} that for all $x \in \Omega \cup \{x_0\}$ the set $\sa(x, N\Z + n)$ is $\delta$-dense in $X$.
\end{proof}

\section{Multiplicative thickness of return-time sets}
\label{sec_comb_structure_in_nilbohr_sets}

In this section, we sharpen the connection between simultaneous approximation and the multiplicative quality of return-time sets described in \cref{sec_intro_thickness}.  In particular, we prove Theorems \ref{maintheorem_ae_return_set_is_mult_thick} \ref{thm_comb_equivalent_to_main_thm}, \ref{mainthm_combinatorial_main}, and \ref{mainthm_nilbohr_sets_are_mult_thick} from the introduction.  These theorems demonstrate that the typical return-time set and all nil-Bohr sets are multiplicatively thick in cosets of certain multiplicative subsemigroups, and they describe some of the combinatorial consequences of those facts.  The main result for nil-Bohr sets, \cref{mainthm_nilbohr_sets_are_mult_thick}, is a necessary building block for the other results, so we attend to it first.

Recall from the introduction that for $N \in \N$, the congruence subsemigroup $\nn$ consists of those natural numbers congruent to $1$ modulo $N$.  Analogously, a \emph{congruence subsemigroup in $(\Z \setminus \{0\},\cdot)$} is a subsemigroup of $(\Z \setminus \{0\},\cdot)$ of the form
\begin{align*}
    \zn \defeq \big\{ n \in \Z \ \big| \ n \equiv 1 \pmod N \big\} = N \Z + 1.
\end{align*}
A \emph{coset of $(\zn, \cdot)$} is a set of the form $I \zn$, where $I \in \Z \setminus \{0\}$.  A set $A \subseteq \Z$ is \emph{multiplicatively thick in $\zn$} if the set $A \cap \zn$ is \emph{thick} as a subset of the semigroup $(\zn, \cdot)$: for all finite $F \subseteq \zn$, there exists $m \in \zn$ such that $mF \subseteq A$.  The set $A$ is \emph{multiplicatively thick in $I \zn$} if any one (all) of the following equivalent conditions hold:
\begin{itemize}
    \item the set $A / I$ is multiplicatively thick in $\zn$ (recall from \cref{sec_top_set_and_systems} that $A / I = \{n \in \Z \ | \ nI \in A \}$);
    \item for all $L \in \N$, there exists $m \in \zn$ such that for all $\ell \in \{-L, \ldots, L\}$, $Im(\ell N + 1) \in A$;
    \item for all finite $F \subseteq I \zn$, there exists $m \in \zn$ such that $mF \subseteq A$;
    \item for all finite $F \subseteq \zn$, there exists $m \in I \zn$ such that $mF \subseteq A$.
\end{itemize}
Replacing $\Z$ with $\N$ leads to the definition of ``multiplicatively thick in a coset of a congruence subsemigroup of $\N$'' as given in \cref{sec_intro_thickness}.  The following lemma describes the relationship between the two notions.

\begin{lemma}
\label{lemma_mult_thick_from_z_to_n}
Let $N, I \in \N$ and $A \subseteq \Z$. If the set $A$ is multiplicatively thick in $I\zn$, then the set $A \cap \N$ is multiplicatively thick in $I\nn$.
\end{lemma}

\begin{proof}
Suppose $A$ is multiplicatively thick in $I\zn$. To see that the set $B \defeq A \cap \N$ is multiplicatively thick in $I\nn$, let $F \subseteq I \nn$ be finite.  Define $F' = I\zn \cap [-N \max F, N \max F]$.  Since $A$ is multiplicatively thick in $I\zn$, there exists $m \in \zn$ such that $mF' \subseteq A$.

If $m > 0$, then we see that $mF \subseteq mF' \cap \N \subseteq B$.  If, on the other hand, $m < 0$, then $-m(N-1)F \subseteq mF'$; indeed, this follows from the fact that $-(N-1) \equiv 1 \pmod N$ and the definition of $F'$.  Note that $-m(N-1) \in \nn$ and that, just as above, $-m(N-1)F \subseteq B$.  We have verified in either case that some element of $\nn$ dilates the set $F$ into the set $B$, verifying that $B$ is multiplicatively thick in $I\nn$.
\end{proof}

Sets that are multiplicatively thick in $\nn$ are ``multiplicatively large'' in that they support a dilation-invariant mean on $\nn$; equivalently, they have multiplicative upper Banach density in $\nn$ equal to 1 (see \cref{def_mult_den_and_synd}).  They are ``multiplicatively rich'' in that they contain an abundance of combinatorial configurations; indeed, they contain (many) dilates of any finite subset of the infinite arithmetic progression $\nn = N \Nz + 1$ and so, for example, can be easily shown to contain the classes of geo-arithmetic configurations described in \cite{bergelson_2005}.  In any infinite semigroup, thick sets are \emph{IP sets}: they contain all finite products of the terms from some sequence in the semigroup.  Thus, a set which is multiplicatively thick in $I\nn$ contains \emph{infinite Hilbert cubes}, configurations of the form
\[I \ \big\{ \ 1, \ n_1, \ n_2, \ n_1 n_2, \ n_3, \ n_1 n_3, \ n_2 n_3, \ n_1n_2n_3, \ n_4, \ \ldots \big\},\]
where the set contains all finite products of terms from some infinite sequence $(n_i)_i \subseteq \nn$.

\subsection{Multiplicative thickness of return-time sets in nilsystems}
\label{sec_mult_thick_in_nilsystems}

We prove in this section that return-time sets in minimal nilsystems and inverse limits of minimal nilsystems are multiplicatively thick in a coset of a congruence subsemigroup.  For certain combinatorial applications, it is useful to have the result for such sets ``up to zero additive density,'' so we begin by making that notion precise.

Recall the additive upper (Banach) density $d^*$ defined in \eqref{eqn_add_up_banach_den} in the introduction.
The quantity $d^*(A)$ is equal to the supremum of the values given to the indicator function $\one_A$ by translation-invariant means on $\Z$. The set $A$ is of \emph{zero additive density} if $d^*(A) = 0$, and it is of \emph{full additive density} if the set $\Z \setminus A$ is of zero additive density.  The set $A$ has a certain property \emph{up to zero additive density} if there exist a set $B \subseteq \Z$ with that property and a set $D \subseteq \Z$ of full additive density such that $A \supseteq B \cap D$.

\begin{lemma}
\label{lemma_dilation_property_of_full_density_sets}
Let $D \subseteq \Z$ be of full additive density.  For all finite $F \subseteq \Z \setminus \{0\}$, the set $\cap_{f \in F} D / f$ is of full additive density.
\end{lemma}

\begin{proof}
The set $E \defeq \Z \setminus D$ is a set of zero additive density.  Thus, for all $f \in \Z$, the set $E / f$ is of zero additive density.  If $F \subseteq \Z \setminus \{0\}$ is finite, the set $\cup_{f \in F} E / f$ is of zero additive density because it is a finite union of zero additive density sets.  The complement of the set $\cup_{f \in F} E / f$, $\cap_{f \in F} D / f$, is therefore a set of full additive density.
\end{proof}

The following is a strengthening of \cref{mainthm_nilbohr_sets_are_mult_thick} that gives information on the coset of the congruence subsemigroup in which a nil-Bohr set (up to zero additive density) is multiplicatively thick.

\begin{theorem}
\label{thm_strengthening_of_theorem_c}
Let $(X,T)$ be a minimal nilsystem. There exists a residual, $T$-invariant set $\Omega \subseteq X$ for which the following holds.  Let $X = X_1 \cup \cdots \cup X_{\cx}$ be a $T$-adapted partition consisting of the connected components of $X$.  For all $x_0 \in X$ and $\delta > 0$, there exists $N \in \N$ such that for all $x \in \Omega \cup \{x_0\}$, $y \in X$, non-zero $I \in \Z$ satisfying $I \equiv \compc(y) - \compc(x) \pmod {\ncc(X)}$, and $D \subseteq \Z$ of full additive density, the set $R\big(x,B(y,\delta) \big) \cap D$ is multiplicatively thick in $I\zn$.
\end{theorem}

\begin{proof}
Let $\Omega$ be the residual set from \cref{thm_main_fullversion}. Let $x_0 \in X$ and $\delta > 0$.  Let $C, N \in \N$ be as guaranteed by \cref{thm_main_fullversion}. Let $x \in \Omega \cup \{x_0\}$, $y \in X$, $I \in \Z \setminus \{0\}$ with $I \equiv \compc(y) - \compc(x) \pmod {\ncc(X)}$, and $D \subseteq \Z$ be of full additive density.

The conditions in \eqref{eqn_conditions_in_main_thm} are satisfied with $\compc(y)$ as $i$, $1$ as $n$, and $I$ as it is.  Therefore, there exists a $\delta$-dense set $Y \subseteq X_{\compc(y)}$ with the properties guaranteed by \cref{thm_main_fullversion}. Since $Y$ is $\delta$-dense, there exists $y' \in Y \cap B(y,\delta)$. Choose $\eps > 0$ such that $B(y',\eps) \subseteq B(y,\delta)$.

To see that the set $R(x,B(y,\delta)) \cap D$ is multiplicatively thick in $I \zn$, let $F \subseteq \zn$ be finite. We must show that there exists $m \in I \zn$ such that $mF \subseteq R(x,B(y,\delta)) \cap D$.  We will accomplish this by showing that
\begin{align}
\label{eqn_target_set_containment_for_thickness}
    I \zn \cap \bigcap_{f \in F} \frac{R(x,B(y,\delta))}{f} \cap \bigcap_{f \in F} \frac{D}{f} \neq \emptyset.
\end{align}

The conclusion of \cref{thm_main_fullversion} (with $y'$ as $y$ and $F$ and $\eps$ as they are) gives the existence of a totally minimal nil-Bohr set $B \subseteq \Z$ such that
\begin{align*}
    \ccd B + I \subseteq \bigcap_{f \in F} \frac{R(x,B(y',\eps))}{f} \subseteq \bigcap_{f \in F} \frac{R(x,B(y,\delta))}{f}.
\end{align*}
Define $B_0 = B / (I N)$.  Since $B_0$ is a totally minimal nil-Bohr set, the set $\ccd I N B_0 + I$ is a set of positive additive density (in fact, a nil-Bohr set) contained in $I \zn \cap (\ccd B + I)$.  Since $D$ is a set of full additive density, \cref{lemma_dilation_property_of_full_density_sets} gives that the set $\cap_{f \in F} D/f$ is a set of full additive density.  Therefore, we have that \eqref{eqn_target_set_containment_for_thickness} holds, as desired.
\end{proof}

\begin{remark}
More can be extracted from the conclusion of \cref{thm_main_fullversion} in the proof of \cref{thm_strengthening_of_theorem_c}.  By taking advantage of the total intersectivity of $B$, it can be shown that the set $R(x,U)$ is multiplicatively thick in other cosets of other congruence subsemigroups.  It is possible, if not tedious, to describe precisely the set of such cosets of congruence subsemigroups.  We do not have a ready application in mind, so we will leave this to the motivated reader.
\end{remark}

The analogue of \cref{thm_strengthening_of_theorem_c} for inverse limits of minimal nilsystems -- \cref{thm_strengthening_of_theorem_c_for_inverse_limits} -- follows readily.  It will be used in the next section in conjunction with topological characteristic factor machinery to prove \cref{maintheorem_ae_return_set_is_mult_thick}.  In this paper, we have need only for \emph{(sequential) inverse limits} of minimal nilsystems, though the inverse limit construction is far more general; see \cite[E.12]{devries_book_1993}.

\begin{theorem}
\label{thm_strengthening_of_theorem_c_for_inverse_limits}
\label{theorem_ae_return_set_in_minimal_nil_is_mult_thick}
Let $(X,T)$ be an inverse limit of minimal nilsystems (of any step). There exists a residual, $T$-invariant set $\Omega \subseteq X$ for which the following holds.  For all $x_0 \in X$ and $\delta > 0$, there exists $C, N \in \N$ and a $T$-adapted, clopen partition $X = X_1 \cup \cdots \cup X_{C}$ such that for all $x \in \Omega \cup \{x_0\}$, $y \in X$, non-zero $I \in \Z$ satisfying $I \equiv \compc(y) - \compc(x) \pmod C$, and $D \subseteq \Z$ of full additive density, the set $R \big(x,B(y,\delta)\big) \cap D$ is multiplicatively thick in $I \zn$.  If $(X,T)$ is totally minimal, then $C$ can be taken to be equal to $1$.
\end{theorem}

\begin{proof}
Let $(X,T)$ be an inverse limit of the family $\{(X_p,T_p) \ | \ p \in \N\}$ of minimal nilsystems, $(X,T) = \varprojlim_{p \in \N} (X_p,T_p)$, and denote the factor maps by $\pi_p: (X,T) \to (X_p,T_p)$. Fix a metric $d_X$ on the inverse limit consistent with its topology.  Let $\Omega_p \subseteq X_p$ be the residual, $T$-invariant set guaranteed by \cref{thm_strengthening_of_theorem_c}.  Define $\Omega = \cap_{p = 1}^\infty \pi_p^{-1} \Omega_p$, and note that by \cref{lem_factor_maps_pass_residuality}, the set $\Omega$ is a residual, $T$-invariant subset of $X$.

Let $x_0 \in X$, and $\delta > 0$. Choose $p \in \N$ sufficiently large (depending on the metric $d_X$) so that for all $y \in X$,
\begin{align*}
    \pi_p^{-1} B(\pi_p y, \delta / 2) \subseteq B(y,\delta).    
\end{align*}
Define $C = \ncc(X_p)$, and let $X_p = X_{p,1} \cup \cdots \cup X_{p,C}$ be a $T$-adapted partition of $X_p$ consisting of the connected components of $X_p$; note that if $(X,T)$ is totally minimal, then so is $(X_p,T_p)$, whence $C=1$ by \cref{lem_total_minimality_of_power_v2}.  Let $N \in \N$ be as guaranteed by \cref{thm_strengthening_of_theorem_c} for the system $(X_p,T_p)$ with $\pi_p x_0$ as $x_0$ and $\delta/2$ as $\delta$.  Define $X_i = \pi_p^{-1} X_{p,i}$ so that $X = X_1 \cup \cdots \cup X_{C}$ is a $T$-adapted, clopen partition of $X$, and note that for all $x \in X$, $\compc(x) = \compc(\pi_p x)$.

Let $x \in \Omega \cup \{x_0\}$, $y \in X$, $I \in \Z \setminus \{0\}$ with $I \equiv \compc(y) - \compc(x) \pmod C$, and $D \subseteq \Z$ be of full additive density.  Note that $\pi_p x \in \Omega_p \cup \{\pi_p x_0\}$ and $I \equiv \compc(\pi_p y) - \compc(\pi_p x) \pmod C$. It follows then from the conclusion of \cref{thm_strengthening_of_theorem_c} that the set $R_{T_p}(\pi_p x, B(\pi_p y, \delta / 2)) \cap D$ is multiplicatively thick in $I\zn$. Since
\[R_{T_p}(\pi_p x, B(\pi_p y, \delta / 2)) = R(x, \pi_p^{-1} B(\pi_p y, \delta / 2)) \subseteq R(x, B(y, \delta)),\]
the set $R(x, B(y, \delta)) \cap D$ is multiplicatively thick in $I\zn$, as was to be shown.
\end{proof}

Since the class of inverse limits of minimal nilsystems encompasses the class of minimal nilsystems, \cref{mainthm_nilbohr_sets_are_mult_thick} from the introduction follows immediately from the next result.

\begin{corollary}
\label{cor_nilbohr_sets_are_mult_thick}
Let $(X,T)$ be an inverse limit of minimal nilsystems (of any step).  For all $x \in X$, non-empty, open $U \subseteq X$, and $D \subseteq \Z$ of full additive density, there exists $N \in \N$ and $I \in \{1, \ldots, N\}$ such that the set $R(x,U) \cap D$ is multiplicatively thick in $I\zn$.  If $(X,T)$ is totally minimal, then $I$ can be taken to be equal to $1$.  In particular, all nil-Bohr sets up to zero additive density are multiplicatively thick in a coset of a congruence subsemigroup in $\N$ and in $\Z$, and all totally minimal nil-Bohr sets are multiplicatively thick in a congruence subsemigroup in $\N$ and in $\Z$.
\end{corollary}

\begin{proof}
Let $x \in X$, $U \subseteq X$ be non-empty, open, and $D \subseteq \Z$ be of full additive density.  Let $\delta > 0$ be such that the set $U$ contains a ball of radius $\delta$.  Let $C, N \in \N$ be as guaranteed by \cref{theorem_ae_return_set_in_minimal_nil_is_mult_thick} with $x$ as $x_0$.  Tracing back to the the application of \cref{thm_main_fullversion}, we see that $N \geq C$. Choosing $y \in U$ such that $B(y,\delta)\subseteq U$, it follows from \cref{theorem_ae_return_set_in_minimal_nil_is_mult_thick} that there exists $I \in \{1, \ldots, N\}$ such that set $R(x,U) \cap D$ is multiplicatively thick in $I\zn$.  If the system $(X,T)$ is totally minimal, then so are all of its factors.  In particular, we see from the proof of \cref{theorem_ae_return_set_in_minimal_nil_is_mult_thick} that $C = 1$, and thus that $I$ can be taken to be equal to $1$.  That the set $R(x,U) \cap D$ is multiplicatively thick in a coset of a multiplicative subsemigroup in $\N$ follows then from \cref{lemma_mult_thick_from_z_to_n}.
\end{proof}

It follows from \cref{thm_strengthening_of_theorem_c} and the proof of \cref{cor_nilbohr_sets_are_mult_thick} that -- provided a nil-Bohr set is described by the times of return of a point belonging to the set $\Omega$ -- the congruence subsemigroup $\zn$ is a function solely of the \emph{diameter} of the nil-Bohr set, the radius of the largest ball contained in the open set $U$ defining the set.  The coset of $\zn$ in which the nil-Bohr set is multiplicatively thick is described explicitly in the statement of \cref{thm_strengthening_of_theorem_c}.

\subsection{Combinatorial and dynamical applications for nil-Bohr sets}
\label{sec_comb_applications}

In this section, we give three applications of \cref{mainthm_nilbohr_sets_are_mult_thick}.  The first two show that certain families of well-studied, combinatorially defined sets have richer multiplicative combinatorial structure than previously known.  The last demonstrates a type of disjointness between minimal additive and multiplicative actions.

\subsubsection{Difference sets and their higher-order generalizations}

The existence of nil-Bohr structure in the families discussed below is well-known by experts but does not seem to be written explicitly in a sufficiently clear and general form in the literature.  Thus, we provide the short derivations here.

When $(X,\mu,T)$ is a measure preserving system and $Y \subseteq X$ is a measurable set of positive measure, the \emph{correlation sequence} $n \mapsto \mu(Y \cap T^{-n}Y)$ is positive definite.  As such, it follows by combining classical results of Herglotz and Weiner that the sequence can be written as a sum of an almost periodic sequence and a null-sequence, defined below.  Analogous decompositions of \emph{multi-correlation sequences} $n \mapsto \mu(T^{-p_1(n)}Y \cap \cdots \cap T^{-p_d(n)}Y)$ as sums of uniform limits of nilsequences (level sets of which are nil-Bohr sets) and null-sequences were achieved by Bergelson, Host, and Kra \cite{bergelson_host_kra_2005} for linear iterates and by Leibman \cite{leibman_2010} for polynomial iterates; see also \cref{rmk_more_general_decomp_applications} below.  We show in the corollary below how combinatorial applications for sets of the form $\{n \in \Z \ | \ (A - p_1(n)) \cap \cdots \cap (A - p_k(n)) \neq \emptyset \}$ can be achieved by combining those decomposition results with a suitable version of the Furstenberg correspondence principle.

We will need two definitions and an observation in the proof of the following theorem.  A \emph{basic nilsequence} is a sequence of the form $n \mapsto f(T^nx)$, where $(X,T)$ is a nilsystem, $x \in X$, and $f: X \to \R$ is continuous.  A sequence $\eta$ is a \emph{null-sequence} if for all $\eps > 0$, the set $\{n \in \Z \ | \ | \eta(n)| > \eps\}$ has zero additive upper density.  It is a short exercise, then, to check that if a non-negative sequence $n \mapsto \phi(n)$ is the sum a uniform limit of basic nilsequences and a null-sequence, and if for some $\eps > 0$, the level set $\{n \in \Z \ | \ \phi(n) > \eps \}$ has positive additive upper density, then the level set $\{n \in \Z \ | \ \phi(n) > \eps / 2 \}$ is a nil-Bohr set up to zero additive density.

\begin{theorem}
\label{cor_correlation_sequences_are_mult_thick}
Let $A \subseteq \Z$ be such that $d^*(A) > 0$.  For all polynomials $p_1, \ldots, p_d \in \Z[x]$ with $p_i(0) = 0$ and all $t \in \Z$, the set
\begin{align}
\label{eqn_times_of_multiple_returns}
    \big\{ n \in \Z \ \big| \ d^* \big((A-p_1(n)) \cap \cdots \cap (A-p_d(n)) \big) > 0 \big\} + t,
\end{align}
is a nil-Bohr set up to zero additive density, and it is multiplicatively thick in a coset of a congruence subsemigroup in $\N$ and in $\Z$.
\end{theorem}

\begin{proof}
The additive upper density of the set $A$, $d^*(A)$, is the supremum of $\lambda(A)$ over the set of translation-invariant means $\lambda$ on $\Z$.  Because $d^*(A) > 0$, there exists a translation-invariant mean $\lambda$ on $\Z$ for which $\lambda(A) > 0$.  According to a suitable version of the Furstenberg correspondence principle (cf. \cite[Theorem 3.1]{bergelson_leibman_2015}), there exists an invertible probability measure preserving system $(X,\mu,T)$ and a measurable set $Y \subseteq X$ with $\mu(Y) = \lambda(A)$ such that for all $n_1, \ldots, n_d \in \Z$,
\begin{align}
\label{eqn_correspond_principle}
    \lambda \big((A-n_1) \cap \cdots \cap (A-n_d) \big) = \mu \big( T^{-n_1} Y \cap \cdots \cap T^{-n_d} Y \big).
\end{align}

Let $p_1, \ldots, p_d \in \Z[x]$ with $p_i(0) = 0$.  By the polynomial \Szemeredi{} theorem \cite[Thm. A$_0$]{bergelson_leibman_1996}, there exists $\eps > 0$ for which the set
\begin{align*}
    \big\{ n \in \Z \ \big| \ \mu \big(T^{-p_1(n)} Y \cap \cdots \cap T^{-p_d(n)} Y \big) > \eps \big\}
\end{align*}
has positive additive upper Banach density.

According to \cite[Theorem 0.1]{Leibman_2015}, the multi-correlation sequence
\[n \mapsto \mu \big( T^{-p_1(n)} Y \cap \cdots \cap T^{-p_d(n)} Y \big)\]
is the sum of a uniform limit of basic nilsequences and a null-sequence.  As discussed above, it follows that the set
\begin{align*}
    \big\{ n \in \Z \ \big| \ \mu \big( T^{-p_1(n)} Y \cap \cdots \cap T^{-p_d(n)} Y \big) > \eps / 2 \big\}
\end{align*}
is a nil-Bohr set up to zero additive density. Translating by $t$ and applying the correspondence principle in \eqref{eqn_correspond_principle}, it follows that the set in \eqref{eqn_times_of_multiple_returns} is a nil-Bohr set up to zero additive density. It follows from \cref{cor_nilbohr_sets_are_mult_thick} that the set in \eqref{eqn_times_of_multiple_returns} is multiplicatively thick in a coset of a congruence subsemigroup in $\N$ and in $\Z$.
\end{proof}

For a set $A \subseteq \Z$, is quick to check that $A - A = \{n \in \Z \ | \ A \cap (A-n) \neq \emptyset\}$.  If $A$ has positive additive upper density, then taking $p_1(n) = 0$ and $p_2(n) = n$, it is a consequence of (the use of \cref{cor_nilbohr_sets_are_mult_thick} in) \cref{cor_correlation_sequences_are_mult_thick} that for all $t \in \Z$, there exists $N \in \N$ and $I \in \{1, \ldots, N\}$ such that the set $A-A+t$ is multiplicatively thick in $I\zn$, as observed above.  In fact, since $A-A$ is a Bohr$_0$ set up to zero additive density, the system underlying it is a 1-step nilsystem and, hence, equicontinuous.  In this case, the set $\Omega$ described in \cref{theorem_periodic_polys_are_uniformly_dense} is equal to the entire space $X$, and so we find that the value of $N$ can be taken to be independent of the translate $t$.

It is an exercise in the definitions to show that the set $\big\{ n \in \Z \ \big| \ (A-p_1(n)) \cap \cdots \cap (A-p_d(n)) \neq \emptyset \big\} + t$ is multiplicatively thick in $I\zn$ if and only if for all $L \in \N$, there exists $m \in \zn$ and $a_{-L}, a_{-L + 1}, \ldots, a_L \in \Z$ such that
\[\big\{ a_\ell + p_i(Im(\ell N + 1) - t) \ \big| \ \ell \in \{-L, \ldots, L\}, \ i \in \{1, \ldots, d\} \big\} \subseteq A.\]
Thus, this application of \cref{mainthm_nilbohr_sets_are_mult_thick} leads to novel affine polynomial configurations in sets of positive additive upper density.

\begin{remark}
\label{rmk_more_general_decomp_applications}
An \emph{approximate nilsequence} (cf. \cite[Ch. 24, Sec. 3.1]{host_kra_book_2018}) is a sequence that, for every $\eps > 0$, can be written as the sum of a nilsequence and a sequence $\eta$ satisfying $\lim_{N \to \infty} \sup_{M \in \Z} N^{-1} \sum_{n=M}^{M+N-1} | \eta(n)|^2 < \eps$.  Approximate nilsequences form a broader class of sequences than the ``nilsequence plus null-sequence'' class discussed above and have been used recently to describe more general classes of multi-correlation sequences; see, for example, \cite{frantzikinakis_2015,ferremoragues_2021,frantzikinakis_host_2018}.

Sufficiently large level sets of sums of nilsequences and null-sequences are nil-Bohr up to zero density, a fact key to the application in \cref{cor_correlation_sequences_are_mult_thick}.  It is natural to ask analogous questions about the level sets of approximate nilsequences.  If it true that positive additive upper density level sets of approximate nilsequences are multiplicatively thick in a coset of a subsemigroup, then we could avail ourselves of those more recent multi-correlation decomposition results for a broader range of combinatorial corollaries.
\end{remark}

\subsubsection{Times that generalized polynomials land in an open set}

The second application of \cref{mainthm_nilbohr_sets_are_mult_thick} concerns \emph{generalized polynomial mappings}, maps $\Z \to \R^d$ which are built from polynomials iteratively by applying addition, multiplication, and the integer part function, $\lfloor \ \cdot \ \rfloor$.  For example, the mapping $g: \Z \to \R$ given by
\[g(n) = 17 n \lfloor \pi n^3 \rfloor^2 - n^2 - \big\lfloor 17 n \lfloor \pi n^3 \rfloor^2 - n^2 \big\rfloor\]
is a generalized polynomial mapping.  It is \emph{bounded} since its image in $\R$ is bounded.  A generalized polynomial mapping into $\R^d$ is such a map in each coordinate.  A precise definition of the class can be found in \cite{bergelson_leibman_2007}, where it is shown that bounded generalized polynomial mappings can be represented dynamically via piecewise polynomial functions and nilrotations on nilmanifolds.

The set of ``return times'' of a bounded generalized polynomial mapping to an open set, if sufficiently large, is a nil-Bohr set and hence is multiplicatively thick in a coset of a congruence subsemigroup.  This is made precise in the following result.

\begin{theorem}
\label{cor_bdd_gen_poly_mult_thick}
Let $g: \Z \to \R^d$ be a bounded generalized polynomial mapping.  For all non-empty, open $U \subseteq \R^d$, if the set
\begin{align}
\label{eqn_gen_poly_return_times}
    \{n \in \Z \ | \ g(n) \in U \}
\end{align}
has positive additive upper Banach density, then it is a nil-Bohr set, and it is multiplicatively thick in a coset of a congruence subsemigroup in $\N$ and in $\Z$.
\end{theorem}

\begin{proof}
By \cite[Theorem A]{bergelson_leibman_2007}, there exists an ergodic, hence minimal \cite[Ch. 11, Cor. 5]{host_kra_book_2018}, nilsystem $(X,T)$, a point $x \in X$, and a piecewise polynomial mapping $f: X \to \R^d$ such that $g(n) = f(T^nx)$.  As explained in the proof of \cite[Theorem C]{bergelson_leibman_2007} (where ``$Y$'' is equal to $X$ since the system $(X,T)$ is minimal), while $f$ is in general discontinuous, there is a dense, open set $V \subseteq X$ on which $f$ is continuous.  According to \cite[Theorem B]{bergelson_leibman_2007}, the bounded piecewise polynomial surface $\mathcal{S} \defeq f(V)$ is such that the sequence $g$ is $\mathcal{S}$-valued and uniformly distributed (with respect to a natural measure on $\mathcal{S}$) on a set of full additive density.  By assumption, the set in \eqref{eqn_gen_poly_return_times} has positive additive upper density, so $U \cap \mathcal{S} \neq \emptyset$.  It follows that the set $f^{-1}(U) \cap V$ contains an open set, and hence that the set
\[\{n \in \Z \ | \ T^n x \in f^{-1}(U) \cap V \}\]
is a nil-Bohr set contained in the set in \eqref{eqn_gen_poly_return_times}.  Now \cref{cor_nilbohr_sets_are_mult_thick} applies to give that the set in \eqref{eqn_gen_poly_return_times} is multiplicatively thick a coset of a congruence subsemigroup in $\N$ and in $\Z$.
\end{proof}

Checking the hypotheses of \cref{cor_bdd_gen_poly_mult_thick} requires knowledge of the distribution of the values of $g$ and so may be difficult to do in general; some concrete examples are given in \cite[Sec. 0.10]{bergelson_leibman_2007}.  A more classical case for which the hypotheses are easier to verify is that of polynomial maps into the torus.  Let $p: \Z \to \R$ be a polynomial, $U \subseteq [0,1)$ be non-empty and open, and $\{ \ \cdot \ \}: \R \to [0,1)$ denote the fractional part map.  The map $n\mapsto  \{p(n)\}$ is a bounded generalized polynomial mapping.  If there exists $n \in \Z$ such that $\{p(n)\} \in U$, then the set in \eqref{eqn_gen_poly_return_times} has positive upper Banach density.  (This can be seen as a consequence of the distality of the systems created by Furstenberg \cite[Ch. 3, Sec. 3]{furstenberg_book_1981} to describe such sets; alternatively, see the discussion in \cite[Sec. 0.7]{bergelson_leibman_2007}.)  Therefore, in this case, the set in \eqref{eqn_gen_poly_return_times} is multiplicatively thick in some coset of a congruence subsemigroup $I \zn$.  For example, if $p(n) = \alpha n^2$, $\alpha \in \R \setminus \Q$, there exists $N \in \N$ and $I \in \Z \setminus \{0\}$ such that for all $L \in \N$, there exists $m \in \zn$ such that for all $\ell \in \{-L, \ldots, L\}$, $\{I^2m^2( \ell N + 1)^2 \alpha\} \in (1/4,1/3)$. In fact, the system in which $\{n \in \Z \ | \ \{ n^2 \alpha \} \in (1/4,1/3)\}$ is a return-time set is totally minimal, so, by \cref{mainthm_nilbohr_sets_are_mult_thick}, the value of $I$ can be taken to be $1$.

\subsubsection{Disjointness of minimal additive and multiplicative systems}
\label{sec_disjointness}

Our final application of \cref{mainthm_nilbohr_sets_are_mult_thick} concerns the disjointness of additive and multiplicative topological dynamical systems. Many famous results and conjectures in classical number theory can be phrased in terms of the degree to which multiplicative functions (eg., the Liouville function) correlate with additively structured functions (eg., periodic functions).  Recent results treating the case of nilpotent additive structure \cite{bergelson_richter_2022,green_tao_2012,frantzikinakis_host_2017} have led to powerful applications in number theory.

For clarity and emphasis in this section, we will call systems \emph{additive systems}.  A \emph{multiplicative system} $(Y,S)$ is an action of the semigroup $(\N, \cdot)$ on a compact metric space $Y$ by continuous maps: for $n, m \in \N$, we have that $S_n: Y \to Y$ is continuous and $S_{nm} = S_n \circ S_m$.  For a set $A \subseteq \N$ and $y \in Y$, we denote by $S_A y$ the set $\{S_n y \ | \ n \in A \}$.  The following definition is a topological analogue of the type of disjointness defined in \cite[Def. 1.22]{bergelson_richter_2022}.

\begin{definition}
\label{def_disjointness}
    An additive system $(X,T)$ and a multiplicative system $(Y,S)$ are \emph{disjoint} if for all $x \in X$ and all $y \in Y$,
    \[\overline{\big\{ (T^n x, S_n y) \ \big| \ n \in \N \big\}} = \overline{T^{\N}x} \times \overline{S_{\N} y}.\]
\end{definition}

We will address the disjointness of minimal additive and multiplicative systems.  To discuss minimality of multiplicative systems, we must define multiplicative syndeticity.

\begin{definition}
\label{def_mult_synd}
Let $I, N \in \N$ and $A \subseteq \N$. The set $A$ is \emph{multiplicatively syndetic in $\nn$} (resp. $\N$) if there exist $n_1, \ldots, n_k \in \nn$ (resp. $\N$) such that $A/n_1 \cup \cdots \cup A/n_k \supseteq \nn$ (resp. $\N$). The set $A$ is \emph{multiplicatively syndetic in $I\nn$} if the set $A/I$ is multiplicatively syndetic in $\nn$.
\end{definition}

The following facts are quick to check from the definitions.
\begin{enumerate}[label=(\Roman*)]
    \item \label{aux_synd_result_one} A set $A \subseteq \N$ is multiplicatively syndetic in $I \nn$ if and only if the set $I \nn \setminus A$ is not multiplicatively thick in $I \nn$.  This holds by an exercise using the definitions and the fact that $\big( (I \nn) \setminus A \big) \big / I = \nn \setminus (A / I)$, recalling the set algebra at the top of \cref{sec_prolong}.  Read another way, the set $A$ is multiplicatively syndetic in $I \nn$ if and only if it has non-empty intersection with all sets that are multiplicatively thick in $I \N_{N,1}$.

    \item \label{aux_synd_result_two} In a multiplicative system $(Y,S)$, every point $y \in Y$ has a dense orbit (ie. $S_{\N} y$ is dense in $Y$) if and only if for all $y \in Y$ and all non-empty, open $V \subseteq Y$, the set $\{n \in \N \ | \ S_n y \in V \}$ is multiplicatively syndetic in $\N$.\\
\end{enumerate}

A multiplicative system $(Y,S)$ is \emph{minimal} if for all $y \in Y$, the set $S_\N y$ is dense in $Y$. We will say that the multiplicative system is \emph{totally minimal} if for all $N \in \N$ and all $y \in Y$, the set $S_{\N_{N,1}} y$ is dense in $Y$.  Thus, by a straightforward analogue of \ref{aux_synd_result_two} above, a multiplicative system $(Y,S)$ is totally minimal if and only if for all $y \in Y$, all non-empty, open $V \subseteq Y$, and all $N \in \N$, the set $\{n \in \N \ | \ S_n y \in V \}$ is multiplicatively syndetic in $\N_{N,1}$. The next lemma records a small upgrade to this.

\begin{lemma}
\label{lemma_totally_min_mult_and_syndeticity}
    Let $(Y,S)$ be a totally minimal multiplicative system.  For all $y \in Y$, all non-empty, open $V \subseteq Y$, and all $N, I \in \N$, the set $\{n \in \N \ | \ S_n y \in V\}$ is multiplicatively syndetic in $I \N_{N,1}$.
\end{lemma}

\begin{proof}
    Let $y \in Y$, $V \subseteq Y$ be non-empty and open, and $N, I \in \N$.  Dividing by $I$, we must show that the set $\{n \in \N \ | \ S_n S_I y \in V\}$ is multiplicatively syndetic in $\N_{N,1}$.  This follows from the fact that the multiplicative system $(Y,S)$ is totally minimal.
\end{proof}

Here is an example of a totally minimal multiplicative system called a \emph{multiplicative rotation on two points} (cf. \cite[Rmk. 1.12]{bergelson_richter_2022}). Put $Y \defeq \{0,1\}$ and $S_n y \defeq y + \Omega(n) \pmod 2$, where $\Omega: \N \to \N$ is defined by
\[\Omega \big(p_1^{e_1} p_2^{e_2} \cdots p_k^{e_k} \big) = \sum_{i=1}^k e_i.\]
That $(Y,S)$ is a multiplicative system follows from the fact that $\Omega$ is multiplicative.  That $(Y,S)$ is totally minimal follows from the fact that $\Omega$ takes both even and odd values on every congruence subsemigroup.  For further examples, it is not hard to show that a multiplicative system that is \emph{aperiodic} (according to \cite[Def. 1.23]{bergelson_richter_2022}) is totally minimal.

We are now in a position to derive a disjointness result from \cref{mainthm_nilbohr_sets_are_mult_thick}.  Since aperiodic systems are totally minimal, in the topological category, this offers a generalization of \cite[Thm. C]{bergelson_richter_2022}.

\begin{theorem}
    Minimal nilsystems and totally minimal multiplicative systems are disjoint.
\end{theorem}

\begin{proof}
    Let $(X,T)$ be a minimal nilsystem and $(Y,S)$ be a totally minimal multiplicative system.  According to \cref{def_disjointness}, by the minimality of the systems, we must show: for all $(x,y) \in X \times Y$,
    \begin{align}
        \label{eqn_target_for_disjointness}
        \overline{\big\{ (T^n x, S_n y) \ \big| \ n \in \N \big\}} = X \times Y.
    \end{align}

    Let $(x,y) \in X \times Y$, and let $U \times V \subseteq X \times Y$ be open and non-empty. By \cref{mainthm_nilbohr_sets_are_mult_thick}, there exists $N, I \in \N$ such that the set $R_T(x,U)$ is multiplicatively thick in $I \N_{N,1}$.  By \cref{lemma_totally_min_mult_and_syndeticity}, the set $\{n \in \N \ | \ S_n y \in V\}$ is multiplicatively syndetic in $I\N_{N,1}$.  By \ref{aux_synd_result_one} above, there exists $n \in R_T(x,U)$ such that $S_n y \in V$.  Therefore, we have that $(T^n x, S_n y) \in U \times V$.  Since $U \times V \subseteq X \times Y$ was arbitrary, we have that \eqref{eqn_target_for_disjointness} holds, as desired.
\end{proof}

\subsection{Multiplicative thickness of return-time sets in general systems}
\label{sec_proof_of_return_sets_mult_thick}

In this section, we prove \cref{maintheorem_ae_return_set_is_mult_thick}, that almost all return-time sets $R(x,U)$ in a minimal system $(X,T)$ are multiplicatively thick in a coset of a congruence subsemigroup.  This is accomplished by lifting the result for inverse limits of nilsystems in \cref{theorem_ae_return_set_in_minimal_nil_is_mult_thick} to arbitrary systems via the recently-developed machinery of topological characteristic factors \cite{glasner_huang_shao_weiss_ye_2020}.

Denote by $\pi^{(d)}$ the map $\pi \times \cdots \times \pi$, and recall from previous sections that $\diag{x} = (x,\ldots, x)$, where the dimension of the embedding should always be clear from context.

\begin{definition}[{cf. \cite[Pg. 3]{glasner_huang_shao_weiss_ye_2020}}]
\label{def_char_factor}
Let $\pi: (X,T) \to (Y,T)$ be a factor map of invertible systems, and let $d \in \N$.  The system $(Y,T)$ is a \emph{$d$-step topological characteristic factor of $(X,T)$} if there exists a residual set $\Omega \subseteq X$ such that for all $x \in \Omega$, the orbit closure $\longorbitz{T \times T^2 \times \cdots \times T^d}{\diag{x}}$ is \emph{$\pi^{(d)}$-saturated}, meaning that
\[\big(\pi^{(d)}\big)^{-1} \big( \pi^{(d)} \longorbitz{T \times T^2 \times \cdots \times T^d}{\diag{x}} \big) = \longorbitz{T \times T^2 \times \cdots \times T^d}{\diag{x}}.\]
\end{definition}

When the orbit closure $\longorbitz{T \times T^2 \times \cdots \times T^d}{\diag{x}}$ is $\pi^{(d)}$-saturated, multiplicative configurations in the return-time set $R(x,U)$ can be ``lifted'' from those in the return-time set $R(\pi x, \pi U)$.  The following lemma makes this idea more precise.  Recall from the introduction that a point $x \in X$ is multiply recurrent under $\{T^n \ | \ n \in \N\}$ if it simultaneously approximates itself under $\{T^n \ | \ n \in \N\}$, that is, if for all open neighborhoods $U \subseteq X$ of $x$, the return-time set $R(x,U)$ is multiplicatively thick in $(\N,\cdot)$.

\begin{lemma}
\label{lem_top_char_for_thickness_of_return_sets}
Let $\pi: (X,T) \to (Y,T)$ be a factor map of invertible systems, and let $x \in X$. Suppose that the point $x$ is multiply recurrent under $\{T^n \ | \ n \in \N\}$ and for all $d \geq 2$, the orbit closure $\longorbitz{T \times T^2 \times \cdots \times T^d}{\diag{x}}$ is $\pi^{(d)}$-saturated.
\begin{enumerate}[label=(\Roman*)]
    \item \label{item_lem_top_factor_reformulation_v2_one} For all $n_1, \ldots, n_k \in \N$ and all non-empty, open $U \subseteq X^k$, if the return-time set $R_{T^{n_1} \times \cdots \times T^{n_k}} ( \diag{\pi x}, \pi^{(k)} U )$ is non-empty, then the set $R_{T^{n_1} \times \cdots \times T^{n_k}} (\diag{x}, U ) \cap \N$ is non-empty.
    
    \item \label{item_lem_top_factor_reformulation_v2_two} For all $I, N \in \N$ and all non-empty, open $U \subseteq X$, if the set $R(\pi x, \pi U)$ is multiplicatively thick in $I \zn$, then the set $R(x, U)$ is multiplicatively thick in $I\nn$.
\end{enumerate}
\end{lemma}

\begin{proof}
\ref{item_lem_top_factor_reformulation_v2_one} Let $n_1, \ldots, n_k \in \N$ and $U \subseteq X^k$ be non-empty and open.  Define $S = T^{n_1} \times \cdots \times T^{n_k}$ and $d = \max_i n_i$.  Since $x$ is multiply recurrent under $\{T^n \ | \ n \in \N\}$, $\orbit{S}{\diag{x}} = \orbitz{S}{\diag{x}}$.  Indeed, there exists an increasing sequence $(n_i)_i \subseteq \N$ such that $\lim_{i \to \infty} S^{n_i} \diag{x} = \diag{x}$.  For all $t \in \Z$, $\lim_{i \to \infty} S^{n_i + t} \diag{x} = S^t \diag{x}$.  Since $n_i + t$ is eventually positive, we see that $S^{\Z} \diag{x} \subseteq \orbit{S}{\diag{x}}$, whereby $\orbitz{S}{\diag{x}} \subseteq \orbit{S}{\diag{x}}$.

Note that $\orbitz{S}{\diag{x}}$ is $\pi^{(k)}$-saturated because, by assumption, $\longorbitz{T \times T^2 \times \cdots \times T^{d}}{\diag{x}}$ is $\pi^{(d)}$-saturated. Since $U$ is open, we see that
\begin{align*}
R_{S} ( \diag{\pi x}, \pi^{(k)} U ) \neq \emptyset & \quad \Longrightarrow \quad \orbitz{S}{\diag{\pi x}} \cap \big( \pi^{(k)} U \big) \neq \emptyset \\
& \quad \Longrightarrow \quad \big(\pi^{(k)} \orbitz{S}{\diag{x}} \big) \cap \big( \pi^{(k)} U \big) \neq \emptyset \\
& \quad \Longrightarrow \quad \orbitz{S}{\diag{x}} \cap U \neq \emptyset \\
& \quad \Longrightarrow \quad \orbit{S}{\diag{x}} \cap U \neq \emptyset \\
& \quad \Longrightarrow \quad R_{S} (\diag{x}, U ) \cap \N \neq \emptyset.
\end{align*}

\ref{item_lem_top_factor_reformulation_v2_two} Let $I, N \in \N$ and $U \subseteq X$ be non-empty, open.  Suppose that $R(\pi x, \pi U)$ is multiplicatively thick in $I\zn$, and let $F \subseteq I\nn$.  We will show that there exists $m \in \nn$ such that $mF \subseteq R(x,U)$.

Write $F = \{n_1, \ldots, n_k\}$, define $S = T^{n_1} \times \cdots \times T^{n_k}$, and write $U^k = U \times \cdots \times U$. By the set algebra from \cref{sec_top_set_and_systems} and \ref{item_lem_top_factor_reformulation_v2_one} above, we see that
\begin{align*}
\begin{aligned}
    \exists m \in \zn, \ mF \subseteq R(\pi x, \pi U) \quad
    &\Longrightarrow \quad \big( N\Z + 1 \big) \cap \bigcap_{f \in F} \frac{R(\pi x, \pi U)}{f} \neq \emptyset \\
    &\Longrightarrow \quad \big( N\Z \big) \cap \big(R_{S}(\diag{\pi x}, \pi^{(k)} U^k) - 1 \big) \neq \emptyset \\
    &\Longrightarrow \quad R_{S^N} \big(\diag{\pi x}, \pi^{(k)} S^{-1} U^k \big) \neq \emptyset\\
    & \Longrightarrow \quad R_{S^N}(\diag{x}, S^{-1} U^k) \cap \N \neq \emptyset \\
    & \Longrightarrow \quad \exists m \in \nn, \ mF \subseteq R(x,U).
\end{aligned}
\end{align*}
The set algebra required in the final implication is the same as that described in the first three implications.  Since $R(\pi x, \pi U)$ is multiplicatively thick in $I\zn$, the first statement holds, so the last statement holds.
\end{proof}

A factor map $\pi: (X,T) \to (Y,T)$ is an \emph{almost 1--1 extension} if there exists a residual set $\Omega \subseteq X$ such that for all $x \in X$, $\pi^{-1} (\pi x) = \{x\}$.  Almost 1--1 extensions are inherent to the topological characteristic factor machinery but do not present a challenge when considering return-time sets, as the following lemma shows.

\begin{lemma}
\label{lem_return_time_property_of_almost_onetoone_extensions}
Let $\pi: (X,T) \to (Y,T)$ be an almost 1--1 extension of invertible systems.  There exists a residual, $T$-invariant set $\Omega \subseteq X$ such that for all $x \in \Omega$ and all sets $U \subseteq X$,
\[R(\pi x, \pi U) = R(x,U).\]
\end{lemma}

\begin{proof}
Let $\Omega \subseteq X$ be the residual set in the definition of an almost 1--1 extension.  By replacing $\Omega$ with $\bigcap_{t \in \Z} T^t \Omega$, we may assume without loss of generality that $\Omega$ is $T$-invariant.

Let $x \in \Omega$ and $U \subseteq X$. By the definition of return-time sets that $R(x,U) \subseteq R(\pi x, \pi U)$.  To see the reverse inclusion, suppose $n \in R(\pi x, \pi U)$.  Since $\pi T^n x \in \pi U$, there exists $y \in U$ such that $\pi y = \pi T^n x$.  But $T^n x\in \Omega$, so $y = T^n x$.  Since $T^n x \in U$, we see that $n \in R(x,U)$, as was to be shown.
\end{proof}

\begin{theorem}
\label{theorem_advanced_main_theorem_on_returntime_sets}
Let $(X,T)$ be a minimal, invertible system.  There exists a residual, $T$-invariant set $\Omega \subseteq X$ for which the following holds.  For all $\delta > 0$, there exists $C, N \in \N$, a clopen, $T$-adapted partition $X = X_1 \cup \cdots \cup X_C$, and a $\delta$-dense set $Y \subseteq X$ such that for all $x \in \Omega$, $y \in Y$, and $I \in \N$ satisfying $I \equiv \compc(y) - \compc(x) \pmod C$, the set $R\big(x,B(y,\delta) \big)$ is multiplicatively thick in $I \nn$.  If $(X,T)$ is totally minimal, then $C$ can be taken to be equal to $1$.
\end{theorem}

\begin{proof}
Denote by $\pi: (X,T) \to (X_\infty,T)$ the maximal infinite-step pro-nilfactor of $(X,T)$ as defined in \cite[Sec. 2.4]{glasner_huang_shao_weiss_ye_2020}.  It follows by \cite[Thm. 3.6]{dong_donoso_maass_shao_ye_2013} that the system $(X_\infty,T)$ is an inverse limit of minimal nilsystems.  According to \cite[Theorem A]{glasner_huang_shao_weiss_ye_2020}, there exist minimal, invertible, almost 1--1 extensions $\sigma^*: (X^*,T) \to (X,T)$ and $\tau^*: (X_\infty^*,T) \to (X_\infty,T)$ and a factor map $\pi^*: (X^*,T) \to (X_\infty^*,T)$ such that $(X_\infty^*,T)$ is a $d$-step topological characteristic factor for $(X^*,T)$ for all $d \geq 2$ and the following diagram commutes.
\[ \begin{tikzcd}
\Omega \subseteq &[-30pt] (X,T) \arrow[swap]{d}{\pi} & \arrow[swap]{l}{\sigma^*} (X^*,T) \arrow{d}{\pi^*} &[-30pt] \supseteq \Omega_{\pi^*} \\%
\Omega_\infty \subseteq &[-30pt] (X_\infty,T) & \arrow[swap]{l}{\tau^*} (X^*_\infty,T) &[-30pt] \supseteq \Omega_{\tau^*}
\end{tikzcd}
\]

Let $\Omega_{\pi^*} \subseteq X^*$ be the set of points that have saturated orbits with respect to $\pi^*$ (a residual set) and which are multiply recurrent under $\{T^n \ | \ n \in \N\}$ (a residual set by \cref{thm_furst_weiss}); $\Omega_{\tau^*} \subseteq X_\infty^*$ be the residual set of one-point fibers for $\tau^*$; and $\Omega_\infty \subseteq X_\infty$ be the residual set guaranteed by \cref{theorem_ae_return_set_in_minimal_nil_is_mult_thick} for the system $(X_\infty,T)$.  Define
\[\Omega \defeq \pi^{-1} \Omega_\infty \cap \sigma^{*} \Omega_{\pi^*} \cap \sigma^{*} (\pi^{*})^{-1} \Omega_{\tau^*}.\]
It follows by the minimality of all the systems involved and \cref{lem_factor_maps_pass_residuality} that $\Omega$ is a residual subset of $X$.  Using the fact that $T: X \to X$ is a homeomorphism, by replacing $\Omega$ with $\cap_{t \in \Z} T^t \Omega$, we may assume that $\Omega$ is $T$-invariant.

Let $\delta > 0$.  According to \cref{lem_factor_maps_pass_residuality} \ref{item_factor_map_fact_four}, there exists $\delta' > 0$ and a $\delta$-dense set $Y \subseteq X$ such that for all $y \in Y$,
\begin{align}
\label{eqn_ball_inclusion_in_thm_b}
    \pi B(y,\delta) \supseteq B(\pi y,\delta').
\end{align}
Let $C, N \in \N$ and $X_\infty = X_{\infty, 1} \cup \cdots \cup X_{\infty, C}$ be the $T$-adapted, clopen partition as guaranteed by \cref{theorem_ae_return_set_in_minimal_nil_is_mult_thick} for the system $(X_\infty,T)$ and with $\delta'$ as $\delta$; note that if $(X,T)$ is totally minimal, then so is $(X_\infty,T)$, whence $C$ can be taken to be equal to 1.  Define $X_i = \pi^{-1} X_{\infty, i}$ so that $X = X_1 \cup \cdots \cup X_{C}$ is a $T$-adapted, clopen partition of $X$, and note that for all $x \in X$, $\compc(x) = \compc(\pi x)$.

Let $x \in \Omega$, $y \in Y$, and $I \in \N$ satisfy $I \equiv \compc(y) - \compc(x) \pmod C$.  We aim to show that the set $R(x,B(y, \delta))$ is multiplicatively thick in $I\nn$.

Let $x^* \in \Omega_{\pi^*}$ be such that $\sigma^* x^* = x$.  By the definition of return-time set,
\[R(x,B(y, \delta)) = R\big(x^*,(\sigma^*)^{-1} B(y, \delta) \big).\]
Since $x^*$ is multiply recurrent under $\{T^n \ | \ n \in \N\}$ and $(X_\infty^*,T)$ is a $d$-step topological characteristic factor for $(X^*,T)$ for all $d \geq 2$, to show that the set $R\big(x^*,(\sigma^*)^{-1} B(y, \delta) \big)$ is multiplicatively thick in $I\nn$, it suffices by \cref{lem_top_char_for_thickness_of_return_sets} \ref{item_lem_top_factor_reformulation_v2_two} to show that the set $R\big(\pi^* x^*,\pi^* (\sigma^*)^{-1} B(y, \delta) \big)$ is multiplicatively thick in $I\zn$.

Appealing to \cref{lem_return_time_property_of_almost_onetoone_extensions} (which holds since $\pi^* x^* \in \Omega_{\tau^*}$), the fact that $\tau^* \circ \pi^* = \pi \circ \sigma^*$, the fact that $y \in Y$, and \eqref{eqn_ball_inclusion_in_thm_b}, we see
\begin{align*}
    R\big(\pi^* x^*,\pi^* (\sigma^*)^{-1} B(y, \delta) \big) &= R\big(\tau^* \pi^* x^*, \tau^*\pi^* (\sigma^*)^{-1} B(y, \delta) \big)\\
    &= R\big( \pi x, \pi B(y, \delta) \big)\\
    &\supseteq R\big( \pi x, B(\pi y,\delta') \big).
\end{align*}
Since $\pi x \in \Omega_\infty$, $\pi y \in X_\infty$, and $I \equiv \compc(\pi y) - \compc(\pi x) \pmod C$, it follows from \cref{theorem_ae_return_set_in_minimal_nil_is_mult_thick} that the set $R\big( \pi x, B(\pi y, \delta') \big)$ is multiplicatively thick in $I\zn$.  This shows that the set $R(x, B(y, \delta))$ is multiplicatively thick in $I\nn$ and finishes the proof of the theorem.
\end{proof}

Using the notion of a \emph{topological natural extension} (see, eg., \cite[Def. 2.8, Lem. 2.9]{glasscock_koutsogiannis_richter_2019}), we can prove an analogue of \cref{theorem_advanced_main_theorem_on_returntime_sets} for non-invertible systems.  Recall that for non-invertible systems, the return-time set $R(x,U)$ is a subset of $\N$.

\begin{theorem}
\label{thm_non_invert_main_result}
Let $(X,T)$ be a minimal system.  There exists a residual, $T$-invariant set $\Omega \subseteq X$ for which the following holds.  For all $\delta > 0$, there exists $C, N \in \N$, a clopen, $T$-adapted partition $X = X_1 \cup \cdots \cup X_C$, and a $\delta$-dense set $Y \subseteq X$ such that for all $x \in \Omega$, $y \in Y$, and $I \in \N$ satisfying $I \equiv \compc(y) - \compc(x) \pmod C$, the set $R\big(x,B(y,\delta) \big)$ is multiplicatively thick in $I\nn$.  If $(X,T)$ is totally minimal, then $C$ can be taken to be equal to $1$.
\end{theorem}

\begin{proof}
Let $\pi: (W,T) \to (X,T)$ be the topological natural extension of $(X,T)$: the space $W \subseteq X^\Z$ is the set of orbits $(w_i)_{i \in \Z}$, $w_{i+1} = Tw_i$; the map $T$ is the left-shift on $W$; and $\pi$ is the projection onto the zero coordinate.  Let $\Omega_W \subseteq W$ be the residual, $T$-invariant set guaranteed by \cref{theorem_advanced_main_theorem_on_returntime_sets}.  Define $\Omega_X \defeq \pi \Omega_W$, and note that $\Omega_X$ is $T$-invariant and, by \cref{lem_factor_maps_pass_residuality}, residual.

Let $\delta > 0$.  Choose $\delta' > 0$ such that the image of a $\delta'$-dense subset of $W$ under $\pi$ is $\delta$-dense in $X$ and so that for all $z \in W$, $\pi B(z,\delta') \subseteq B(\pi z, \delta).$  Let $C, N \in \N$, $W = W_1 \cup \cdots \cup W_C$ be the clopen, $T$-adapted partition, and $Z \subseteq W$ be the $\delta'$-dense set guaranteed by \cref{theorem_advanced_main_theorem_on_returntime_sets} for the system $(W,T)$ with $\delta'$ as $\delta$.

Define $X_i = \pi W_i$ and $Y = \pi Z$.  We claim that $X = X_1 \cup \cdots \cup X_C$ is a clopen, $T$-adapted partition of $X$ and that $Y$ is $\delta$-dense.  That $Y$ is $\delta$-dense follows from our choice of $\delta'$.  To see the first claim, it suffices to show that the sets $\pi W_1$, \dots, $\pi W_C$ are disjoint.  If $w \defeq (w_i)_{i \in \Z}$ and $v \defeq (v_i)_{i \in \Z}$ are elements of $W$ for which $\pi w = \pi v$, then for all $n \in \Nz$, $w_n = v_n$, whereby the distance between $T^{Cn} w$ and $T^{Cn} v$ approaches zero as $n \to \infty$.  Since $\compc(T^{Cn}w) = \compc(w)$ and $\compc(T^{Cn}v) = \compc(v)$, it follows that $w$ and $v$ belong to the same component of the partition of $W$, and hence that the sets $\pi W_1$, \dots, $\pi W_C$ are disjoint.

Let $x \in \Omega$, $y \in Y$, and $I \in \N$ satisfy $I \equiv \compc(y) - \compc(x) \pmod C$.  Let $w \in \Omega_W$ and $z \in Z$ be such that $\pi w = x$ and $\pi z = y$.    It follows from \cref{theorem_advanced_main_theorem_on_returntime_sets} that the set $ R(w,B(z,\delta'))$ is multiplicatively thick in $I\nn$.  Since
\[R(w,B(z,\delta')) \cap \N \subseteq R(\pi w,\pi B(z,\delta')) \cap \N \subseteq R(x,B(y,\delta)),\]
we have that the set $R(x,B(y,\delta))$ is multiplicatively thick in $I\nn$, as desired.
\end{proof}

We can finally prove \cref{maintheorem_ae_return_set_is_mult_thick}.

\begin{proof}[Proof of \cref{maintheorem_ae_return_set_is_mult_thick}]
Let $\Omega$ be as in \cref{thm_non_invert_main_result}.  Let $U \subseteq X$ be non-empty and open.  Choose $\delta > 0$ to be such that for all $\delta$-dense sets $Y \subseteq X$, there exists $y \in Y$ such that $B(y,\delta) \subseteq U$.  Let $C, N \in \N$, $X = X_1 \cup \cdots \cup X_C$, and $Y \subseteq X$ be as guaranteed by \cref{thm_non_invert_main_result}.  Tracing the origin of the parameters $C$ and $N$ all the way back to \cref{thm_main_fullversion}, we find that $C \mid N$.  Thus, by the choice of $\delta$ and because $\{1, \ldots, N\}$ contains a complete set of residues modulo $C$, we see that for all $x \in \Omega$, there exists $I \in \{1, \ldots, N\}$ such that the set $R(x,U)$ is multiplicatively thick in $I\nn$.
\end{proof}

The improvements to \cref{maintheorem_ae_return_set_is_mult_thick} mentioned in the introduction follow immediately from \cref{thm_non_invert_main_result}.  We call attention in particular to the one regarding the translated return-time set $R(x,U) - t = R(T^t x,U)$.  Since the set $\Omega$ in \cref{thm_non_invert_main_result} is $T$-invariant, it follows that for all $t \in \N$, the set $R(x,U) - t$ is multiplicatively thick in $[I-t]\nn$, where $I \in \{1, \ldots, C\}$ is as guaranteed for the set $R(x,U)$ and $[I-t]$ denotes the least positive modulo $C$ residue of $I-t$.  It follows, in particular, that all translates of $R(x,U)$ are multiplicatively thick in at least one of only finitely many cosets of a particular congruence subsemigroup.  It is not known whether all nil-Bohr sets satisfy this property; see \cref{quest_finer_nil_uniformity_questions} in \cref{sec_improvements_to_main_results}.

\subsection{Enhancements to van der Waerden's theorem}
\label{sec_refinement_on_vdw}

In this section, we demonstrate the equivalence between Theorems \ref{maintheorem_ae_return_set_is_mult_thick} and \ref{thm_comb_equivalent_to_main_thm}, and we derive the combinatorial refinements of van der Waerden's theorem in \cref{main_vdw_type_theorem}.  Note that \cref{maintheorem_ae_return_set_is_mult_thick} is proved in \cref{sec_proof_of_return_sets_mult_thick}, so the following constitutes a proof of \cref{thm_comb_equivalent_to_main_thm}.

\begin{proof}[Proof of the equivalence of Theorems \ref{maintheorem_ae_return_set_is_mult_thick} and \ref{thm_comb_equivalent_to_main_thm}]

Let $(X,T)$ be a minimal system.  To see that \cref{maintheorem_ae_return_set_is_mult_thick} implies \cref{thm_comb_equivalent_to_main_thm}, let $\Omega \subseteq X$ be as guaranteed by \cref{maintheorem_ae_return_set_is_mult_thick}, and let $U \subseteq X$ be non-empty, open.  Let $N \in \N$ be as guaranteed by \cref{maintheorem_ae_return_set_is_mult_thick}.  Let $V \subseteq X$ be non-empty, open, and let $x \in \Omega \cap V$.  By \cref{maintheorem_ae_return_set_is_mult_thick}, there exists $I \in \{1, \ldots, N\}$ such that the set $R(x,U)$ is multiplicatively thick in $I\nn$.  Let $n_1, \ldots, n_k \in I\nn$.  There exists $m \in \nn$ such that $m\{n_1, \ldots, n_k\} \subseteq R(x,U)$.  It follows that $x$ is a member of the set in \eqref{eqn_set_eqn_in_comb_equiv_thm}, proving that the set is non-empty, as desired.

To see that \cref{thm_comb_equivalent_to_main_thm} implies \cref{maintheorem_ae_return_set_is_mult_thick}, we will first define the residual, $T$-invariant set $\Omega \subseteq X$.  Let $\mathcal{U}$ be a countable collection of non-empty, open sets that form a base for the topology of $X$.  For $U \in \mathcal{U}$, let $N = N(U) \in \N$ be as guaranteed by \cref{thm_comb_equivalent_to_main_thm}.  For all $k \in \N$, $n_1, \ldots, n_k \in \nn$, and non-empty, open sets $V \subseteq X$, define
\[V' = V'(n_1, \ldots, n_k, V) \defeq \bigcup_{\substack{ m \in \nn \\ I \in \{1, \ldots, N\} }} \big( V \cap T^{-In_1m}U \cap \cdots \cap T^{-In_km}U \big).\]
We claim that $V'$ is an open, dense subset of $V$.  It is open and contained in $V$ by definition, and it is non-empty by \cref{thm_comb_equivalent_to_main_thm}.  To see that it is dense in $V$, we will show that $V' \cap W$ is non-empty for every non-empty, open $W \subseteq V$.  Let $W \subseteq V$ be non-empty, open.  By the same reasoning as before, the set $W'$ is an open, non-empty subset of $W$.  Since $W' \subseteq V' \cap W$, the set $V' \cap W$ is non-empty.

Define
\[\Omega \defeq \bigcap_{t \in \Nz} T^{-t} \bigcap_{U \in \mathcal{U}} \bigcap_{\substack{k \in \N \\ n_1, \ldots, n_k \in \nn}} X'.\]
Since $X'$ is an open, dense subset of $X$, the set $\Omega$ is residual and $T$-invariant.  This is the set $\Omega$ for which we intend to verify the statement in \cref{maintheorem_ae_return_set_is_mult_thick}.

Let $U \subseteq X$ be non-empty, open.  Without loss of generality, by replacing $U$ with a member of $\mathcal{U}$ that it contains, we may assume that $U \in \mathcal{U}$.  Let $N = N(U) \in \N$ be as guaranteed by \cref{thm_comb_equivalent_to_main_thm}.  Let $x \in \Omega$.  To finish the proof, we must show that there exists $I \in \{1, \ldots, N\}$ such that for all $n_1, \ldots, n_k \in \nn$, there exists $m \in \nn$ such that 
\begin{align}
    \label{eqn_to_show_for_thickness}
    I\{n_1, \ldots, n_k\} m \subseteq R(x,U).
\end{align}
For this, we claim that it suffices to show that for all $n_1, \ldots, n_k \in \nn$, there exists $I \in \{1, \ldots, N\}$ and $m \in \nn$ such that \eqref{eqn_to_show_for_thickness} holds.  Indeed, since $\nn$ is countable and the function $(n_1, \ldots, n_k) \mapsto I$ has finite range, the order of the quantifiers on $n_1, \ldots, n_k$ and $I$ can be exchanged.

Let $n_1, \ldots, n_k \in \nn$.  Since $x \in \Omega$, we have that $x \in X'$, where we recall that $X'$ depends on $U$, $N$, $X$, and $n_1, \ldots, n_k$.  Therefore, there exists $I \in \{1, \ldots, N\}$ and $m \in \nn$ such that $x \in T^{-In_1m}U \cap \cdots \cap T^{-In_km}U$.  This implies that \eqref{eqn_to_show_for_thickness} holds, as was to be shown.
\end{proof}

Now we combine \cref{thm_comb_equivalent_to_main_thm} with a standard topological correspondence principle to deduce \cref{main_vdw_type_theorem}.

\begin{proof}[Proof of \cref{main_vdw_type_theorem}]
Let $Z = \{0,1\}^{\Nz}$ endowed with the product topology.  Denote the $n^{\text{th}}$ coordinate of $z \in Z$ by $z(n)$.  Let $T: Z \to Z$ be the left shift, so that $(Tz)(n) = z(n+1)$.  Let $C_0 = \{z \in Z \ | \ z(0) = 1\}$, and note that $T^n z \in C_0$ if and only if $z(n) = 1$.

Let $\one_A \in Z$ be the indicator function of the set $A$.  Because $A$ is piecewise syndetic, the orbit closure $\overline{T^{\Nz} \one_A}$ contains a closed set, $T$-invariant set $X$ for which the system $(X,T)$ is minimal and with the property that the set $U \defeq X \cap C_0$ is non-empty. Indeed, note first that $\overline{T^{\Nz} \one_A}$ contains a point of the form $\one_S$, where $S \subseteq \Nz$ is syndetic.  By a theorem of Auslander and Ellis (cf. \cite[Thm. 8.7]{furstenberg_book_1981}), the point $\one_S$ is proximal to a uniformly recurrent point $x \in Z$.\footnote{In a system $(X,T)$, points $x$ and $y$ are \emph{proximal} if $\inf_{n \in \N} d_X(T^n x, T^n y) = 0$, and the point $x$ is \emph{uniformly recurrent} if for all open neighborhoods $U$ of $x$, the return-time set $R(x,U)$ is additively syndetic.}  Define $X = \overline{T^{\Nz} x}$.  Note that $X \subseteq \overline{T^{\Nz} \one_S} \subseteq \overline{T^{\Nz} \one_A}$ and that $(X,T)$ is minimal because $x$ is uniformly recurrent.  Since $S$ is syndetic and $x$ and $\one_S$ are proximal, the set $\{n \in \Nz \ | \ x(n) = 1\}$ is syndetic, and hence $X \cap C_0$ is non-empty.

Let $N \in \N$ be as guaranteed by \cref{thm_comb_equivalent_to_main_thm} for the system $(X,T)$ and the set $U$.  Let $B \subseteq \N$ be syndetic, and define $Y = \overline{T^{B} \one_A}$.  We claim that $Y \cap X$ has non-empty interior in $X$.  Indeed, since $B$ is syndetic, there exists $\ell \in \N$ such that $\cup_{i=0}^\ell (B-i) \supseteq \Nz$.  It follows that $\cup_{i=0}^\ell T^{-i} Y = Z$, and hence by the Baire Category Theorem that there exists $i \in \{0, \ldots, \ell\}$ such that the set $T^{-i} Y \cap X$ has non-empty interior in $X$.  Since minimal maps are semiopen \cite[Thm. 2.4]{kolyada_snoha_trofimchuk_2001}, by applying the map $T^i$ and using that $T^i X \subseteq X$, we see that $Y \cap X$ has non-empty interior in $X$.

Let $V \subseteq X$ be the interior of $Y \cap X$, and let $I \in \{1, \ldots, N\}$ be as guaranteed by \cref{thm_comb_equivalent_to_main_thm}. Let $n_1, \ldots, n_k \in I \nn$.  By \cref{thm_comb_equivalent_to_main_thm}, there exists $m \in \nn$ such that \eqref{eqn_set_eqn_in_comb_equiv_thm} holds.  It follows that
\[\overline{T^{B} \one_A} \cap T^{-n_1m} C_0 \cap \cdots \cap T^{-n_km} C_0 \neq \emptyset.\]
Since $C_0$ and its pre-images under $T$ are open, there exists $b \in B$ such that $T^b \one_A \in \cap_{i=1}^k T^{-n_im} C_0$.  Therefore, for all $i \in \{1, \ldots, k\}$, $T^{b+ n_im} \one_A \in C_0$, meaning $b+n_im \in A$.  This demonstrates the non-emptiness of the set in \eqref{eqn_set_intersection_in_combo_corollary}.

Let $k \in \N$, and let $D \subseteq \N$ be the set of $\nn$-starters described in \eqref{eqn_nn_starters}.  To prove that $D$ is additively thick, we will show that it has non-empty intersection with all additively syndetic subsets of $\N$.  Let $B \subseteq \N$ be syndetic, and let $I \in \{1, \ldots, N\}$ be as guaranteed by the work above.  For $i \in \{1, \ldots, k+1\}$, define $n_i = I((i-1)N+1)$.  The work above guarantees that there exists $m \in \nn$ such that \eqref{eqn_set_intersection_in_combo_corollary} holds.  Note that $\cap_{i=1}^{k+1} (A-n_im) \subseteq D$ (with the ``$m$'' in the set-builder notation equal to $Im$), whereby it follows from \eqref{eqn_set_intersection_in_combo_corollary} that $B \cap D \neq \emptyset$, as was to be shown.
\end{proof}

\section{Related results and open ends}
\label{sec_final_section}

There are a number of natural directions and open questions suggested by the main results in this paper.  We collect some of those questions and directions in this section in an effort to spur further investigation.

\subsection{Dense simultaneous approximation in the typical skew product}
\label{sec_dense_approx_in_skew_prods}

Denote by $C(\T,\R)$ the space of continuous maps from the 1-torus to the real numbers equipped with the uniform norm, $\| \cdot \|_\infty$. For $\alpha \in \T$ and $h \in C(\T,\R)$, define the skew product map $T_h: \T^2 \to \T^2$ by
\[T_h(x,y) = \big( x + \alpha, y + \pi(h(x)) \big),\]
where $\pi: \R \to \T$ is the quotient map. For $n \in \Z$, define \[h_n(x) = \begin{cases} \sum_{i=0}^{n-1} h(x + i \alpha) & \text{ if $n > 0$} \\ 0 & \text{ if $n=0$} \\ - \sum_{i=n}^{-1} h(x + i \alpha) & \text{ if $n < 0$} \end{cases}\]
so that $T_h^n(x,y) = \big( x + n \alpha, y + \pi(h_n(x)) \big)$.

Skew product systems are distal but not, in general, nilsystems.  Nevertheless, the main result in this section, \cref{theorem_simult_approx_is_typical}, shows that at least some points in skew product systems satisfy the dense simultaneous approximation property.

\begin{theorem}
\label{theorem_simult_approx_is_typical}
Let $\alpha \in \T \setminus \Q$. For all $(x,y) \in \T^2$, for almost every $h \in C(\T,\R)$,
\begin{enumerate}[label=(\Roman*)]
    \item the point $(x,y)$ exhibits the dense simultaneous approximation property in the system $(\T^2, T_h)$, and
    \item for all non-empty, open $U \subseteq \T^2$, the return-time set $R((x,y),U)$ is multiplicatively thick in a congruence subsemigroup.
\end{enumerate}
\end{theorem}

\begin{proof}
Let $(x,y) \in \T^2$.  For $\delta > 0$, $N \in \N$, $n \in \Z$, $w \in \T^2$, $\eps > 0$, and $F \subseteq \Z$, define
\[H = H(\delta, N, n, w, \eps, F) = \big\{ h \in C(\T,\R) \ \big| \ \exists m \in \zn, \ mF \subseteq R\big((x,y),B(w,\eps)\big) \big\}.\]
We will show that for all $\delta > 0$, there exists $N \in \N$ such that for all $n \in \Z$ coprime to $N$, there exists a $\delta$-dense set $W \subseteq \T^2$ such that for all $w \in W$, $\eps > 0$, and finite $F \subseteq N \N + n$, the set $H$ is open and dense in $C(\T,\R)$.  It will follow that the set
\[\bigcap_{e \in \N} \ \bigcup_{N} \ \bigcap_{n} \ \bigcup_{W} \ \bigcap_{w, \eps, F} \ H(1/e, N, n, w, \eps, F)\]
is a residual subset of $C(\T,\R)$. Unpacking the definitions, any element $h$ of this set is such that the point $(x,y)$ exhibits the dense simultaneous approximation property in the system $(\T^2, T_h)$.  Moreover, for any non-empty, open set $U \subseteq \T^2$, choosing $\delta > 0$ sufficiently small, then $w \in U$, then $\eps > 0$ sufficiently small shows that the return-time set $R((x,y),U)$ is multiplicatively thick in $\zn$, proving the theorem.

Let $\delta  >0$.  Choose $N \in \N$ such that $N \geq \delta^{-1}$.  Let $n \in \Z$ be coprime to $N$.  Put $W = \{0, 1/N, \ldots, (N-1) / N\}^2 + (x,y) \subseteq \T^2$.  Let $w = (w_1,w_2) \in W$, $\eps > 0$, and $F \subseteq N\Z + n$ be finite.

The set $H$ is open because for all $m, f \in \Z$, the map $h \mapsto T_h^{mf}(x,y)$ is continuous. To see that $H$ is dense in $C(\T,\R)$, fix $g \in C(\T,\R)$ and $\sigma > 0$.  We will show that there is $h \in H$ such that $\|g - h\|_\infty < \sigma$.  The idea is to choose $m \in \zn$ large and make small adjustments (in the uniform norm) to $g$ in order to ensure that the points $T_h^{mf}(x,y)$, $f \in F$, approximate $w$.  Since we endow finite Cartesian product spaces with the maximum metric, note that $mF \subseteq R\big((x,y),B(w,\eps)\big)$ if and only if for all $f \in F$, $d_{\T}(mf \alpha, w_1 - x) < \eps$ and $d_{\T}\big(\pi(h_{mf}(x)), w_2 - y \big) < \eps$.

Since $\alpha$ is irrational, the set
\[\big\{ m \in \Z \ \big| \ \forall f \in F, \ d_{\T}(mf \alpha, w_1 - x) < \eps \big\}\]
is a totally minimal Bohr set.  Choose $m > \sigma^{-1}$ from the intersection of this set with $\zn$. Write $F = \{f_k^- < \cdots < f_1^- < 0 < f_1^+ < \cdots < f_\ell^+\}$, and denote by $J$ the interval $[mf_k^-, mf_\ell^+] \cap \Z$. Note that for any finite sequence $(\eta_i)_{i \in J} \subseteq (-\sigma, \sigma)$, there exists $h \in C(\T,\R)$ with $\|g - h\|_\infty < \sigma$ such that for all $i \in J$,
\begin{align}
\label{eqn_property_for_small_pert_of_g}
    h(x + i \alpha) = g( x + i \alpha) + \eta_i.
\end{align}
Indeed, the points $x + i \alpha$, $i \in J$, are distinct, so one needs only to make small continuous adjustments to $g$ to arrive at $h$.

We will show that there exists a sequence $(\eta_i)_{i \in J} \subseteq (-\sigma, \sigma)$ so that any function $h$ with the property described in \eqref{eqn_property_for_small_pert_of_g} belongs to $H$. Note first that
\[    h_{m f_1^-}(x) = g_{m f_1^-}(x) - \sum_{i=mf_1^-}^{-1} \eta_i \quad \text{ and } \quad h_{m f_1^+}(x) = g_{m f_1^+}(x) + \sum_{i=0}^{mf_1^+ - 1} \eta_i.\]
Since $m > \sigma^{-1}$, we have that $|mf_1^\pm| \sigma > 1$, so there is a choice of $\eta_{mf_1^-}, \ldots, \eta_{m f_1^+ - 1} \in (-\sigma , \sigma)$ so that $d_{\T}\big(\pi(h_{m f_1^-}(x)), w_2 - y\big) < \eps$ and $d_{\T}\big(\pi(h_{m f_1^+}(x)), w_2 - y\big) < \eps$. Then
\begin{align*}
h_{m f_2^-}(x) &= h_{m f_1^-}(x) + g_{m(f_2^--f_1^-)}(x + mf_1^- \alpha) - \sum_{i=mf_2^-}^{mf_1^- - 1} \eta_i \\
h_{m f_2^+}(x) &= h_{m f_1^+}(x) + g_{m(f_2^+-f_1^+)}(x + mf_1^+ \alpha) + \sum_{i=mf_1^+}^{mf_2^+ - 1} \eta_i
\end{align*}
so there exists a choice of $\eta_{-mf_2^-}, \ldots, \eta_{m f_1^- - 1} \in (-\sigma , \sigma)$ so that $d_{\T}\big(\pi(h_{m f_2^-}(x)), w_2 - y\big) < \eps$ and a choice of $\eta_{mf_1^+}, \ldots, \eta_{m f_2^+ - 1} \in (-\sigma , \sigma)$ so that $d_{\T}\big(\pi(h_{m f_2^+}(x)), w_2 - y\big) < \eps$.  We continue inductively to define the rest of the sequence $(\eta_i)_{i \in J}$.  A function $h \in C(\T,\R)$ with the property described in \eqref{eqn_property_for_small_pert_of_g} will then be a member of $H$ and will be $\sigma$-close to $g$.
\end{proof}

The result in \cref{theorem_simult_approx_is_typical} is surely not optimal.  We invite the reader to address any of the myriad questions raised by it.

\begin{question}
\label{quest_extensions_of_skew_product_example} \leavevmode
\begin{enumerate}[label=(\Roman*)]
    \item Can the quantifiers on $(x,y)$ and $h$ be swapped and improved so that all or almost all systems $(\T^2, T_h)$ are such that all or almost all points $(x,y) \in \T^2$ exhibit the dense simultaneous approximation property?
    \item Do all or almost all systems $(\T^2, T_h)$ satisfy the a.e. dense simultaneous approximation property?
    \item Can a similar result be achieved in the class of $C^1$-skewing functions $h$?  (The technique employed in \cref{theorem_simult_approx_is_typical} seems likely to fail when smoothness is required.)
    \item Can similar results be achieved for higher-order skew product systems?
\end{enumerate}
\end{question}

\subsection{An additively syndetic set that is nowhere multiplicatively thick}
\label{sec_add_synd_not_mult_thick}

In a minimal system, return-time sets are additively syndetic: a union of finitely many translates of $R(x,U)$ covers $\N$.  By \cref{maintheorem_ae_return_set_is_mult_thick}, almost all return-time sets are multiplicatively thick in a coset of a congruence subsemigroup.  It is natural to wonder, then, whether all additively syndetic sets are multiplicatively thick in a coset of a congruence subsemigroup.  This is not the case, as we demonstrate in this section.

In \cref{theorem_add_syndetic_not_mult_thick}, we construct a set $A \subseteq \N$ which is additively syndetic in $\N$ but not multiplicatively thick in any coset of any congruence subsemigroup of $\N$. By \cref{lemma_mult_thick_from_z_to_n}, this yields a similar example in $\Z$: the set $A \cup (-\N)$ is additively syndetic in $\Z$ and not multiplicatively thick in any congruence subsemigroup of $\Z$.

\begin{definition}
\label{def_mult_den_and_synd}
Let $I, N \in \N$ and $A \subseteq \N$.  The \emph{multiplicative upper (Banach) density in $\nn$} of the set $A$ is the quantity $\sup\big\{ \alpha \geq 0 \ \big| \ \forall \text{ finite } F \subseteq \nn, \ \exists m \in \nn, \ |mF \cap A| \geq \alpha |F| \big\}$.  (The terminology can be justified by showing that this is equal to a supremum of the values given to $\one_A$ by the dilation-invariant means on $\nn$; see \cite[Sec. 3]{bergelson_glasscock_2020}.)  The set $A$ has \emph{full multiplicative upper density in $\nn$} if the set $\N \setminus A$ has zero multiplicative upper density in $\nn$.
\end{definition}

We will make use of the following auxiliary results in the construction below.  Recall the definition of multiplicatively syndetic in \cref{def_mult_synd}.
\begin{enumerate}[label=(\Roman*)]
    \item \label{aux_result_two} If $A$ is multiplicatively syndetic in $I \nn$ and a set $D \subseteq \N$ has full multiplicative upper density in $I\nn$, then the set $A \cap D$ is multiplicatively syndetic in $I \nn$.  Indeed, there exist $n_1, \ldots, n_k \in \nn$ for which $\cup_{i=1}^k A/(In_i) \supseteq \nn$.  Thus,
    \[\bigcup_{i=1}^k \frac{A \cap D}{In_i} \ \supseteq \ \bigcup_{i=1}^k \frac{A}{In_i} \  \cap \ \bigcap_{i=1}^k \frac{D}{In_i} \ \supseteq \ \nn \cap \bigcap_{i=1}^k \frac{D}{In_i}.\]
    The set $\cap_{i=1}^k D/(In_i)$ is of full multiplicative upper density in $\nn$, hence multiplicatively syndetic in $\nn$.  Thus, there exist $m_1, \ldots, m_\ell \in \nn$ such that $\cup_{i=1}^k \cup_{j=1}^\ell (A \cap D) / (I n_i m_j) \supseteq \nn$, as desired.
    \item \label{aux_result_three} The set $4^{\N} - 1 \defeq \{4^n - 1 \ | \ n \in \N\}$ has zero multiplicative upper density in all cosets of all congruence subsemigroups; that is, the set $I\nn \setminus (4^{\N} - 1)$ is of full multiplicative upper density in $I \nn$.  This is straightforward to show by hand and is left as an exercise left to the reader.
\end{enumerate}

\begin{theorem}
\label{theorem_add_syndetic_not_mult_thick}
There exists a set $A \subseteq \N$ satisfying $A \cup (A-1) = \N$ that is not multiplicatively thick in any coset of any congruence subsemigroup.
\end{theorem}

\begin{proof}
Define
\begin{align*}
    R_0 &= \bigcup_{n = 0}^\infty \big[2^{2n},2^{2n+1}\big) & \quad R_1 &= \bigcup_{n = 0}^\infty \big[2^{2n+1},2^{2n+2}\big)\\
    J_0 &= R_0 \cap 2\N & \quad J_1 &= R_1 \cap (2\N - 1)\\
    B &= J_0 \cup J_1 & \quad B' &= B \big\setminus \big(4^{\N} - 1 \big). 
\end{align*}
Finally, define $A = \N \setminus B'$.  We will verify that $A \cup (A-1) = \N$ and that $A$ is not multiplicatively thick in any coset of any congruence subsemigroup by showing that $B' \cap (B' - 1) = \emptyset$ and that $B'$ is multiplicatively syndetic in all cosets of all congruence subsemigroups.  This suffices by fact \ref{aux_synd_result_one} in \cref{sec_disjointness}.

That $B' \cap (B' - 1) = \emptyset$ follows immediately from the definition of $B'$ and the fact that $B \cap (B - 1) \subseteq 4^{\N} - 1$. To show that $B'$ is multiplicatively syndetic in all cosets of all congruence subsemigroups, by auxiliary results \ref{aux_result_two} and \ref{aux_result_three} above, it suffices to show that $B$ is multiplicatively syndetic in all cosets of all congruence subsemigroups.  We will show in fact that for all $I \nn$, at least one of the sets $J_0$ or $J_1$ is multiplicatively syndetic in $I \nn$.

Let $I, N \in \N$. We will show that:
\begin{enumerate}[label=(\Roman*)] \setcounter{enumi}{2}
    \item \label{item_either_even_or_odd_zero_mult_den} either $2\N$ or $2\N - 1$ has full multiplicative upper density in $I\nn$; and
    \item \label{item_both_r1_r2_mult_synd} both $R_0 \cap \N$ and $R_1 \cap \N$ are multiplicatively syndetic in $I \nn$.
\end{enumerate}
Taken together, these facts imply that either $J_0$ or $J_1$ is multiplicatively syndetic in $I \nn$.  Indeed, if $2\N$ has full multiplicative upper density in $I\nn$, then by auxiliary result \ref{aux_result_two} above, the set $J_0$ is multiplicatively syndetic in $I \nn$.  If $2\N - 1$ has full multiplicative upper density in $I \nn$, then $J_1$ is multiplicatively syndetic in $I \nn$.

To see \ref{item_either_even_or_odd_zero_mult_den}, we consider three cases: a) $\nn \cap 2\N \neq \emptyset$; b) $I \in 2\N$; and c) $I \nn \subseteq 2\N - 1$.  If a) holds, let $e \in \nn \cap 2\N$.  The set $Ie\nn$ is contained entirely in $2\N \cap I \nn$ and has full multiplicative upper density in $I\nn$.  If b) holds, then $I \nn \subseteq 2\N$, in which case $2\N$ has full multiplicative upper density in $I \nn$.  If c) holds, then $2\N - 1$ has full multiplicative upper density in $I\nn$.

To see \ref{item_both_r1_r2_mult_synd}, let $L \in \N$ be so large that neither $R_1$ nor $R_2$ contain a real-valued, length $L$ arithmetic progression starting in the interval $[I,2I)$ of step size contained in the interval $[NI,2NI)$.  We will show that $\cup_{\ell = 0}^L R_0 / I(\ell N + 1) \supseteq \nn$ and $\cup_{\ell = 0}^L R_1 / I(\ell N + 1) \supseteq \nn$, from which \ref{item_both_r1_r2_mult_synd} follows. 

Let $m \in \nn$; we must show that the set $I m \{1, N + 1, \ldots, LN + 1\}$ intersects both $R_0$ and $R_1$. Let $k$ be such that $2^k \leq m < 2^{k+1}$, and put $c = m / 2^k \in [1,2)$.  Note that if $k$ is even, then $2^k R_1 \subseteq R_1$ and $2^k R_2 \subseteq R_2$, while if $k$ is odd, then $2^k R_1 \subseteq R_2$ and $2^k R_2 \subseteq R_1$. Therefore, to show that the set $2^k c I \{1, N + 1, \ldots, LN + 1\}$ intersects both $R_0$ and $R_1$, it suffices to prove that the set $c I \{1, N + 1, \ldots, LN + 1\}$ intersects both $R_1$ and $R_2$.  But this set is a length $L+1$ arithmetic progression starting in the interval $[I,2I)$ with step size contained in the interval $[NI,2NI)$, so by the choice of $L$, it must intersect both $R_1$ and $R_2$.  This finishes the verification of \ref{item_both_r1_r2_mult_synd} and the proof of the theorem.
\end{proof}

While the set $A$ constructed in \cref{theorem_add_syndetic_not_mult_thick} is nowhere multiplicatively thick, it is still multiplicatively quite large: it is multiplicatively syndetic in all cosets of all congruence subsemigroups.  This begs the question as to just how much such combinatorial constructions can be improved by making $A$ as ``multiplicatively small'' as possible while retaining additive syndeticity; a concrete realization of this question is given in \cref{quest_mult_piecewise_synd} in the next section.

\subsection{Open ends}
\label{sec_open_ends}

There are several ways in which the main results may be improved upon within the scope of this work.  We begin by highlighting those, then we move toward the more speculative open ends that the main results suggest.

\subsubsection{Improvements to the main results}
\label{sec_improvements_to_main_results}

The first question pertains to the degree to which the modulus $N$ in the dense simultaneous approximation property (\cref{def_dense_simult_approx}) can be chosen to be independent of the approximating point in a nilsystem.

\begin{question}
\label{quest_all_x_have_simult_approx}
Let $(X,T)$ be a minimal nilsystem.  Is it true that for all $\delta > 0$, there exists $N \in \N$ such that for all $x \in X$, the set $\sa(x,N\Z + 1)$ is $\delta$-dense in $X$?
\end{question}

In the case that ``for all $x \in X$'' is replaced with ``for almost all $x \in X$'', \cref{maintheorem_nilsystems_have_simult_approx} yields a positive answer: minimal nilsystems possess the a.e. dense simultaneous approximation property. The residual set in that theorem originates with the set of vectors in the product nilmanifold $X^{d+1}$ whose $(T^0 \times T^1 \times \cdots \times T^d)$-orbit closure is equal to its $(T^0 \times T^1 \times \cdots \times T^d)$-prolongation class.  This equality allowed us to show that the set of $N$-periodic polynomial orbits was uniformly dense in most $(T^0 \times T^1 \times \cdots \times T^d)$-orbit closures.

In equicontinuous systems (and, hence, in $1$-step nilsystems), orbit closures coincide with prolongation classes, but this is already not the case in $2$-step nilsystems: in $\big(\T^2, T: (x,y) \mapsto (x, y + x) \big)$, for example, we see that $\orbitz{T}{(0,0)} = \{(0,0)\}$ while $\prolong_T((0,0)) = \{0\} \times \T$.  Thus, it seems that using the tools in this paper to give a positive answer to \cref{quest_all_x_have_simult_approx} would require a finer understanding of the denseness of $N$-periodic polynomial orbits in $(T^0 \times T^1 \times \cdots \times T^d)$-orbit closures.

A negative answer to \cref{quest_all_x_have_simult_approx} may be achievable via careful computation in a specific example.  Such an example would be interesting but may leave open the following modification of \cref{quest_all_x_have_simult_approx} and one of its related combinatorial consequences.

\begin{question}
\label{quest_finer_nil_uniformity_questions}
\begin{enumerate}[label=(\Roman*)]
    \item \label{item_quest_nilsystem_orbit_indepent} Let $(X,T)$ be a minimal nilsystem and $x_0 \in X$.  Is it true that for all $\delta > 0$, there exists $N \in \N$ such that for all $t \in \Z$, the set $\sa(T^tx_0,N\Z + 1)$ is $\delta$-dense in $X$?
    \item \label{item_quest_bohr_set_translate_indepent} If $A \subseteq \Z$ is a nil-Bohr set, does there exists $N \in \N$ such that for all $t \in \Z$, the set $A + t$ is multiplicatively thick in a coset of $\zn$?
\end{enumerate}
\end{question}

A second question begged by the main results in this work is whether or not all return-time sets are multiplicatively thick in a coset of a congruence subsemigroup.

\begin{question}
\label{quest_all_returns_mult_thick}
Let $(X,T)$ be a minimal system.  Is it true that for all $x \in X$ and all non-empty, open $U \subseteq X$, the return-time set $R(x,U)$ is multiplicatively thick in a coset of a congruence subsemigroup?  Note that this is equivalent to asking: Does the conclusion in \cref{thm_comb_equivalent_to_main_thm} still hold if the set $V$ is allowed to be a singleton $\{x\}$?
\end{question}

Residual sets seem likely to be an inherent part of the characteristic factor machinery used in \cref{sec_proof_of_return_sets_mult_thick}.  Thus, that machinery is not well-suited to address \cref{quest_all_returns_mult_thick}, at least in the manner that it is used in this work.  Since return-time sets are additively syndetic, one could approach the question from a much broader combinatorial framework; see \cref{subsec_comb_explorations} below. Should some approach yield a positive answer to \cref{quest_all_returns_mult_thick}, a number of follow-up questions along the lines of those in Questions \ref{quest_all_x_have_simult_approx} and \ref{quest_finer_nil_uniformity_questions} suggest themselves: can the modulus $N$ be chosen independently of the point $x$, or for $x_0 \in X$, can $N$ be chosen so that for all $t \in \Z$, the set $R(x_0,U) + t$ is multiplicatively thick in a coset of $\nn$?

By \cref{maintheorem_nilsystems_have_simult_approx}, all points in minimal nilsystems exhibit the dense simultaneous approximation property. By \cref{thm_glasner}, almost every point of a minimal, weakly mixing system exhibits the property.  And, it is demonstrated in \cref{sec_dense_approx_in_skew_prods} that there are distal systems outside the class of nilsystems with points that satisfy the property.  Thus, it is natural to ask about the extent to which this property holds.

\begin{question}
\label{quest_dense_simult_in_all_systems}
Let $(X,T)$ be a minimal system.  Does the dense simultaneous approximation property hold for some point in $X$?  Does it hold for all points in $X$?
\end{question}

At issue in this question is how the modulus $N$ in the dense simultaneous approximation property is affected by the degree of approximation $\eps$.  In a nilsystem, the modulus can be chosen independently of $\eps$ by appealing to the fact that polynomial orbits are determined by finitely many values, \cref{lem_topology_on_poly_space}.  Already in an inverse limit of minimal nilsystems, it is not clear whether or not this independence can be ensured.  A positive answer to \cref{quest_dense_simult_in_all_systems} for inverse limits of nilsystems would seem likely to combine with the topological characteristic factor machinery from \cref{sec_proof_of_return_sets_mult_thick} to yield a positive answer for distal systems (for which the factor map to the infinite-step pro-nilfactor is open) or even all minimal systems.

\subsubsection{Broader simultaneous approximation phenomena}
\label{sec_broader_simult_approx}

Given a collection $\mapspace$ of continuous maps of a compact metric space $X$, we could define $\sa(x,\mapspace)$ to be the set of points that are simultaneously approximated by $x$ under $\mapspace$ according to the definition from the first lines of the introduction.  It is natural to wonder what it is about the set of maps $\mapspace = \{T^i \ | \ i \in N\Z + n\}$ in \cref{maintheorem_nilsystems_have_simult_approx} that suffices for the dense simultaneous approximation property to hold.  Can this property be meaningfully abstracted to collections of commuting nilrotations, or perhaps even more generally?

\begin{question}
\label{quest_more_abstract_conditions_on_mapspace}
Let $(X,T)$ be a minimal nilsystem.  Let $x \in X$ and $\delta > 0$.  What necessary and sufficient conditions can be placed on a collection of commuting nilrotations $\mapspace$ so that the set $\sa(x, \mapspace)$ is $\delta$-dense in $X$?  To what degree can those conditions be made independent of the point $x$ (up to a residual set)?
\end{question}

In a related direction, there seems to be nothing particularly relevant about $\Z$-actions in the definition of simultaneous approximation.  A more general definition could be as follows.  Let $G$ be a group, and suppose $\mapspace$ is a collection of $G$-actions on $(X,d)$ by homeomorphisms.  The point \emph{$x$ simultaneously approximates $y$ under $\mapspace$} if for all finite $F \subseteq \mapspace$ and all $\eps > 0$, there exists $g \in G \setminus \{e_G\}$ such that for all $T \in F$, $d_X(T_g x, y) < \eps$.  Defining $\sa(x,\mapspace)$ as above, the analogue of \cref{quest_more_abstract_conditions_on_mapspace} in this more general setting of $G$-actions makes sense.

This is potentially already interesting in the case of $\Z^2$-nilsystems $(X,T,S)$.  For $A \subseteq \Z^2$, the set $\sa(x,A)$ consists of those points $y \in X$ such that for all $F \subseteq A$ and all $\eps > 0$, there exists $m \in \Z^2$ such that for all $f \in F$, $d_X(T^{m_1f_1}S^{m_2f_2} x, y) < \eps$.

\begin{question}
\label{quest_z2_formulations}
Does every point in a minimal $\Z^2$-nilsystems satisfy the dense simultaneous approximation property?  Do minimal $\Z^2$-nilsystems satisfy the a.e. dense simultaneous approximation property?  Are return-time sets in minimal $\Z^2$-systems multiplicatively thick in a coset of a multiplicative subsemigroup of $(\Z^2 \setminus \{(0,0)\},\cdot)$?
\end{question}

To answer these questions, it must first be decided what the ``dense simultaneous approximation property'' should mean in a $\Z^2$-nilsystem $(X,T,S)$ and what the analogue of a ``congruence subsemigroup'' should be in $(\Z^2 \setminus \{(0,0)\},\cdot)$ (where $(n,m) \cdot (\ell, k) = (n\ell,mk)$).  Developing results parallel to those in \cref{section_simul_approx_for_irrational_rotations} would likely be useful for that purpose.

\subsubsection{Combinatorial approaches}
\label{subsec_comb_explorations}

As discussed in the introduction, the precise extent to which additively syndetic sets are ``multiplicatively large'' is not known.  It follows from \cref{maintheorem_ae_return_set_is_mult_thick} that almost all ``dynamically syndetic'' sets are multiplicatively thick in a coset of a congruence subsemigroup, but the example in \cref{sec_add_synd_not_mult_thick} shows that there are additively syndetic sets that are not multiplicatively thick in any coset of any congruence subsemigroup.  Nevertheless, the set constructed in that example is multiplicatively syndetic, and hence has positive multiplicative upper density (cf. \cref{def_mult_den_and_synd}) in all cosets of all congruence subsemigroups.

\begin{question}
\label{quest_mult_piecewise_synd}
Do additively syndetic sets have positive multiplicative upper density in a coset of a congruence subsemigroup?
\end{question}

A number of related questions pertaining to the multiplicative largeness of additively syndetic sets were posed in \cite[Sec. 9]{glasscock_koutsogiannis_richter_2019}, and most of them remain open. A positive answer to \cref{quest_mult_piecewise_synd} would yield, by a straightforward application of \Szemeredi{}'s theorem (cf. \cite[Thm. 2.4]{glasscock_koutsogiannis_richter_2019}), arbitrarily long geometric progressions in additively syndetic sets.  Thus, it seems that without a breakthrough, we are very far from a positive answer: it is still an open problem to determine whether or not additively syndetic sets contain a three-term geometric progression or even two distinct integers whose ratio is a perfect square.

The class of $\text{IP}_0^*$ sets and their translates form a subclass of additively syndetic sets closely related to nilsystems in which more is known.  An $\text{IP}_0$ subset of $\Z$ is one that contains arbitrarily large sets of the form
\begin{align*}
\left\{ \sum_{f \in F} x_f \ \middle| \ \emptyset \neq F \subseteq \{1, \ldots, r\} \right\}, \quad x_1, \ldots, x_r \in \Z.
\end{align*}
An $\text{IP}_0^*$ subset of $\Z$ is one which has non-empty intersection with all $\text{IP}_0$ sets.  Such sets were used recently by Bergelson and Leibman \cite{bergelson_leibman_2018} to characterize nilsystems: roughly speaking, a system $(X,T)$ is a nilsystem if and only if all of its return-time sets are translated $\text{IP}_0^*$ sets.

\begin{question}
\label{quest_ipzstar_sets_mult_thick}
Are additive translates of $\text{IP}_0^*$ sets multiplicatively thick in a coset of a congruence subsemigroup?
\end{question}

There is some evidence for a positive answer to \cref{quest_ipzstar_sets_mult_thick}.  It was shown in \cite[Thm. 7.2]{glasscock_koutsogiannis_richter_2019} that $\text{IP}_0^*$ subsets of $\N$ have positive multiplicative upper density in a coset of a subsemigroup of the form $\{n \in \N \ | \ \gcd(n,N)= 1\}$, $N \in \N$.  It was shown in \cite[Thm. II]{bergelson_furstenberg_weiss_2006} that $\Delta_0^*$ sets, a subclass of $\text{IP}_0^*$ sets, are \emph{piecewise Bohr$_0$}, ie. contain the intersection of a Bohr$_0$ set with a set containing arbitrarily long intervals.  More recent work \cite{host_kra_nilbohr_2011} demonstrates that $\text{SG}_d^*$ sets, a combinatorially defined class of sets related to $\text{IP}_0^*$ sets, are \emph{piecewise nil-Bohr$_0$}.  While we cannot immediately deduce anything multiplicatively about piecewise nil-Bohr sets from \cref{mainthm_nilbohr_sets_are_mult_thick}, it seems plausible that progress could be made on \cref{quest_ipzstar_sets_mult_thick} by combining tools and techniques from those works with the ones presented here.

\bibliographystyle{alphanum}
\bibliography{thickbib}

\newcommand{\etalchar}[1]{$^{#1}$}
\begin{thebibliography}{DDM{\etalchar{+}}}

\bibitem[AG1]{akin_glasner_1998}
Ethan Akin and Eli Glasner.
\newblock Topological ergodic decomposition and homogeneous flows.
\newblock In {\em Topological dynamics and applications ({M}inneapolis, {MN},
  1995)}, volume 215 of {\em Contemp. Math.},  43--52. Amer. Math. Soc.,
  Providence, RI, 1998.

\bibitem[AG2]{auslander_guerin_1997}
Joseph Auslander and Marianne Guerin.
\newblock Regional proximality and the prolongation.
\newblock {\em Forum Math.}, 9(6):761--774, 1997.

\bibitem[AGH]{auslander_green_hahn_1963_book}
Louis Auslander, Leon Green, and Frank Hahn.
\newblock {\em Flows on homogeneous spaces}.
\newblock Annals of Mathematics Studies, No. 53. Princeton University Press,
  Princeton, N.J., 1963.
\newblock With the assistance of L. Markus and W. Massey, and an appendix by L.
  Greenberg.

\bibitem[AM]{auslander_markley_1998}
Joseph Auslander and Nelson Markley.
\newblock Dynamical properties of automorphisms of minimal flows.
\newblock {\em Bull. Polish Acad. Sci. Math.}, 46(2):105--112, 1998.

\bibitem[BBHS]{beiglbock_bergelson_hindman_strauss_2006}
Mathias Beiglb\"{o}ck, Vitaly Bergelson, Neil Hindman, and Dona Strauss.
\newblock Multiplicative structures in additively large sets.
\newblock {\em J. Combin. Theory Ser. A}, 113(7):1219--1242, 2006.

\bibitem[Ber1]{berg_1971}
Kenneth Berg.
\newblock Quasi-disjointness in ergodic theory.
\newblock {\em Trans. Amer. Math. Soc.}, 162:71--87, 1971.

\bibitem[Ber2]{bergelson_2000}
Vitaly Bergelson.
\newblock The multifarious {P}oincar\'{e} recurrence theorem.
\newblock In {\em Descriptive set theory and dynamical systems
  ({M}arseille-{L}uminy, 1996)}, volume 277 of {\em London Math. Soc. Lecture
  Note Ser.},  31--57. Cambridge Univ. Press, Cambridge, 2000.

\bibitem[Ber3]{bergelson_2005}
Vitaly Bergelson.
\newblock Multiplicatively large sets and ergodic {R}amsey theory.
\newblock {\em Israel J. Math.}, 148:23--40, 2005.

\bibitem[BFW]{bergelson_furstenberg_weiss_2006}
Vitaly Bergelson, Hillel Furstenberg, and Benjamin Weiss.
\newblock Piecewise-{B}ohr sets of integers and combinatorial number theory.
\newblock In {\em Topics in discrete mathematics}, volume~26 of {\em Algorithms
  Combin.},  13--37. Springer, Berlin, 2006.

\bibitem[BG1]{bergelson_glasscock_2018}
Vitaly Bergelson and Daniel Glasscock.
\newblock Multiplicative richness of additively large sets in {$\mathbb{Z}^d$}.
\newblock {\em J. Algebra}, 503:67--103, 2018.

\bibitem[BG2]{bergelson_glasscock_2020}
Vitaly Bergelson and Daniel Glasscock.
\newblock On the interplay between additive and multiplicative largeness and
  its combinatorial applications.
\newblock {\em J. Combin. Theory Ser. A}, 172:105203, 60, 2020.

\bibitem[BHK]{bergelson_host_kra_2005}
Vitaly Bergelson, Bernard Host, and Bryna Kra.
\newblock Multiple recurrence and nilsequences.
\newblock {\em Invent. Math.}, 160(2):261--303, 2005.
\newblock With an appendix by Imre Ruzsa.

\bibitem[BL1]{bergelson_leibman_1996}
Vitaly Bergelson and Alexander Leibman.
\newblock Polynomial extensions of van der {W}aerden's and {S}zemer\'edi's
  theorems.
\newblock {\em J. Amer. Math. Soc.}, 9(3):725--753, 1996.

\bibitem[BL2]{bergelson_leibman_2007}
Vitaly Bergelson and Alexander Leibman.
\newblock Distribution of values of bounded generalized polynomials.
\newblock {\em Acta Math.}, 198(2):155--230, 2007.

\bibitem[BL3]{bergelson_leibman_2015}
Vitaly Bergelson and Alexander Leibman.
\newblock Cubic averages and large intersections.
\newblock In {\em Recent trends in ergodic theory and dynamical systems},
  volume 631 of {\em Contemp. Math.},  5--19. Amer. Math. Soc., Providence, RI,
  2015.

\bibitem[BL4]{bergelson_leibman_2018}
Vitaly Bergelson and Alexander Leibman.
\newblock {${\rm IP}_r^\ast$}-recurrence and nilsystems.
\newblock {\em Adv. Math.}, 339:642--656, 2018.

\bibitem[Bog]{bogolyubov_1939}
Nikolay~N. Bogolyubov.
\newblock Some algebraical properties of almost periods (in {R}ussian).
\newblock {\em Zap. kafedry mat. fiziki Kiev}, 4:185--194, 1939.

\bibitem[BR]{bergelson_richter_2022}
Vitaly Bergelson and Florian~K. Richter.
\newblock Dynamical generalizations of the prime number theorem and
  disjointness of additive and multiplicative semigroup actions.
\newblock {\em Duke Math. J.}, 171(15):3133--3200, 2022.

\bibitem[CG]{corwin_greenleaf_1990}
Lawrence~J. Corwin and Frederick~P. Greenleaf.
\newblock {\em Representations of nilpotent {L}ie groups and their
  applications. {P}art {I}}, volume~18 of {\em Cambridge Studies in Advanced
  Mathematics}.
\newblock Cambridge University Press, Cambridge, 1990.
\newblock Basic theory and examples.

\bibitem[DDM{\etalchar{+}}]{dong_donoso_maass_shao_ye_2013}
Pandeng Dong, Sebasti\'{a}n Donoso, Alejandro Maass, Song Shao, and Xiangdong
  Ye.
\newblock Infinite-step nilsystems, independence and complexity.
\newblock {\em Ergodic Theory Dynam. Systems}, 33(1):118--143, 2013.

\bibitem[dV]{devries_book_1993}
Jan de~Vries.
\newblock {\em Elements of topological dynamics}, volume 257 of {\em
  Mathematics and its Applications}.
\newblock Kluwer Academic Publishers Group, Dordrecht, 1993.

\bibitem[FH1]{frantzikinakis_host_2017}
Nikos Frantzikinakis and Bernard Host.
\newblock Higher order {F}ourier analysis of multiplicative functions and
  applications.
\newblock {\em J. Amer. Math. Soc.}, 30(1):67--157, 2017.

\bibitem[FH2]{frantzikinakis_host_2018}
Nikos Frantzikinakis and Bernard Host.
\newblock Weighted multiple ergodic averages and correlation sequences.
\newblock {\em Ergodic Theory Dynam. Systems}, 38(1):81--142, 2018.

\bibitem[FM1]{ferremoragues_2021}
Andreu Ferr\'{e}~Moragues.
\newblock Properties of multicorrelation sequences and large returns under some
  ergodicity assumptions.
\newblock {\em Discrete Contin. Dyn. Syst.}, 41(6):2809--2828, 2021.

\bibitem[FM2]{frantzikinakis_mccutcheon_2009}
Nikos Frantzikinakis and Randall McCutcheon.
\newblock Ergodic theory: Recurrence.
\newblock In Robert~A. Meyers, editor, {\em Encyclopedia of Complexity and
  Systems Science},  3083--3095. Springer New York, New York, NY, 2009.

\bibitem[F{\o}l1]{folnerbogoliouboff1954}
Erling F{\o}lner.
\newblock Generalization of a theorem of {B}ogolio\`uboff to topological
  abelian groups. {W}ith an appendix on {B}anach mean values in non-abelian
  groups.
\newblock {\em Math. Scand.}, 2:5--18, 1954.

\bibitem[F{\o}l2]{folnernoteonbogoliouboff1954}
Erling F{\o}lner.
\newblock Note on a generalization of a theorem of {B}ogolio\`uboff.
\newblock {\em Math. Scand.}, 2:224--226, 1954.

\bibitem[For]{fort_1951}
Marion~Kirkland Fort, Jr.
\newblock Points of continuity of semi-continuous functions.
\newblock {\em Publ. Math. Debrecen}, 2:100--102, 1951.

\bibitem[Fra]{frantzikinakis_2015}
Nikos Frantzikinakis.
\newblock Multiple correlation sequences and nilsequences.
\newblock {\em Invent. Math.}, 202(2):875--892, 2015.

\bibitem[Fur1]{furstenberg_1963}
Harry Furstenberg.
\newblock The structure of distal flows.
\newblock {\em Amer. J. Math.}, 85:477--515, 1963.

\bibitem[Fur2]{furstenberg_1977}
Harry Furstenberg.
\newblock Ergodic behavior of diagonal measures and a theorem of {S}zemer\'edi
  on arithmetic progressions.
\newblock {\em J. Analyse Math.}, 31:204--256, 1977.

\bibitem[Fur3]{furstenberg_1981}
Harry Furstenberg.
\newblock Poincar\'{e} recurrence and number theory.
\newblock {\em Bull. Amer. Math. Soc. (N.S.)}, 5(3):211--234, 1981.

\bibitem[Fur4]{furstenberg_book_1981}
Harry Furstenberg.
\newblock {\em Recurrence in ergodic theory and combinatorial number theory}.
\newblock Princeton University Press, Princeton, N.J., 1981.
\newblock M. B. Porter Lectures.

\bibitem[FW]{furstenberg_weiss_1978}
Harry Furstenberg and Benjamin Weiss.
\newblock Topological dynamics and combinatorial number theory.
\newblock {\em J. Analyse Math.}, 34:61--85 (1979), 1978.

\bibitem[GH]{gottschalk_hedlund_1955}
Walter~Helbig Gottschalk and Gustav~Arnold Hedlund.
\newblock {\em Topological dynamics}.
\newblock American Mathematical Society Colloquium Publications, Vol. 36.
  American Mathematical Society, Providence, R.I., 1955.

\bibitem[GHS{\etalchar{+}}]{glasner_huang_shao_weiss_ye_2020}
Eli Glasner, Wen Huang, Song Shao, Benjamin Weiss, and Xiangdong Ye.
\newblock Topological characteristic factors and nilsystems, 2020.

\bibitem[GKR]{glasscock_koutsogiannis_richter_2019}
Daniel Glasscock, Andreas Koutsogiannis, and Florian~Karl Richter.
\newblock Multiplicative combinatorial properties of return time sets in
  minimal dynamical systems.
\newblock {\em Discrete Contin. Dyn. Syst.}, 39(10):5891--5921, 2019.

\bibitem[Gla]{glasner_1994}
Eli Glasner.
\newblock Topological ergodic decompositions and applications to products of
  powers of a minimal transformation.
\newblock {\em J. Anal. Math.}, 64:241--262, 1994.

\bibitem[GT1]{green_tao_2012}
Ben Green and Terence Tao.
\newblock The {M}\"obius function is strongly orthogonal to nilsequences.
\newblock {\em Ann. of Math. (2)}, 175(2):541--566, 2012.

\bibitem[GT2]{green_tao_quantitative_behaviour_2012}
Ben Green and Terence Tao.
\newblock The quantitative behaviour of polynomial orbits on nilmanifolds.
\newblock {\em Ann. of Math. (2)}, 175(2):465--540, 2012.

\bibitem[HJS]{hindman_jones_strauss_2019}
Neil Hindman, Lakeshia~Legette Jones, and Dona Strauss.
\newblock The relationships among many notions of largeness for subsets of a
  semigroup.
\newblock {\em Semigroup Forum}, 99(1):9--31, 2019.

\bibitem[HK1]{host_kra_nilbohr_2011}
Bernard Host and Bryna Kra.
\newblock Nil-{B}ohr sets of integers.
\newblock {\em Ergodic Theory Dynam. Systems}, 31(1):113--142, 2011.

\bibitem[HK2]{host_kra_book_2018}
Bernard Host and Bryna Kra.
\newblock {\em Nilpotent structures in ergodic theory}, volume 236 of {\em
  Mathematical Surveys and Monographs}.
\newblock American Mathematical Society, Providence, RI, 2018.

\bibitem[JST]{jamneshan_shalom_tao_2021}
Asgar Jamneshan, Or~Shalom, and Terence Tao.
\newblock The structure of arbitrary conze-lesigne systems, 2021.

\bibitem[Key]{keynes_1967}
Harvey~B. Keynes.
\newblock A study of the proximal relation in coset transformation groups.
\newblock {\em Trans. Amer. Math. Soc.}, 128:389--402, 1967.

\bibitem[KM]{Kakeya_Morimoto_1930}
S\^{o}ichi Kakeya and Seigo Morimoto.
\newblock On a theorem of {M}{M}. {B}andet and van der {W}aerden.
\newblock {\em Japanese journal of mathematics: transactions and abstracts},
  7:163--165, 1930.

\bibitem[KST]{kolyada_snoha_trofimchuk_2001}
Sergi\u{i} Kolyada, L'ubom\'ir Snoha, and Serge\u{i} Trofimchuk.
\newblock Noninvertible minimal maps.
\newblock {\em Fund. Math.}, 168(2):141--163, 2001.

\bibitem[Kur]{kuratowski_1966}
Kazimierz Kuratowski.
\newblock {\em Topology. {V}ol. {I}}.
\newblock Academic Press, New York-London; Pa\'{n}stwowe Wydawnictwo Naukowe
  [Polish Scientific Publishers], Warsaw, 1966.
\newblock New edition, revised and augmented, Translated from the French by J.
  Jaworowski.

\bibitem[Lei1]{leibman_rational_2006}
Alexander Leibman.
\newblock Rational sub-nilmanifolds of a compact nilmanifold.
\newblock {\em Ergodic Theory Dynam. Systems}, 26(3):787--798, 2006.

\bibitem[Lei2]{leibman_2010}
Alexander Leibman.
\newblock Multiple polynomial correlation sequences and nilsequences.
\newblock {\em Ergodic Theory Dynam. Systems}, 30(3):841--854, 2010.

\bibitem[Lei3]{Leibman_2015}
Alexander Leibman.
\newblock Nilsequences, null-sequences, and multiple correlation sequences.
\newblock {\em Ergodic Theory Dynam. Systems}, 35(1):176--191, 2015.

\bibitem[Tao]{tao_hofa_book}
Terence Tao.
\newblock {\em Higher order {F}ourier analysis}, volume 142 of {\em Graduate
  Studies in Mathematics}.
\newblock American Mathematical Society, Providence, RI, 2012.

\bibitem[vdW]{vanderWaerden_1927}
Bartel~Leendert van~der Waerden.
\newblock Beweis einer baudetschen vermutung.
\newblock {\em Nieuw Arch. Wiskd., II. Ser.}, 15:212--216, 1927.

\bibitem[Ye]{ye_1992}
Xiangdong Ye.
\newblock D-function of a minimal set and an extension of {S}harkovskii's
  theorem to minimal sets.
\newblock {\em Ergodic Theory and Dynamical Systems}, 12(2):365--376, 1992.

\end{thebibliography}

\bigskip
\footnotesize
\noindent
Daniel Glasscock\\
\textsc{University of Massachusetts Lowell}\par\nopagebreak
\noindent
\href{mailto:daniel_glasscock@uml.edu}
{\texttt{daniel{\_}glasscock@uml.edu}}

\end{document}